\documentclass[reqno,12pt]{amsart}%
\usepackage{amsfonts}
\usepackage{amssymb}
\usepackage{amsmath}
\usepackage{lscape}
\usepackage{graphicx}%
\providecommand{\U}[1]{\protect\rule{.1in}{.1in}}
\renewcommand{\theequation}{\thesection.\arabic{equation}}
\renewcommand{\baselinestretch}{1}

\newfont{\bbf}{cmbx12 scaled 1435}

\newtheorem{thm}{Theorem}
\newtheorem{prop}{Proposition}
\newtheorem{lem}{Lemma}

\newtheorem{hp}{Assumption}

\renewcommand{\theequation}{\thesection.\arabic{equation}}

\renewcommand{\baselinestretch}{1.7}
\ifx\pdfoutput\relax\let\pdfoutput=\undefined\fi
\newcount\msipdfoutput
\ifx\pdfoutput\undefined\else
\ifcase\pdfoutput\else
\ifx\paperwidth\undefined\else
\ifdim\paperheight=0pt\relax\else\pdfpageheight\paperheight\fi
\ifdim\paperwidth=0pt\relax\else\pdfpagewidth\paperwidth\fi
\fi\fi\fi
\begin{document}
\parindent=0.1in

\thispagestyle{empty}

\vspace*{1cm}

\begin{center}
Robust Adaptive Rate-Optimal Testing for the White Noise Hypothesis
\footnote{This is a revised version of a paper previously entitled `Adaptive
Rate-Optimal Detection of Small Correlation Coefficients'. We would like to
thank Richard Davis, Marcelo Fernandes, Liudas Giraitis, Bruce Hansen, George
Kapetanios, Remigijus Leipus, Peter Phillips and Aris Spanos for stimulating
questions and suggestions. We would also like to thank participants of the
Queen Mary Econometric Reading Group, Economic Seminar at York, Southampton
and Warwick Universities, 2008 Oxbridge Time Series Workshop, 2008 Vienna
Model Selection Workshop, 2008 North American Econometric Society Conference
and 2008 European Econometric Society Conference, Journ\'{e}es 2008 de
Statistiques de Rennes, 2009 Bristol Econometrics Worshop and 2011 Brazilian
Time Series and Econometrics School. Last but not least, we would like to
thank an anonymous associate editor and two referees whose feedback helped to
improve the paper. All remaining errors are ours. The first author gratefully
acknowledges financial support of Fonds de recherche sur la soci\'{e}t\'{e} et
la culture (FQRSC) and the Social Sciences and Humanities Research Council of
Canada (SSHRC). The last two authors are thankful for the financial support
from the School of Economics and Finance of Queen Mary, University of London.}

(with supplementary material)

\vspace*{0.5 cm}

Alain Guay\footnote{CIRP\'{E}E and CIREQ, Universit\'{e} du Qu\'{e}bec \`{a}
Montr\'{e}al, e-mail: \texttt{guay.alain@uqam.ca}} 

\medskip

Emmanuel Guerre\footnote{School of Economics and Finance, Queen Mary,
University of London, e-mail: \texttt{e.guerre@qmul.ac.uk }}\\[0pt]

\medskip

\v{S}t\v{e}p\'{a}na Lazarov\'{a}\footnote{School of Economics and Finance,
Queen Mary, University of London, e-mail: \texttt{s.lazarova@qmul.ac.uk }%
}\\[0pt]\vspace{0.5cm}

This version: 2nd November 2012
\end{center}

\vfill

\begin{center}
\pagebreak
\end{center}

\thispagestyle{empty}

\begin{center}
\textbf{Abstract}

\begin{minipage}[t]{12cm}
\begin{small}
A new  test is proposed for the weak white noise null hypothesis.  The test is based on a new automatic selection of the order for a Box-Pierce (1970) test statistic or the test statistic of Hong  (1996).
The heteroskedasticity and autocorrelation-consistent (HAC) critical values from Lee (2007) are used, allowing for estimation of the error term. The data-driven order selection is tailored to
detect a new class of alternatives with autocorrelation coefficients which can be
$o(n^{-1/2})$ provided there is sufficiently many of such coefficients. A simulation experiment illustrates the good statistical properties of the test both under
the weak white noise null and the alternative.
\bigskip
\linebreak
\bigskip
\linebreak
\end{small}
\end{minipage}

\end{center}

\renewcommand{\thefootnote}{\arabic{footnote}} \setcounter{footnote}{0}
\setcounter{page}{0} \newpage

\section{Introduction}

Testing for white noise is important in many econometric contexts. Ignoring
autocorrelation of the error terms in a linear regression model can lead to
erroneous confidence intervals and tests. Correlation of residuals from an
ARMA model or of the squared residuals from an ARCH model can indicate an
improper choice of the order. Investigation of autocorrelation function is
also a popular diagnostic tool in macroeconomics and finance, see e.g. Durlauf
(1991) and Campbell, Lo and Craig MacKinlay (1997). Earliest tests of the
white noise hypothesis were based on confidence intervals for autocorrelation
coefficients as described by Fan and Yao (2005). See also Xiao and Wu (2011)
who have recently derived the asymptotic distribution of the maximum
standardized sample covariance of weak white noise, that is an stationary
process which is uncorrelated but possibly dependent. A second approach was
established by Grenander and Rosenblatt (1952) who extended goodness-of-fit
tests such as Kolmogorov and Cram\'{e}r-von Mises tests to tests of white
noise hypothesis. Grenander and Rosenblatt (1952) has been refined by Durlauf
(1991), Anderson (1993) and Deo (2000). Delgado, Hidalgo and Velasco (2005)
have studied a modified test statistic to be used with residuals. Shao (2011a)
has recently extended this setup to cover the weak white noise null
hypothesis. A third approach, pioneered by Box and Pierce (1970), is based on
the sum of squared sample autocorrelation coefficients up to a given order
$p$. Delgado and Velasco (2012), Francq, Roy and Zakoian (2005), Kuan and Lee
(2006) and Lobato (2001) have considered the weak white noise hypothesis. The
case where $p$ grows with the sample size $n$ has been considered by Hong
(1996) in a strong white noise setup and recently extended to the weak white
noise null hypothesis by Shao (2011b) and Xiao and Wu (2011).

This paper contributes to the literature by proposing a data-driven choice
$\hat{p}$ of the order $p$ used in a Box-Pierce type statistic for a test of
the weak white noise null hypothesis. Under this null, $\widehat{p}$ tends to
$1$ in probability so that the null limit behavior of the test statistic is
driven by the first-order sample autocovariance. It is shown that the test can
be implemented using robust critical values of Lee (2007) who extends the work
of Lobato (2001) for the case of observed variables and of Kuan and Lee (2006)
for the case of residuals. The general framework of Lee (2007) includes as a
specific case standardization using steep origin kernels proposed by Phillips,
Sun and Jin (2006) which can improve the power of the resulting test. Under
the alternative, the data-driven $\widehat{p}$ can be as large as necessary.

An appealing feature of Cram\'{e}r-von Mises type of tests is the ability to
detect Pitman local directional alternatives converging to the null with the
parametric rate $n^{-1/2}$. This contrasts with detection results for
Box-Pierce\ type test by Hong (1996) which is only consistent under slower
rates of convergence for local alternatives defined through the spectral
density function. The conclusions of Hong (1996) suggest that Cram\'{e}r-von
Mises tests are more powerful than Box-Pierce tests. One of the contributions
of the present paper is to point out that this ranking of two types of tests
is not universal and there exist classes of alternatives against which
Box-Pierce tests are more powerful than Cram\'{e}r-von Mises tests.

We illustrate this point using a new class of alternatives defined through the
autocovariance function. The new class of alternatives formalizes the idea
that small autocorrelation coefficients of magnitude $\rho_{n}$ can be
detected provided that there are sufficiently many coefficients present at
smaller lags. An important finding of the paper is that detection is still
possible for very small $\rho_{n}$, namely for $\rho_{n}=o\left(
n^{-1/2}\right)  $. As described in Section \ref{Optimal}, this type of
alternatives includes moving average processes with a significant long term
multiplier but $o\left(  n^{-1/2}\right)  $ impulse response coefficients.
Such processes therefore correspond to a macroeconomic scenario where short
term policies have no significant effects whereas long term policies may have
an impact. For such alternatives, the conditional expectation of the present
given the past gives $o\left(  n^{-1/2}\right)  $ weights to each lagged
observations. Therefore this process is hard to predict since it is very close
to a martingale difference process. These alternatives can be of interest in
finance where arbitrage could forbid strong deviations from martingale difference.

Why such alternatives can be detected by Box-Pierce tests can be intuitively
explained as follows. Let $\widehat{R}_{j}$ and $R_{j}$ be respectively the
sample and population covariance at lag $j$. Following Hong (1996), Shao
(2011b) and Xiao and Wu (2011), the nonrobust critical region of the
Box-Pierce test of order $p_{n}\rightarrow\infty$ is%
\begin{equation}
\frac{n\sum_{j=1}^{p_{n}}\left(  \widehat{R}_{j}^{2}/\widehat{R}_{0}%
^{2}-1\right)  }{\left(  2p_{n}\right)  ^{1/2}}\geq c_{\alpha}\text{,}%
\label{BPcrit}%
\end{equation}
where $c_{\alpha}$ is a standard normal critical value. Arguing as Shao
(2011b, Theorem 2.2) suggests that
\begin{equation}
\frac{n\sum_{j=1}^{p_{n}}\left(  \widehat{R}_{j}^{2}/\widehat{R}_{0}%
^{2}-1\right)  }{\left(  2p_{n}\right)  ^{1/2}}=\frac{n\sum_{j=1}^{p_{n}}%
R_{j}^{2}/R_{0}^{2}}{\left(  2p_{n}\right)  ^{1/2}}+O_{\mathbb{P}}\left(
1\right)  .\label{BPasymp}%
\end{equation}
\textbf{ }(\ref{BPasymp}) suggests that the Box-Pierce test is consistent
provided $\left(  n/\left(  2p_{n}\right)  ^{1/2}\right)  \sum_{j=1}^{p_{n}%
}R_{j}^{2}/R_{0}^{2}$ is large enough. Let $N_{n}$ be the number of
correlation coefficients $R_{j}^{2}/R_{0}^{2}\geq\rho_{n}^{2}$ for
$j\in\left[  1,p_{n}\right]  $, so that $\left(  n/\left(  2p_{n}\right)
^{1/2}\right)  \sum_{j=1}^{p_{n}}R_{j}^{2}/R_{0}^{2}\geq nN_{n}\rho_{n}%
^{2}/\left(  2p_{n}\right)  ^{1/2}$. The Box-Pierce test is consistent if
\begin{equation}
n^{1/2}\left(  \frac{N_{n}}{p_{n}^{1/2}}\right)  ^{1/2}\rho_{n}\rightarrow
\infty,\label{BPcons}%
\end{equation}
a condition which allows for $\rho_{n}=o\left(  n^{-1/2}\right)  $ provided
there are enough correlation coefficients larger than $\rho_{n}$, that is,
$N_{n}/p_{n}^{1/2}\rightarrow\infty$, which holds in particular when the exact
order of $N_{n}$ is $p_{n}$. In other words, summing squared sample
correlations in the Box-Pierce statistic allows us to detect very small
population correlations provided they are not too sparse and are concentrated
at lags smaller than $p_{n}$. As shown in this paper, such alternatives are
not detected by Cram\'{e}r-von Mises tests.

An important limitation of the critical region (\ref{BPcrit}) is the use of an
ad hoc order $p_{n}$. Many authors consider a deterministic $p_{n}$ such that
$p_{n}\rightarrow\infty$. This choice of order is inadequate for detecting
alternatives with correlations at low lags: taking $p_{n}=30$ for instance is
unlikely to give a test with power against popular $AR(1)$ or $MA(1)$
alternatives on samples of moderate size. Conversely, taking a fixed $p_{n}$
is not suitable for detecting higher order alternatives. The need to properly
address the tuning of a smoothing parameter with a role similar to $p_{n}$ has
spurred the development of data-driven approaches for various nonparametric
testing problems. The so-called adaptive approach, focuses on data-driven
tests which detect alternatives in a smoothness class converging to the null
at the fastest possible rate given that the smoothness class is unknown to the
test user. See in particular Fan (1996), Spokoiny (1996), Horowitz and
Spokoiny (2001), Guerre and Lavergne (2005), Guay and Guerre (2006) and Chen
and Gao (2007) for various nonparametric models and related null hypotheses of
theoretical or practical relevance. Golubev, Nussbaum and Zhou (2010)\ has
proved Le Cam equivalence of Gaussian time series with spectral densities in a
Besov space with the continuous-time Gaussian white noise model considered in
Spokoiny (1996). This result is limited to Gaussian time series and is not
useful in practice since it does not deliver ready-to-apply white noise tests.
In fact, most of the data-driven choices of $p_{n}$ proposed in the white
noise testing literature are not adaptive rate-optimal. As an exception, Fan
and Yao (2005) extend the work of Fan (1996), outlining but not analyzing a
data-driven test which is based on the maximum of a set of standardized
Box-Pierce statistics, see also Golubev et al. (2010).

A popular data-driven method of choosing the order is the selection procedure
proposed by Newey and West (1994) in the context of long run variance
estimation. See, among other, the simulation section of Hong and Lee (2005).
This selection procedure is however difficult to justify theoretically. Newey
and West selection method, although being optimal for long-run variance
estimation, does not produce a rate-optimal test because the optimal order for
testing differs from the optimal order for estimation, see e.g. Guerre and
Lavergne (2002) and the references therein. Escanciano and Lobato (2009) study
a data-driven choice of order based on an AIC/BIC criterion which is suitable
for estimation but is not adaptive rate-optimal for tests of the white noise
hypothesis. This contrasts with the new data-driven tests proposed here.

The paper is organized as follows. Section \ref{Construction of the test}
describes the penalty approach leading to the data-driven order $\widehat{p}$
and the construction of the rejection region of the test. Section
\ref{Main results} studies the statistical properties of the test under the
general weak white noise null hypothesis and under the new class of
alternatives mentioned above. It illustrates the importance of the choice of a
suitable penalty both under the null and the alternative. Section
\ref{Optimal} states our adaptive rate-optimality results and compares the new
test with the Cram\'{e}r-von Mises test of Deo (2000), the data-driven test of
Escanciano and Lobato (2009) and the maximum test of Xiao and Wu (2011).
Section \ref{Simulation experiments} reports a simulation experiment that
proposes a calibration of the penalty term and compares our automatic test
with other data-driven tests, including tests of Deo (2000) or Escanciano and
Lobato (2009) and a test that uses the Newey and West (1994) plug-in order
selection procedure. Section \ref{Concluding remarks} concludes. Proofs can be
found in the supplementary material.

\section{Construction of the test and choice of the critical
values\label{Construction of the test}}

\setcounter{equation}{0}

Consider a variable $u_{t}$, $t=1,...,n$, which is either directly observed or
defined as the error of a parametric model $m(X_{t};\theta)=u_{t}$ with some
observed covariate $X_{t}$. In the later case $u_{t}$ is not observed but can
be estimated using the residuals $\widehat{u}_{t}=u_{t}(\widehat{\theta})$
where $\widehat{\theta}$\ is an estimator of $\theta$. We are interested in
testing that $u_{t}$\ is uncorrelated. Suppose $\{u_{t}\}$\ is a stationary
process with zero mean and covariance function $R_{j}=\mathrm{Cov}%
(u_{t},u_{t+j})$. The null and alternative hypotheses are then
\[
\mathcal{H}_{0}:R_{j}=0\text{ for all }j\neq0\quad\quad\text{versus}\quad
\quad\mathcal{H}_{1}:R_{j}\neq0\text{ for some }j\neq0.
\]
A natural estimator of the covariance is $\widehat{R}_{j}=\sum_{t=1}%
^{n-|j|}\widehat{u}_{t}\widehat{u}_{t+|j|}/n$, $j=0,\pm1,\ldots,\pm(n-1)$,
which uses the residuals as if they were the true error terms. Given the
kernel spectral density estimator
\[
\hat{f}_{n}(\lambda;p)=\frac{1}{2\pi}\sum_{j=-\infty}^{\infty}K\left(
\frac{\left\vert j\right\vert }{p}\right)  \widehat{R}_{j}\exp\left(
-ij\lambda\right)  \text{,\quad\quad}K\left(  0\right)  =1\quad\quad\text{and
}\int_{0}^{\infty}K\left(  x\right)  dx=1,
\]
where the support of $K$ is $\left[  0,1\right]  $, Hong (1996) has proposed a
test statistic
\begin{equation}
\widehat{S}_{p}=n\pi\int_{-\pi}^{\pi}\left\vert \hat{f}_{n}(\lambda
;p)-\frac{\widehat{R}_{0}}{2\pi}\right\vert ^{2}d\lambda=n\sum_{j=1}%
^{n-1}K^{2}\left(  \frac{j}{p}\right)  \widehat{R}_{j}^{2}.\label{Sp hat}%
\end{equation}
For the uniform kernel $K(t)=\mathbb{I}(t\in\lbrack0,1])$ and up to a division
by $\widehat{R}_{0}^{2}$,$\ \widehat{S}_{p}$ is the Box-Pierce statistic
$\widehat{BP}_{p}/\widehat{R}_{0}^{2}=n\sum_{j=1}^{p}\widehat{R}_{j}%
^{2}/\widehat{R}_{0}^{2}$. Large values of $\widehat{S}_{p}$\ indicate
evidence against the null. Under certain weak dependence conditions on the
weak white noise $\left\{  u_{t}\right\}  $ and for $p=p_{n}\rightarrow\infty$
growing with a suitable rate, Shao (2011b) shows that $\left(  \left(
\widehat{S}_{p}-\widehat{S}_{1}\right)  /R_{0}^{2}-E_{\Delta}(p)\right)
/V_{\Delta}(p)$ converges to a standard normal where%
\begin{align*}
E_{\Delta}(p) &  =\sum_{j=1}^{n-1}\left(  1-\frac{j}{n}\right)  \left(
K^{2}\left(  \frac{j}{p}\right)  -K^{2}\left(  j\right)  \right)  ,\\
V_{\Delta}^{2}(p) &  =2\sum_{j=1}^{n-1}\left(  1-\frac{j}{n}\right)
^{2}\left(  K^{2}\left(  \frac{j}{p}\right)  -K^{2}\left(  j\right)  \right)
^{2},
\end{align*}
and we shall use accordingly $E_{\Delta}(p)$ and $V_{\Delta}^{2}(p)$ as a
standardization for\textbf{ }$\left(  \widehat{S}_{p}-\widehat{S}_{1}\right)
/R_{0}^{2}$. In this notation, the subscript \textquotedblleft$\Delta
$\textquotedblright\ indicates difference $\widehat{S}_{p}-\widehat{S}_{1}$.
For the Box-Pierce statistic, $E_{\Delta}(p)=\left(  p-1\right)  \left(
1+O\left(  p/n\right)  \right)  $ and $V_{\Delta}^{2}(p)=2\left(  p-1\right)
\left(  1+O\left(  p/n\right)  \right)  $ and these approximations remain
valid for other kernels up to a multiplicative constant. We propose to select
$\widehat{p}$\ as the smallest integer number maximizing the penalized
statistic,\textbf{\ }%
\begin{align}
\widehat{p} &  =\arg\max_{p\in\left[  1,\overline{p}_{n}\right]  }\left(
\frac{\widehat{S}_{p}}{\widehat{R}_{0}^{2}}-E\left(  p\right)  -\gamma
_{n}V_{\Delta}(p)\right)  \nonumber\\
&  =\arg\max_{p\in\left[  1,\overline{p}_{n}\right]  }\left(  \frac
{\widehat{S}_{p}-\widehat{S}_{1}}{\widehat{R}_{0}^{2}}-E_{\Delta}%
(p)-\gamma_{n}V_{\Delta}(p)\right)  ,\label{Hatp}%
\end{align}
where $E(p)=\sum_{j=1}^{n-1}\left(  1-j/n\right)  K^{2}\left(  j/p\right)  $
and $\overline{p}_{n}\leq n-1$. This penalization procedure is similar to
penalization proposed by Guay and Guerre (2006) or Guerre and Lavergne (2005).
It differs from the penalization used in the AIC or BIC procedures which use a
higher penalty term $\gamma_{n}E\left(  p\right)  $ in place of $E\left(
p\right)  +\gamma_{n}V_{\Delta}(p)$. Escanciano and Lobato (2009) similarly
use penalty term $\widehat{\gamma}_{n}E\left(  p\right)  $ for $p$ in a
bounded finite set.

The intuition for $\widehat{p}$ is as follows. Note first that (\ref{Hatp})
uses the difference $\widehat{S}_{p}-\widehat{S}_{1}$. The idea here is that
the test should be based on $\widehat{S}_{1}$ unless $\widehat{S}_{p}%
-\widehat{S}_{1}$ is large enough for some $p$. Since the criterion maximized
in (\ref{Hatp}) is equal to $0$ for $p=1$, $\widehat{p}$ differs from $1$
whenever there is a $p$ such that $\left(  \widehat{S}_{p}-\widehat{S}%
_{1}\right)  /\widehat{R}_{0}^{2}-E_{\Delta}(p)-\gamma_{n}V_{\Delta}(p)>0$ or
equivalently%
\begin{equation}
\frac{\left(  \widehat{S}_{p}-\widehat{S}_{1}\right)  /\widehat{R}_{0}%
^{2}-E_{\Delta}(p)}{V_{\Delta}(p)}>\gamma_{n},\label{Rej}%
\end{equation}
an inequality which, in view of the asymptotic normality established by Shao
(2011b) under the null, has the flavour of a one-sided significance test using
a critical value $\gamma_{n}$. Such a construction suggests that the
data-driven statistic $\widehat{S}_{\widehat{p}}$ better captures higher order
covariances than $\widehat{S}_{1}$. Therefore, rejecting the null when
$\widehat{S}_{\widehat{p}}\geq z$ should give a more powerful test than the
test $\widehat{S}_{1}\geq z$ based on $\widehat{S}_{1}$ and the same critical
value $z$ as recommended below. See (\ref{Improv}) in Theorem \ref{AltCV} for
a more formal statement. Why the chosen $\widehat{p}$ should have certain
optimality properties can be seen by viewing (\ref{Hatp}) as a bias variance
trade-off. Theorem 2.2 in Shao (2011b) suggests that $\left(  \widehat{S}%
_{p}-\widehat{S}_{1}\right)  /\widehat{R}_{0}^{2}-E_{\Delta}(p)$ is an
estimator of $n\sum_{j=2}^{\infty}R_{j}^{2}$ with a bias $n\sum_{j=p+1}%
^{\infty}R_{j}^{2}$ and a standard deviation $V_{\Delta}(p)$. Hence
(\ref{Hatp}) choses a $p$ which maximizes $-n\sum_{j=p+1}^{\infty}R_{j}%
^{2}-\gamma_{n}V_{\Delta}(p)$ and therefore achieves the so called bias
variance trade-off, leading to a data-driven test statistic $\widehat
{S}_{\widehat{p}}=\widehat{S}_{1}+\widehat{S}_{\widehat{p}}-\widehat{S}_{1}$
with the best potential to detect an alternative.

Under $\mathcal{H}_{0}$, it is expected that $\widehat{p}=1$ with a high
probability provided $\gamma_{n}$ is large enough since all the $\widehat
{S}_{p}-\widehat{S}_{1}$ estimate 0. Since $\widehat{S}_{\widehat{p}}%
=\widehat{S}_{1}+o_{\mathbb{P}}\left(  1\right)  $ under the null, the
critical values of the test can be taken to be the same as the critical values
of the test based upon the simple statistic $\widehat{S}_{1}$. A HAC-robust
standardization of $\widehat{S}_{1}$ is given in Lee (2007). In the case where
$u_{t}$ is observed, an inconsistent \textquotedblleft
estimator\textquotedblright\ of the long run variance of $\sum_{t=1}%
^{n-1}u_{t}u_{t+1}/(n-1)$ is, for a kernel $k\left(  \cdot\right)  $,
$k_{ij}=k\left(  \left\vert i-j\right\vert /n\right)  $ and $\varphi_{i}%
=\sum_{t=1}^{i-1}\left(  u_{t}u_{t+1}-\widehat{R}_{1}\right)  /n^{1/2}$,%
\[
\widetilde{\Gamma}_{1}=\sum_{i=1}^{n-1}\sum_{j=1}^{n-1}\left(  \left(
k_{ij}-k_{i,j+1}\right)  -\left(  k_{i+1,j}-k_{i+1,j+1}\right)  \right)
\varphi_{i}\varphi_{j}.
\]
For residuals $\hat{u}_{t}$, let $\widehat{\theta}_{i}$ be the estimator
$\widehat{\theta}$\ computed with the first $i$\ observations and estimate
$\varphi_{i}$ recursively by $\widehat{\varphi}_{i}=\sum_{t=1}^{i-1}\left(
u_{t}\left(  \widehat{\theta}_{i}\right)  u_{t+1}\left(  \widehat{\theta}%
_{i}\right)  -\widehat{R}_{1}\right)  /n^{1/2}$. Let%
\[
\widehat{\Gamma}_{1}=\sum_{i=1}^{n-1}\sum_{j=1}^{n-1}\left(  \left(
k_{ij}-k_{i,j+1}\right)  -\left(  k_{i+1,j}-k_{i+1,j+1}\right)  \right)
\widehat{\varphi}_{i}\widehat{\varphi}_{j}.
\]
It follows from Lee (2007) that the limit distribution of $n\widehat{R}%
_{1}/\widetilde{\Gamma}_{1}$ when $u_{t}$ is observed and of $n\widehat{R}%
_{1}/\widehat{\Gamma}_{1}$ when $u_{t}$ is is estimated by residuals $\hat
{u}_{t}$ is, assuming that $k\left(  \cdot\right)  $ is twice continuously
differentiable%
\begin{equation}
\frac{W^{2}\left(  1\right)  }{-\int_{0}^{1}\int_{0}^{1}k^{\prime\prime
}\left(  r-s\right)  \left(  W\left(  r\right)  -rW\left(  1\right)  \right)
\left(  W\left(  s\right)  -sW\left(  1\right)  \right)  drds}\label{StdW}%
\end{equation}
where $W$\ is a standard Brownian motion. Let $z_{L}\left(  \alpha\right)  $
be be the $\left(  1-\alpha\right)  $th quantile of (\ref{StdW}). The critical
values and rejection region of the test are
\begin{align}
\widehat{z}_{L}(\alpha) &  =K^{2}\left(  1\right)  \widetilde{\Gamma}_{1}%
z_{L}\left(  \alpha\right)  ,\label{Zlob}\\
\widehat{z}_{KL}(\alpha) &  =K^{2}\left(  1\right)  \widehat{\Gamma}_{1}%
z_{L}\left(  \alpha\right)  ,\label{Zest}%
\end{align}%
\begin{equation}
\widehat{S}_{\widehat{p}}\geq\widehat{z}(\alpha)\text{\quad where }\widehat
{z}(\alpha)=\left\{
\begin{array}
[c]{ll}%
\widehat{z}_{L}(\alpha) & \text{for observed }\left\{  u_{t}\right\}
\text{,}\\
\widehat{z}_{KL}(\alpha) & \text{for residuals }\left\{  \widehat{u}%
_{t}\right\}  \text{.}%
\end{array}
\right.  \label{Test}%
\end{equation}
We also consider a modified version of the test which employs a
standardization of the sample covariances as used by Deo (2000) or Escanciano
and Lobato (2009),%
\begin{equation}
\widehat{S}_{p}^{\ast}=n\sum_{j=1}^{n-1}K^{2}\left(  \frac{j}{p}\right)
\left(  \frac{\widehat{R}_{j}}{\widehat{\tau}_{j}}\right)  ^{2}\text{\quad
where }\widehat{\tau}_{j}^{2}=\frac{1}{n-j}\sum_{t=1}^{n-j}\widehat{u}_{t}%
^{2}\widehat{u}_{t+j}^{2}-\left(  \frac{n}{n-j}\widehat{R}_{j}\right)
^{2}.\label{Sstar}%
\end{equation}
The sample variance $\widehat{\tau}_{j}^{2}$ is an estimator of $\tau_{j}%
^{2}=\operatorname*{Var}\left(  u_{t}u_{t+j}\right)  $ which, for observed
$u_{t}$, is the asymptotic variance of $n^{1/2}\left(  \widehat{R}_{j}%
-R_{j}\right)  $ in the case of uncorrelated $u_{t}u_{t+j}$ or for martingale
difference. The corresponding data-driven order $p$ and critical values are%
\begin{align}
\widehat{p}^{\ast} &  =\arg\max_{p\in\left[  1,\overline{p}_{n}\right]
}\left(  \widehat{S}_{p}^{\ast}-E\left(  p\right)  -\gamma_{n}V_{\Delta
}(p)\right)  ,\label{Hatpstar}\\
\widehat{z}^{\ast}(\alpha) &  =\frac{\widehat{z}(\alpha)}{\widehat{\tau}%
_{1}^{2}}.\label{Zstar}%
\end{align}
While the test (\ref{Test}) is studied in Theorems \ref{Level} and
\ref{Sparse}, the test with rejection region $\widehat{S}_{\widehat{p}^{\ast}%
}^{\ast}\geq\widehat{z}^{\ast}(\alpha)$ is studied in Theorem \ref{Extension}.

Let us now turn to notations and our main assumptions. In what follows,
$a_{n}\asymp b_{n}$ means that the sequences $\left\{  a_{n}\right\}  $ and
$\left\{  b_{n}\right\}  $ have the same order, i.e. that $a_{n}/b_{n}$ and
$b_{n}/a_{n}$ are both $O\left(  1\right)  $. For a real random variable $Z$
and a positive real number $a$, $\left\Vert Z\right\Vert _{a}=\mathbb{E}%
^{1/a}\left[  \left\vert Z\right\vert ^{a}\right]  $. Consider first the case
of observed $u_{t}$. When studying the performance of the test under the
alternative, we consider a sequence $\left\{  u_{t,n}\right\}  $ of stationary
alternatives with autocovariance coefficients $\left\{  R_{j,n}\right\}  $.
This means that for each given $n$, the process $\left\{  u_{t,n}%
,t\in\mathbb{N}\right\}  $ is stationary. This type of sequences includes for
instance local $MA\left(  \infty\right)  $ alternatives $u_{t,n}%
=\varepsilon_{t}+\sum_{i=1}^{\infty}a_{i,n}\varepsilon_{t-i}$ where
$a_{i,n}\rightarrow0$ when $n$ grows. Further, for residuals $\widehat{u}%
_{t}=u_{t}\left(  \widehat{\theta}\right)  $, we assume that $\sqrt{n}\left(
\widehat{\theta}-\theta_{n}\right)  $ is asymptotically centered with
$\theta_{n}$ is a pseudo-true value and set $u_{t}\left(  \theta_{n}\right)
=u_{t,n}$. For the sake of brevity, $\left\{  u_{t,n}\right\}  $ and $\left\{
R_{j,n}\right\}  $ are abbreviated to $\left\{  u_{t}\right\}  $ and $\left\{
R_{j}\right\}  $ in the rest of the paper but we maintain the dependence with
respect to $n$ when stating our main assumptions. Under the null and the
alternative, we follow Shao (2011b), Xiao and Wu (2011), and restrict
ourselves to stationary processes satisfying a moment contraction condition by
Wu (2005). We assume that $u_{t,n}=F_{n}\left(  \ldots,e_{t-1},e_{t}\right)  $
for some measurable $F$, where $e_{t}$, $t=-\infty,\ldots,+\infty$, are i.i.d.
(univariate or vector) random variables. Consider an independent copy
$\left\{  e_{t}^{\prime}\right\}  $ of $\left\{  e_{t}\right\}  $ and define%
\[
u_{t,n}^{\tau}=F_{n}\left(  \ldots,e_{\tau-1},e_{\tau}^{\prime},e_{\tau
+1},\ldots,e_{t-1},e_{t}\right)  \quad\quad\quad\tau\leq t\leq n,
\]
where $e_{\tau}$ is changed to $e_{\tau}^{\prime}$. Assume that for some $a>0$
and for all $j\geq0$,%
\[
\left\Vert u_{t,n}-u_{t,n}^{t-j}\right\Vert _{a}\leq\delta_{a}\left(
j\right)  \text{\quad\quad where }\delta_{a}\left(  j\right)  \rightarrow
0\text{ when }j\rightarrow\infty,
\]
a condition meaning that shocks cannot have a long run impact. A fast decrease
of $\delta_{a}\left(  j\right)  $ also ensures that $u_{t}=u_{t,n}$ becomes
independent of $u_{t-j}$ when $j$ grows as the $\alpha$-mixing assumption used
in Francq et al. (2005) or Delgado and Velasco (2012). Shao (2011b) assumes
that $\delta_{a}\left(  j\right)  $ decreases at an exponential rate, a
condition which is satisfied by many linear and nonlinear time series models,
including threshold, stochastic volatility, bilinear or GARCH models, see Shao
(2011b), Wu (2005, 2007) and the references therein. Our main assumptions are
given below.\setcounter{hp}{10}

\begin{hp}
[Kernel]The kernel function $K\left(  \cdot\right)  $ in (\ref{Sp hat}) from
$\mathbb{R}^{+}$ to $\left[  0,\infty\right)  $ is nonincreasing, bounded away
from $0$ on $[0,1/2]$ and continuous differentiable over its support $[0,1]$.
The kernel $k\left(  \cdot\right)  $ used for the critical values is twice
continuously differentiable over its compact support.\label{Kernel}
\end{hp}

\setcounter{hp}{17}

\begin{hp}
[Regularity]Under $\mathcal{H}_{0}$ and $\mathcal{H}_{1}$, $\sup_{t}\left\Vert
u_{t,n}\right\Vert _{12a}<C_{0}R_{0,n}^{1/2}$ for some $a>1$ and, for some
$b>0$, $\delta_{12a}\left(  j\right)  \leq C_{1}j^{-7-b}$. Moreover
$1/C_{2}\leq R_{0,n}\leq C_{2}$, and \newline$\max_{j\in\left[  1,\overline
{p}_{n}\right]  }R_{0,n}^{2}/\operatorname*{Var}\left(  u_{t,n}u_{t+j,n}%
\right)  \leq C_{3}$.\label{Reg}
\end{hp}

\setcounter{hp}{15}

\begin{hp}
[Order $p$]The maximal order $\overline{p}_{n}$ diverges faster than some
power of $n$ with $\overline{p}_{n}=o(n^{1/\left(  2\left(  1+3/a\right)
\right)  })$ as $n\rightarrow\infty$, where $a>1$ is the same constant as in
Assumption \ref{Reg} above. The penalty sequence $\gamma_{n}$ satisfies
$\gamma_{n}>0$, $\gamma_{n}\rightarrow\infty$ and $\gamma_{n}=o\left(
n^{1/4}\right)  $ as $n\rightarrow\infty$.\label{P}
\end{hp}

\setcounter{hp}{12}

\begin{hp}
[Model]\label{M}The processes $\{u_{t,n}\}$, the model $m(X_{t};\theta)=u_{t}$
and the estimators $\left\{  \widehat{\theta}_{t}\right\}  $ satisfy the
following conditions:

(i) There is a sequence $\left\{  \theta_{n}\right\}  $, with $\theta
_{n}=\theta_{0}$ for all $n$ under $\mathcal{H}_{0}$, such that
\begin{equation}
\left\{  \left(  n^{1/2}\left(  \widehat{\theta}_{[ns]}-\theta_{n}\right)
^{\prime},n^{-1/2}\sum_{t=1}^{[ns]}\left(  u_{t,n}u_{t-1,n}-\mathbb{E}\left[
u_{t,n}u_{t-1,n}\right]  \right)  \right)  ^{\prime},s\in\left[  0,1\right]
\right\}  \label{M1}%
\end{equation}
$D_{[0,1]}$-converges in distribution to a Brownian motion with a full rank
covariance matrix.

(ii) The residual function admits a second order expansion $u_{t}\left(
\theta\right)  =u_{t,n}+(\theta-\theta_{n})^{\prime}u_{t,n}^{(1)}+\left(
\theta-\theta_{n}\right)  ^{\prime}u_{t,n}^{(2)}\left(  \theta-\theta
_{n}\right)  +\mathfrak{r}_{t,n}\left(  \theta\right)  $ where, for any
$C>0$,
\begin{equation}
\sup_{t\in\left[  1,n\right]  }\sup_{\theta;\left\Vert \theta-\theta
_{n}\right\Vert \leq Cn^{-1/2}}\left\vert \mathfrak{r}_{t,n}\left(
\theta\right)  \right\vert =o_{\mathbb{P}}\left(  \frac{1}{n}\right)
\label{M2}%
\end{equation}
and, for each $n$, $\{u_{t,n},u_{t,n}^{(1)},u_{t,n}^{(2)}\}$ is a stationary
process with $\mathbb{E}^{1/2}\left[  \left\Vert a_{t}\right\Vert ^{2}\right]
\leq C_{4}$, $\left\{  a_{t}\right\}  $ being successively $\left\{
u_{t,n}^{(1)}\right\}  $, $\left\{  u_{t,n}^{(2)}\right\}  $ $\left\{
u_{t,n}^{2}\right\}  $, $\left\{  u_{t,n}u_{t,n}^{(1)}\right\}  $, $\left\{
u_{t,n}^{\left(  1\right)  }u_{t,n}^{(1)^{\prime}}\right\}  $, $\left\{
u_{t,n}u_{t,n}^{(2)}\right\}  $, and where \newline$\sum_{j=-\infty}^{\infty
}\mathbb{E}\left[  \left\Vert u_{t-j,n}^{\left(  1\right)  }u_{t,n}\right\Vert
^{2}\right]  \leq C_{5}$, $\sup_{j\in\mathbb{Z}}\mathbb{E}\left[  \left\Vert
n^{-1/2}\sum_{t=j+1}^{n}\left(  u_{t-j,n}^{\left(  1\right)  }u_{t,n}%
-\mathbb{E}[u_{t-j,n}^{\left(  1\right)  }u_{t,n}]\right)  \right\Vert
^{2}\right]  \leq C_{6}$, $\sup_{j\in\mathbb{Z}}\mathbb{E}\left[  \left\Vert
u_{t,n}^{\left(  1\right)  }u_{t,n}u_{t-j,n}^{2}\right\Vert \right]  \leq
C_{7}$, and \newline$\sup_{j\in\mathbb{Z}}\mathbb{E}\left[  \left\Vert
n^{-1/2}\sum_{t=j+1}^{n}\left(  u_{t,n}^{\left(  1\right)  }u_{t,n}%
u_{t-j,n}^{2}-\mathbb{E}[u_{t,n}^{\left(  1\right)  }u_{t,n}u_{t-j,n}%
^{2}]\right)  \right\Vert ^{2}\right]  \leq C_{8}.$
\end{hp}

The compact sets $\left[  0,1/2\right]  $ and $\left[  0,1\right]  $ in
Assumption \ref{Kernel} are somehow arbitrary and can be replaced by any
nested compact intervals. Note however that Assumption \ref{Kernel} forbids
the use of the Daniell kernel $K\left(  x\right)  =\sin\left(  x\right)  /x$
due to the nonincreasing function and bounded support conditions.

Assumption \ref{Reg} imposes a polynomial decay on the coefficients
$\delta_{12a}\left(  j\right)  $, a condition which is weaker than the
exponential rate assumed in Shao (2011b). Note that in Assumption \ref{P} the
order of $\overline{p}_{n}$ can come closer to $n^{1/2}$ when $a$ is high,
that is when $u_{t}$ has finite moments of higher order. Under Assumption
\ref{Reg}, $\left\{  u_{t,n}\right\}  $\ must have finite moments of order
twelve at least. This is mostly needed for a proof of Theorem \ref{Level}
below based on Lindeberg substitution method, see Pollard (2002, p. 179),
which uses moment bounds as the Cauchy-Schwarz inequality $\mathbb{E}\left[
\left(  u_{t}^{2}u_{t+j}^{2}\right)  ^{3}\right]  \leq\mathbb{E}\left[
u_{t}^{12}\right]  $. Since implementing the proposed data-driven tests with a
large $\overline{p}_{n}$ would in principle allow us to detect a wider class
of alternatives, Assumption \ref{P}, which plays an important role under the
null in our proofs, may be too restrictive. Our simulation experiments indeed
suggest that Assumption \ref{P} can be weakened when focusing on white noise
processes of practical relevance since the order $\overline{p}_{n}\asymp n$
gives good results for various white noise processes of practical interest. On
the other hand, choosing a smaller $\overline{p}_{n}$ still gives a good
power, see comments on Table 5 at the end of the simulation experiments section.

When $\left\{  u_{t}\right\}  $ is observed, Assumption \ref{M} is equivalent
to Assumption 1 of Lobato (2001) and the FCLT for $n^{-1/2}\sum_{t=1}%
^{[ns]}\left(  u_{t}u_{t-1}-\mathbb{E}\left[  u_{t}u_{t-1}\right]  \right)  $
is a consequence of Assumption \ref{Reg} and the FCLT of Wu (2007). Assumption
\ref{M} is easily verified for simple linear models and OLS estimation where
$u_{t,n}^{(2)}$ and $\mathfrak{r}_{t,n}$ can be set to $0$. Assumption
\ref{M}-(i) is a shortened version of Assumptions B1 and A2 of Kuan and Lee
(2006) who employ a standard linear expansion $n^{1/2}\left(  \widehat{\theta
}-\theta_{n}\right)  =n^{-1/2}\sum_{t=1}^{n}\ell_{t}+o_{\mathbb{P}}\left(
1\right)  $ to show that (\ref{M1}) satisfies a functional central limit
theorem (FCLT) called for in \ref{M}-(i). The FCLT is mostly used under
$\mathcal{H}_{0}$ to show that $\mathbb{P}\left(  \widehat{S}_{1}\geq
\widehat{z}\left(  \alpha\right)  \right)  \rightarrow\alpha$ and
$\mathbb{P}\left(  \widehat{S}_{1}^{\ast}\geq\widehat{z}^{\ast}\left(
\alpha\right)  \right)  \rightarrow\alpha$ in the case of residuals. The
full-rank FCLT condition in Assumption \ref{M}-(i) implies certain
restrictions. For example, for a correctly specified $AR(1)$ model
$X_{t}-\theta X_{t-1}=u_{t}$, the case of $\theta=0$ is ruled out, a value of
the parameter which would in principle be excluded when considering such an
$AR(1)$\ specification. Theorem \ref{AltCV} at the end of the next section
explains how to overcome this issue with an alternative choice of critical
values when Assumption \ref{M}-(i) is too restrictive. The next section
describes some suitable theoretical requirements for the penalty sequence
$\gamma_{n}$ while the simulation section proposes a calibration of
$\gamma_{n}$ which gives good results for various white noise processes and alternatives.

\section{Asymptotic level and consistency\label{Main results}}

\setcounter{equation}{0}An important issue in the construction of the test
(\ref{Test}) is the choice of the penalty sequence. Choosing $\gamma_{n}%
$\ large enough implies that $\widehat{p}$ stays close to $1$ and so the test
statistic $\widehat{S}_{\widehat{p}}$ remains close to $\widehat{S}_{1}$.
Hence, on the one hand, using large $\gamma_{n}$ ensures that the level of the
test is close to its nominal size. On the other hand, a large $\gamma_{n}$ may
substantially limit the power of the test since the statistic $\widehat
{S}_{\widehat{p}}$ would not differ from $\widehat{S}_{1}$. The trade-off
between size and power is addressed by Theorem \ref{Penalty lower bound}\ and
Theorem \ref{Sparse}.

Consider first the properties of the test under the null hypothesis. The
following theorem gives a lower bound for $\gamma_{n}$ which ensures that
$\widehat{p}=1$ asymptotically so that the test is asymptotically of level
$\alpha$.

\begin{thm}
\label{Penalty lower bound}Let Assumptions \ref{Kernel} , \ref{M}, \ref{P} and
\ref{Reg} hold. If the penalty sequence $\{\gamma_{n},n\geq1\}$ satisfies
\begin{equation}
\gamma_{n}\geq\left(  1+\epsilon\right)  \left(  2\ln\ln n\right)
^{1/2}\text{\quad for some }\epsilon>0, \label{Gam}%
\end{equation}
then under $\mathcal{H}_{0}$, $\lim_{n\rightarrow\infty}\mathbb{P}\left(
\widehat{p}=1\right)  =1$ and the test (\ref{Test}) is asymptotically of level
$\alpha$. \label{Level}
\end{thm}

\noindent Under the null hypothesis, the selected order $\widehat{p}$ is
asymptotically equal to $1$. It follows that $\widehat{S}_{\widehat{p}%
}=\widehat{S}_{1}+o_{\mathbb{P}}\left(  1\right)  $ and that critical values
(\ref{Zlob}) or (\ref{Zest}) guarantee that the test is asymptotically of
level $\alpha$. A key result is therefore that $\lim_{n\rightarrow\infty
}\mathbb{P}\left(  \widehat{p}=1\right)  =1$ holds under various white noise
models and observed $u_{t}$ or residuals $\widehat{u}_{t}$. That the
estimation has no impact asymptotically follows from (\ref{Gam}) which imposes
$\gamma_{n}\rightarrow\infty$. When $\widehat{\theta}$ is $\sqrt{n}%
$-consistent, estimating the residuals gives test statistics satisfying
\[
\widehat{S}_{p}=n\sum_{j=1}^{p}\left(  \frac{1}{n}\sum_{t=1}^{n-j}u_{t}%
u_{t+j}\right)  ^{2}+O_{\mathbb{P}}\left(  1\right)
\]
uniformly in $p$. The fact that the remainder term $O_{\mathbb{P}}\left(
1\right)  $ is negligible compared to $\gamma_{n}$ is a crucial element in
showing that the asymptotic behavior of $\widehat{p}$ is not affected by the
estimation under the null. The divergence of $\gamma_{n}$ is also important to
account for the fact that the standardization $E_{\Delta}\left(  p\right)  $
and $V_{\Delta}\left(  p\right)  $ are only valid when $p\rightarrow\infty$
since $\gamma_{n}\rightarrow\infty$ imposes that either $\widehat{p}=1$ or
$\widehat{p}$ diverges because (\ref{Rej}) cannot hold for finite $p>1$.
Compared to the existing adaptive results of Horowitz and Spokoiny (2001),
Guerre and Lavergne (2005), Guay and Guerre (2006) or Chen and Gao (2007), an
important technical contribution of our paper is that Theorem \ref{Level}
holds without assuming that the set of admissible $p$ is a power set $\left\{
a^{j},j\in\mathbb{N}\right\}  $, $a>1$.

Another important finding is that the penalty sequence $\gamma_{n}$ can
diverge with the low order $\left(  \ln\ln n\right)  ^{1/2}$ allowed by
(\ref{Gam}). This contrasts with the larger order $\ln n$ used in the BIC
selection procedure and in the corresponding data-driven tests. In view of the
potential negative impact of a large $\gamma_{n}$ on the power of the test, it
is worth asking if the lower bound (\ref{Gam}) can be improved, that is if
$\mathbb{P}\left(  \widehat{p}=1\right)  \rightarrow1$ would be ensured for
even lower values of penalty term $\gamma_{n}$. The proof suggests that this
is not the case. The main argument is based on expression%
\begin{equation}
\mathbb{P}\left(  \widehat{p}\neq1\right)  =\mathbb{P}\left(  \max
_{p\in\left[  2,\overline{p}_{n}\right]  }\left(  \frac{\left(  \hat{S}%
_{p}-\hat{S}_{1}\right)  /\hat{R}_{0}^{2}-E_{\Delta}(p)}{V_{\Delta}%
(p)}\right)  \geq\gamma_{n}\right)  \label{Hatpnot}%
\end{equation}
for the probability of not selecting $1$. It can be seen from the proof of
Theorem \ref{Penalty lower bound}\ that, for the Box-Pierce version of the
test, the right-hand side of (\ref{Hatpnot}) asymptotically behaves like the
maximum of standardized partial sums whose exact order is $\left(  2\ln\ln
n\right)  ^{1/2}$, see (B.38) in the Supplementary Material. Hence the bound
(\ref{Gam}) is optimal to achieve $\mathbb{P}\left(  \widehat{p}=1\right)
\rightarrow1$.

Let us now turn to the detection properties of the test. Recall that the
covariance of the alternative may depend on the sample size so that
$R_{j}=R_{j,n}$ may go to $0$ when $n$ increases. The new class of
alternatives is defined similarly to (\ref{BPcons}) in the introduction
section. Consider first a sequence $\rho_{n}\rightarrow0$ and a lag order
$P_{n}$. An important indicator for detection of alternatives is the number of
correlations above $\rho_{n}$,%
\begin{equation}
N_{n}=N_{n}\left(  P_{n},\rho_{n}\right)  =\#\left\{  |R_{j}/R_{0}|\geq
\rho_{n},\quad1\leq j\leq P_{n}\right\}  . \label{Sparse2}%
\end{equation}
The next theorem gives a detection condition on $N_{n}$, $P_{n}$ and $\rho
_{n}$.

\begin{thm}
\label{Sparse} Suppose Assumptions \ref{Kernel}, \ref{M}, \ref{Reg} and
\ref{P} hold. There exists a constant $\kappa_{\ast}>0$ such that the test
(\ref{Test}) is consistent against all alternatives $\{u_{t}\}$ satisfying,
for some $\rho_{n}>0$ and $P_{n}\in\left[  1,\overline{p}_{n}/2\right]  $,%
\begin{equation}
n^{1/2}\left(  \frac{N_{n}}{\gamma_{n}P_{n}^{1/2}}\right)  ^{1/2}\rho_{n}%
\geq\kappa_{\ast}. \label{Sparse3}%
\end{equation}

\end{thm}

\noindent Condition (\ref{Sparse3}) is similar to the detection condition
(\ref{BPcons}) required for consistency of the Box-Pierce test (\ref{BPcrit}).
However a key difference between the two conditions is that while in
(\ref{BPcons}) the lag order $p_{n}$ is assumed known and is used in the
construction of the test statistic, in (\ref{Sparse3}) the lag order $P_{n}$
in (\ref{Sparse3}) is unknown. This illustrates the adaptive capability of the
new test. A second important difference between (\ref{BPcons}) and
(\ref{Sparse3}) is that the latter involves penalty sequence $\gamma_{n}$. For
given $P_{n}$ and $N_{n}$ detection condition (\ref{Sparse3}) admits rate
$\rho_{n}^{\ast}$ satisfying%
\begin{equation}
\rho_{n}^{\ast}\asymp\frac{1}{n^{1/2}}\left(  \frac{\gamma_{n}P_{n}^{1/2}%
}{N_{n}}\right)  ^{1/2}. \label{Bestrho}%
\end{equation}
Rate $\rho_{n}^{\ast}$ in (\ref{Bestrho}) deteriorates with the penalty
sequence. Condition (\ref{Sparse3}) thus demonstrates the potential negative
impact of the penalty sequence on the power of the test. This impact can also
be seen from proof of Theorem \ref{Sparse} which uses the fact that the test
(\ref{Test}) rejects the null whenever%
\begin{equation}
\frac{\widehat{S}_{p}-\widehat{R}_{0}^{2}E\left(  p\right)  }{\widehat{R}%
_{0}^{2}V_{\Delta}\left(  p\right)  }\geq\gamma_{n}+\frac{\widehat{z}\left(
\alpha\right)  }{\widehat{R}_{0}^{2}V_{\Delta}\left(  p\right)  }\text{ for
some }p\in\left[  2,\overline{p}_{n}\right]  . \label{Detect}%
\end{equation}
For the alternatives for which (\ref{Detect}) only holds for $p\rightarrow
\infty$ so that $V_{\Delta}\left(  p\right)  \rightarrow\infty$,
(\ref{Detect}) suggests that $\gamma_{n}$ may be more important than the
critical value $\widehat{z}\left(  \alpha\right)  $ regarding detection.

Two special cases of (\ref{Bestrho}) are worth mentioning. First, the
situation where $\lim_{n\rightarrow\infty}\gamma_{n}P_{n}^{1/2}/N_{n}=0$\ is
of special interest since (\ref{Bestrho}) shows that the test can detect
correlation coefficients converging to $0$ at a rate that is faster than the
parametric rate $n^{-1/2}$. The best possible rate in this case is $\rho
_{n}^{\ast}\asymp\gamma_{n}^{1/2}/\left(  nP_{n}^{1/2}\right)  ^{1/2}$ which
is achieved for \textquotedblleft saturated\textquotedblright\ alternatives
with $N_{n}\asymp P_{n}$. Second, a less favorable case corresponds to more
sparse correlation coefficients satisfying $\lim_{n\rightarrow\infty}%
\gamma_{n}P_{n}^{1/2}/N_{n}=\infty$. In this case (\ref{Bestrho}) does not
allow for correlation coefficients converging to $0$ at the rate of $n^{-1/2}%
$. This case has been covered by Donoho and Jin (2004) for a theoretical model
where a known number $P_{n}$ of independent Gaussian variables with mean
$n\left(  R_{j}/R_{0}\right)  ^{2}$ and variance $1$ is observed. These
authors show that in such a setup the best possible detection rate is
$\rho_{n}=\left(  \ln n/n\right)  ^{1/2}$, a rate which is achieved by the
maximum white noise test of Xiao and Wu (2011). This suggests that our test
may not be optimal when $\lim_{n\rightarrow\infty}\gamma_{n}P_{n}^{1/2}%
/N_{n}=\infty$. However, it is shown in Proposition \ref{WildCvM} in Section
\ref{Optimal} below that the test of Xiao and Wu (2011), unlike our test, does
not detect moderately sparse alternatives satisfying (\ref{Bestrho}) with
$\lim_{n\rightarrow\infty}\gamma_{n}P_{n}^{1/2}/N_{n}=0$ and $\gamma_{n}%
\asymp\left(  2\ln\ln n\right)  ^{1/2}$.

We conclude this section with two extensions of our main results. The first
extension shows that the test derived from (\ref{Sstar}) and (\ref{Hatpstar})
has similar properties as the test (\ref{Test}).

\begin{thm}
\label{Extension}Suppose Assumptions \ref{Kernel}, \ref{M}, \ref{Reg} and
\ref{P} hold. Then $\mathbb{P}\left(  \widehat{p}^{\ast}=1\right)
\rightarrow1$ under $\mathcal{H}_{0}$ and the test which rejects the null when
$\widehat{S}_{\widehat{p}^{\ast}}^{\ast}\geq\widehat{z}^{\ast}\left(
\alpha\right)  $ is asymptotically of level $\alpha$. It also detects the
alternatives satisfying (\ref{Sparse3}) in Theorem \ref{Sparse} for a large
enough $\kappa_{\ast}$.
\end{thm}

The second extension is useful in the case of residuals when the full-rank
FCLT condition in Assumption \ref{M}-(i) is too restrictive so that the
critical value $\widehat{z}_{KL}\left(  \alpha\right)  $ in (\ref{Zest})
cannot be used. Suppose that an additional test statistic $\widehat{T}_{n}%
$\ with critical values $\widehat{t}_{n}\left(  \alpha\right)  $ satisfying
$\lim_{n\rightarrow\infty}\mathbb{P}\left(  \widehat{T}_{n}\geq\widehat{t}%
_{n}\left(  \alpha\right)  \right)  =\alpha$ under the null is available.
Consider the critical value%
\begin{equation}
\widehat{c}_{n}^{\ast}\left(  \alpha\right)  =\widehat{S}_{1}^{\ast}%
-\widehat{T}_{n}+\widehat{t}_{n}\left(  \alpha\right)  . \label{CVc}%
\end{equation}

\begin{thm}
\label{AltCV} Suppose that Assumptions K, R and P hold, as Assumption M-(ii)
with $\sqrt{n}\left(  \widehat{\theta}-\theta_{n}\right)  =O_{\mathbb{P}%
}\left(  1\right)  $ where the deterministic sequence $\left\{  \theta
_{n}\right\}  $ is such that $\theta_{n}=\theta_{0}$ for all $n$ under $H_{0}%
$. Suppose also that (A0) $\lim_{n\rightarrow\infty}\mathbb{P}\left(
\widehat{T}_{n}\geq\widehat{t}_{n}\left(  \alpha\right)  \right)  =\alpha$
under $\mathcal{H}_{0}$ and (A1) $\widehat{c}_{n}\left(  \alpha\right)  \leq
O_{\mathbb{P}}\left(  \gamma_{n}\right)  $ under the considered alternative.
Then the test which rejects the null when $\widehat{S}_{\widehat{p}^{\ast}%
}^{\ast}\geq\widehat{c}_{n}\left(  \alpha\right)  $ is asymptotically of level
$\alpha$ and detects the alternatives satisfying the condition (\ref{Sparse3})
of Theorem \ref{Sparse} for a sufficiently large $\kappa_{\ast}$. Moreover and
even if (A1) does not hold, we have under the alternative and for any sample
size $n$,%
\begin{equation}
\mathbb{P}\left(  \widehat{S}_{\widehat{p}^{\ast}}^{\ast}\geq\widehat{c}%
_{n}\left(  \alpha\right)  \right)  \geq\mathbb{P}\left(  \widehat{T}_{n}%
\geq\widehat{t}_{n}\left(  \alpha\right)  \right)  . \label{Improv}%
\end{equation}

\end{thm}

Condition (A1), which allows for $\widehat{c}_{n}\left(  \alpha\right)
\overset{\mathbb{P}}{\mathbb{\rightarrow}}-\infty$, means, when $\widehat
{t}_{n}\left(  \alpha\right)  =O_{\mathbb{P}}\left(  1\right)  $ as usual,
that $\widehat{T}_{n}$ diverges at least as fast as $\widehat{S}_{1}^{\ast}$
or that both lack power against the considered alternative and are
$O_{\mathbb{P}}\left(  1\right)  $. The bound (\ref{Improv}) means that the
data-driven test is at least as powerful than the test based on $\widehat
{T}_{n}$. As a consequence of (\ref{Improv}), the test $\widehat{S}%
_{\widehat{p}^{\ast}}^{\ast}\geq\widehat{z}^{\ast}\left(  \alpha\right)  $ is
as least as powerful as $\widehat{S}_{1}^{\ast}\geq\widehat{z}^{\ast}\left(
\alpha\right)  $, $\widehat{z}^{\ast}\left(  \alpha\right)  $ as in
(\ref{Zstar}). The use of the critical value (\ref{CVc}) can give a
data-driven test whose power properties can be tailored to be optimal against
some specific alternatives by a proper choice of a corresponding optimal
$\widehat{T}_{n}$. Examples of test statistic $\widehat{T}_{n}$\ which does
not require Assumption \ref{M}-(i) can be found in Delgado and Velasco (2012)
and Francq et al. (2005). Delgado and Velasco (2012) propose a Box-Pierce
statistic corrected for estimation with an elegant general approach and some
parametric optimality properties under Gaussianity whereas Francq et al.
(2005) is more specific to ARMA specifications.

\setcounter{equation}{0}

\section{Adaptive rate-optimality and comparisons with other tests
\label{Optimal}}

\setcounter{equation}{0} While Theorem \ref{Level} gives the lower bound
(\ref{Gam}) of order $\left(  2\ln\ln n\right)  ^{1/2}$ for the penalty
sequence $\gamma_{n}$ that is necessary to ensure that the test is
asymptotically of level $\alpha$, Theorem \ref{Sparse} suggests that
increasing $\gamma_{n}$ can impair the power of the test. Hence a good
compromise for the choice of the penalty sequence suitable both under
$\mathcal{H}_{0}$ and $\mathcal{H}_{1}$ is $\gamma_{n}\asymp\left(  2\ln\ln
n\right)  ^{1/2}$. Once this choice is made one may ask if the resulting test
is the best possible in the sense that there is no other test that can detect
alternatives satisfying a condition less restrictive than (\ref{Sparse3}),
when $\kappa_{\ast}=\kappa_{n}\rightarrow0$ is allowed. The absence of a
better test is the so called adaptive rate-optimality. The next theorem
establishes adaptive rate-optimality for alternatives satisfying
$\lim_{n\rightarrow\infty}\gamma_{n}P_{n}^{1/2}/N_{n}=0$.\footnote{As
discussed when introducing approximation (\ref{Bestrho}), the test
(\ref{Test}) is not optimal for detection of sparse alternatives with
$\lim_{n\rightarrow\infty}\gamma_{n}P_{n}^{1/2}/N_{n}=\infty$ which are not
considered here.}

\begin{thm}
\label{Sparseopt} Let $u_{t}$ be observed. For any sequence $\kappa
_{n}\rightarrow0$, there exists a sequence of alternatives $\left\{
u_{t}\right\}  $ such that, for some $P_{n}\in\left[  1,\overline{p}%
_{n}\right]  $ and $\rho_{n}>0$ with
\[
\rho_{n}\geq\frac{\kappa_{n}}{n^{1/2}}\left(  \frac{\left(  2\ln\ln n\right)
^{1/2}P_{n}^{1/2}}{N_{n}}\right)  ^{1/2},\quad\quad\lim_{n\rightarrow\infty
}\frac{\left(  2\ln\ln n\right)  ^{1/2}P_{n}^{1/2}}{N_{n}}=0,
\]
such that the other assumptions of Theorem \ref{Sparse} are satisfied, but
that cannot be detected by any possible asymptotically $\alpha$-level test.
\end{thm}

%

\noindent
Hence, when $\gamma_{n}\asymp\left(  2\ln\ln n\right)  ^{1/2}$, it is not
possible to improve on the detection condition (\ref{Sparse3}) and the rate
$\rho_{n}^{\ast}$ in (\ref{Bestrho}) is optimal. We now give an example of
alternatives which are detected by the test (\ref{Test}) but not by other
popular tests.\ Consider the following high-order moving average process,%
\begin{equation}
u_{t}=u_{t,n}=\varepsilon_{t}+\frac{\nu\gamma_{n}^{1/2}}{n^{1/2}P_{n}^{1/4}%
}\sum_{k=1}^{P_{n}}\psi_{k}\varepsilon_{t-k},\text{\quad}\sum_{k=1}^{P_{n}%
}\psi_{k}^{2}=O(P_{n}),\quad\lim_{n\rightarrow\infty}P_{n}=\infty,
\label{Wildalt}%
\end{equation}
where $\left\{  \varepsilon_{t}\right\}  $ is a strong white noise with
variance $\sigma^{2}$, $\nu$ is a scaling constant and $\gamma_{n}%
\asymp\left(  2\ln\ln n\right)  ^{1/2}$. This alternative has moving average
coefficients of order $\gamma_{n}^{1/2}/\left(  n^{1/2}P_{n}^{1/4}\right)
=o\left(  n^{-1/2}\right)  $ provided $P_{n}$ diverges at a polynomial rate.
Hence short term shocks have statistically negligible impact. However when
$\psi_{k}=1$ for all $k$, the long term multiplier of (\ref{Wildalt}) is equal
to $\nu\left(  \gamma_{n}P_{n}^{3/2}/n\right)  ^{1/2}$ which is of larger
order than $n^{-1/2}$. The following lemma describes the covariance function
and conditional expectation of the alternative (\ref{Wildalt}).

\begin{lem}
\label{Wildlem} If $P_{n}=o((n/\gamma_{n})^{2/3})$ and $\lim_{n\rightarrow
\infty}\left(  \gamma_{n}/n\right)  =0$ then the alternative $\left\{
u_{t}\right\}  $ in (\ref{Wildalt}) satisfies $R_{0}=\sigma^{2}\left(
1+O\left(  \gamma_{n}P_{n}^{1/2}/n\right)  \right)  $ and, uniformly in
$j\in\left[  1,P_{n}\right]  $,%
\[
R_{j}=\frac{\nu\gamma_{n}^{1/2}}{n^{1/2}P_{n}^{1/4}}\psi_{j}\sigma
^{2}+o\left(  \frac{\gamma_{n}^{1/2}}{n^{1/2}P_{n}^{1/4}}\right)  .
\]
Moreover%
\[
\mathbb{E}\left[  u_{t}|u_{t-k},k\geq1\right]  =\frac{\nu\gamma_{n}^{1/2}%
}{n^{1/2}P_{n}^{1/4}}\sum_{k=1}^{P_{n}}\psi_{k}u_{t-k}+O_{\mathbb{P}}\left(
\frac{\gamma_{n}P_{n}}{n}\right)  .
\]

\end{lem}

\noindent Hence a distinctive feature of the alternative (\ref{Wildalt})
when\textbf{\ }$\max_{1\leq k\leq P_{n}}\left\vert \psi_{k}\right\vert
=O\left(  1\right)  $\textbf{\ }is that $\max_{j\geq1}\left\vert
R_{j}\right\vert =o\left(  n^{-1/2}\right)  $\ provided $P_{n}/\gamma_{n}%
^{2}\rightarrow\infty$. The expression of $\mathbb{E}\left[  u_{t}%
|u_{t-k},k\geq1\right]  $ reveals that $u_{t}$ can be very difficult to
forecast since the coefficients of the lagged variables are all $o\left(
n^{-1/2}\right)  $ provided $P_{n}=o\left(  n^{1/2}/\gamma_{n}\right)  $. This
suggests that such a process will be seen in practice as a martingale
difference when using standard statistical tools. This may be a relevant
example of alternatives in economical or financial contexts where arbitrage occurs.

We show in Proposition \ref{WildCvM}\ below that the new tests detect these
alternatives but that this is not the case for three tests based on the
following test statistics,%
\begin{equation}
W_{n}=b_{n}\left(  n^{1/2}\max_{j\in\left[  1,J_{n}\right]  }\left\vert
\frac{\widehat{R_{j}}}{\widehat{\tau}_{j}}\right\vert -b_{n}\right)
,\quad\text{where }b_{n}=\left(  2\ln J_{n}-\ln\ln J_{n}-\ln\left(
4\pi\right)  \right)  ^{1/2}, \label{Wu}%
\end{equation}%
\begin{equation}
CvM_{n}=\frac{n}{\pi^{2}}\sum_{j=1}^{J_{n}}\frac{\widehat{R}_{j}^{2}}%
{j^{2}\widehat{\tau}_{j}^{2}}, \label{CvM}%
\end{equation}%
\begin{align}
EL_{n}  &  =\widehat{BP}_{\widehat{p}_{EL}^{\ast}}^{\ast},\quad\widehat
{p}_{EL}^{\ast}=\arg\max_{p\in\left[  1,J_{n}\right]  }\left\{  \widehat
{BP}_{p}^{\ast}-\widehat{\gamma}_{EL}^{\ast}p\right\}  \text{ where}%
\label{EL}\\
&  \widehat{\gamma}_{EL}^{\ast}=\left\{
\begin{array}
[c]{ll}%
\ln n & \text{if }n^{1/2}\max_{j\in\left[  1,J_{n}\right]  }\left\vert
\frac{\widehat{R}_{j}}{\widehat{\tau}_{j}}\right\vert \leq\left(  2.4\ln
n\right)  ^{1/2},\\
2 & \text{otherwise.}%
\end{array}
\right. \nonumber
\end{align}
Statistic $W_{n}$ in (\ref{Wu}) is studied in Xiao and Wu (2011) who show that
$W_{n}$ asymptotically has an extreme value distribution. The statistic
$CvM_{n}$ in (\ref{CvM}), due to Deo (2000) for observed $u_{t}$, is a version
of the Cram\'{e}r-von Mises test of Durlauf (1991) partially corrected for
heteroskedasticity. Test statistic $EL_{n}$ has been introduced by Escanciano
and Lobato (2009) for observed $u_{t}$ and a fixed $J_{n}$. As in our test,
the order $\widehat{p}_{EL}^{\ast}$\ selected by Escanciano and Lobato (2009)
is asymptotically equal to $1$ under $\mathcal{H}_{0}$ and similar critical
values can be used. To show that tests (\ref{Wu})--(\ref{EL})\ do not detect
alternatives with small correlation coefficients, it is sufficient to consider
a Gaussian null hypothesis $G_{0}$ under which $\left\{  u_{t}\right\}  $ is a
Gaussian white noise process $\left\{  \varepsilon_{t}\right\}  $ with
variance $\sigma^{2}$ against an alternative $G_{1}$ under which $\left\{
u_{t}\right\}  $ is given by (\ref{Wildalt}) with Gaussian i.i.d. $\left\{
\varepsilon_{t}\right\}  $, $\sum_{k=1}^{P_{n}}\psi_{k}^{2}=O(P_{n})$,
$\max_{1\leq k\leq P_{n}}\left\vert \psi_{k}\right\vert =O\left(  1\right)  ,$
$\min_{1\leq k\leq P_{n}}\left\vert \psi_{k}\sigma^{2}\right\vert \geq1$,
$\nu>0$, $\gamma_{n}$ and $P_{n}\rightarrow\infty$ with $\gamma_{n}%
/P_{n}^{1/2}=o\left(  1/\ln n\right)  $ and $P_{n}=O\left(  \left(
n/\gamma_{n}\right)  ^{1/14}\right)  \leq\overline{p}_{n}/2$ and $\gamma
_{n}\asymp\left(  2\ln\ln n\right)  ^{1/2}$ satisfies (\ref{Gam}). We assume
that $J_{n}=O\left(  n^{1/2}\right)  $.

\begin{prop}
Let $u_{t}$ be observed. Suppose that Assumptions \ref{Kernel} and \ref{P}
hold. For $\nu$ large enough, the alternative $G_{1}$ as above satisfies
(\ref{Sparse3}) and

\textit{(i)} the test (\ref{Test}) and its $\widehat{S}_{\widehat{p}^{\ast}%
}^{\ast}$ version consistently detect $G_{1}$. By contrast,

\textit{(ii)} statistics $W_{n}$, $CvM_{n}$ and $EL_{n}$ have the same
asymptotic distribution under $G_{0}$ and $G_{1}$ and the corresponding tests
are therefore not consistent.\label{WildCvM}
\end{prop}

\noindent Proposition \ref{WildCvM}-(ii) implies that tests based on $W_{n}$,
$CvM_{n}$ or $EL_{n}$ are not adaptive rate-optimal. Let $\widehat{R}%
_{0,j}/\widehat{\tau}_{0,j}$\ and $\widehat{R}_{1,j}/\widehat{\tau}_{1,j}$\ be
the standardized sample covariance computed under $G_{0}$\ and $G_{1}%
$\ respectively. It is established in the proof of Proposition \ref{WildCvM}
that%
\begin{equation}
\max_{j\in\left[  1,J_{n}\right]  }\left\vert \frac{\widehat{R}_{0,j}%
}{\widehat{\tau}_{0,j}}-\frac{\widehat{R}_{1,j}}{\widehat{\tau}_{1,j}%
}\right\vert =o_{\mathbb{P}}\left(  \frac{1}{\left(  n\log n\right)  ^{1/2}%
}\right)  , \label{G01cov}%
\end{equation}
which implies that tests based on $W_{n}$\ and $CvM_{n}$\ are not consistent.
The case of $EL_{n}$\ test is a bit more involved but, due to its penalty
scheme, this test statistic is asymptotically equal to $\widehat{BP}_{1}%
^{\ast}$ under the null and the alternative so that it cannot detect $G_{1}$
by (\ref{G01cov}).

\section{Simulation experiments\label{Simulation experiments}}

\setcounter{equation}{0} Our simulation experiments aim to propose a valid
penalty sequence $\gamma_{n}$ to be tested under various strong and weak white
noise processes and under various alternatives. Since preliminary experiments
have shown that the test statistic $\widehat{S}_{\widehat{p}}$ may yield an
oversized test for some practically relevant white noise processes, we
consider the test based on $\widehat{S}_{\widehat{p}^{\ast}}^{\ast}$ as in
(\ref{Sstar}) and (\ref{Hatpstar}). To investigate the impact of choosing a
large $\overline{p}_{n}$ latter on we allow for all possible orders, setting
$\overline{p}_{n}=n-1.$ We consider two kernels. The first is $K\left(
t\right)  =\mathbb{I}\left(  t\in\left[  0,1\right]  \right)  $ which gives
the Box-Pierce statistic so that the corresponding tests are labelled $BP$.
The second uses the Parzen kernel%
\[
\mathsf{k}(t)=\left\{
\begin{array}
[c]{lll}%
1-6t^{2}+6|t|^{3}, &  & \left\vert t\right\vert \leq1/2,\\
2(1-|t|)^{3}, &  & 1/2<\left\vert t\right\vert \leq1,\\
0 &  & \text{otherwise.}%
\end{array}
\right.
\]
However since $\mathsf{k}\left(  1\right)  =0$ which would give a meaningless
$\widehat{S}_{1}^{\ast}=0$, we change $\mathsf{k}\left(  t\right)  $ into
$K\left(  t\right)  =\mathsf{k}\left(  t/2\right)  /\mathsf{k}\left(
1/2\right)  $ and label the corresponding tests as $Parz$. The critical values
(\ref{Zstar}) $\widehat{z}^{\ast}\left(  \alpha\right)  $, see also
(\ref{Zlob}) and (\ref{Zest}), use a power Parzen kernel $k\left(  t\right)
=\mathsf{k}^{32}\left(  t\right)  $, where the exponent 32 is has been
proposed by Lee (2007) whose simulations show that such a choice ensures that
the test with rejection region $n\widehat{R}_{1}^{2}\geq\widehat{z}^{\ast
}\left(  \alpha\right)  $ has good power properties. We consider $10\%$, $5\%$
and $1\%$ significance levels. A preliminary simulation experiment with
$100,000$ replications gives that the corresponding quantiles $z_{L}\left(
\alpha\right)  $ of (\ref{StdW}) used in $\widehat{z}^{\ast}\left(
\alpha\right)  $ are approximately $3.73$, $5.58$ and $10.97$ respectively,
which are in line with the critical values tabulated by Phillips et al. (2006,
Table 6).

The first experiment analyzes the sensitivity of the test to the penalty term
and aims to calibrate the proportionality constant for the penalty sequence.
The experiment investigates the behavior of the test under the null for
$\gamma_{n}=\gamma\left(  2\ln\ln\left(  n-2\right)  \right)  ^{1/2}$ where
the proportionality coefficient $\gamma$ ranges from $2.8$ to $3.8$. The
process $u_{t}$ is a white noise with the standard normal distribution. The
next table reports the simulated levels for $50,000$ replications and the
percentage $\%\left\{  \widehat{p}^{\ast}\neq1\right\}  $ of simulation draws
for which $\widehat{p}^{\ast}\neq1$, an important indicator in deciding
whether a difference between nominal and observed levels is due to a too small
$\gamma_{n}$ or improper critical values. In Table 1, `*' indicates an
oversized test, i.e. such that the null of a level smaller than the nominal
size is rejected at 1\% level by the one-sided test using the simulated level.%
\[
\text{\textbf{[INSERT TABLE 1 HERE]}}%
\]
A threshold value for the $BP$ test is $\gamma=3.4$ which ensures that the
observed sizes are close to the nominal sizes for $n=1,000$. The $Parz$ test
is slightly less oversized. Both tests have very similar value of $\%\left\{
\widehat{p}^{\ast}\neq1\right\}  $, well below $1\%$ for $\gamma=3.4$. In the
remaining simulation experiments $\gamma=3.4$ is used.

We introduce some benchmark tests. We compare our $BP$ and $Parz$ tests with
the data-driven test $EL$ based on the statistic $EL_{n}$ in (\ref{EL}) with
$J_{n}=n-1$ and the critical values of Lee (2007) in (\ref{Zstar}). We also
consider the Newey-West data-driven order $\widehat{p}_{IMSE}$ used by Hong
and Lee (2005) and the test statistic%
\[
\widehat{p}_{IMSE}=\left(  1\vee\widehat{C}^{1/5}\left(  f\right)  \right)
n^{1/5}\text{,\quad where\quad}\widehat{C}\left(  f\right)  =\frac
{144\sum_{j=-(n-1)}^{n-1}\mathsf{k}\left(  j/\widetilde{p}\right)
j^{4}\widehat{R}_{j}^{2}/\widehat{\tau}_{j}^{2}}{0.539285\sum_{j=-(n-1)}%
^{n-1}\mathsf{k}\left(  j/\widetilde{p}\right)  \widehat{R}_{j}^{2}%
/\widehat{\tau}_{j}^{2}},
\]%
\[
IMSE=\frac{\sum_{j=1}^{\widehat{p}_{IMSE}}\mathsf{k}^{2}\left(  j/\widehat
{p}_{IMSE}\right)  \left\{  \widehat{R}_{j}^{2}/\widehat{\tau}_{j}^{2}-\left(
1-\frac{j}{n}\right)  \right\}  }{\left(  2\sum_{j=1}^{\widehat{p}_{IMSE}%
}\mathsf{k}^{4}\left(  j/\widehat{p}_{IMSE}\right)  \left(  1-\frac{j}%
{n}\right)  ^{2}\right)  ^{1/2}},
\]
where $\mathsf{k}\left(  \cdot\right)  $ is the Parzen kernel and
$\widehat{\tau}_{j}^{2}$ is defined as in (\ref{Sstar}). In the definition of
$\widehat{p}_{IMSE}$, $\widetilde{p}$ is a pilot bandwidth that is set to
$\widetilde{p}=4(n/100)^{4/25}$. Note that $\widehat{C}\left(  f\right)  $
remains potentially stochastic under the null so that the null limit
distribution of\ $IMSE$ may differ from the standard normal distribution valid
for deterministic $p_{n}\rightarrow\infty$. We however follow common practice
and use standard normal critical values for the $IMSE$ test. The last
benchmark test, $CvM$, is based on Deo's (2000) Cram\'{e}r-von Mises statistic
$CvM_{n}$ in (\ref{CvM}) and uses the critical values tabulated by Anderson
and Darling (1952).

The first comparison under $\mathcal{H}_{0}$ is based on i.i.d. $\left\{
u_{t}\right\}  $ with the following distributions: standard normal (`Nor' in
Table 2), Student with three degrees of freedom (`Stud'), and centered chi
square with one degree of freedom (`Chi'). The Student distribution is used to
test the sensitivity of our test to the lack of higher-order moments while the
chi square distribution can reveal sensitivity to skewness.%
\[
\text{\textbf{[INSERT TABLE 2 HERE]}}%
\]
As in Table 1, the size of the $Parz$ test is slightly better than the size of
the $BP$ test but both perform well here, although $BP$ is slightly oversized
under the `Chi' white noise. The $EL$ and $IMSE$ are generally oversized with
strong size distortions for `Chi'. The $CvM$ test performs well except for the
`Chi' experiment.

The next experiment considers observed weak white noise $u_{t}$ or residuals
$\hat{u}_{t}$. Two conditional heteroskedastic martingale difference processes
are examined. The first is a GARCH(1,1) process with $u_{t}=s_{t}\zeta_{t}$
and $s_{t}^{2}=0.001+0.90s_{t-1}^{2}+0.05u_{t-1}^{2}$ where $\zeta_{t}$ are
i.i.d. standard normal innovations. The second process is an ARCH(1) process
with $u_{t}=s_{t}\zeta_{t}$ and $s_{t}^{2}=0.001+0.9u_{t-1}^{2}$. Due to the
ARCH coefficient larger than $\sqrt{1/3}=0.577$, $\mathbb{E}\left[  u_{t}%
^{4}\right]  =\infty$ and the tests are, in principle, not expected to behave
well in this experiment. The next three processes are uncorrelated but are not
martingale differences, so that the $CvM$ test is not expected to have a
correct size and is only reported here as a benchmark. The first, labelled
`Bilinear' in Table 3 below, is a bilinear model $u_{t}=\zeta_{t}%
+0.9\zeta_{t-1}u_{t-2}$. The second, labelled `No-MDS', is given by
$u_{t}=\zeta_{t-1}\zeta_{t-2}\left(  1+\zeta_{t-2}+\zeta_{t}\right)  $ and has
been examined by Lobato (2001). The third, `All-Pass', is an All-Pass
ARMA(1,1) process examined by Lobato, Nankervis and Savin (2002),
$u_{t}-0.5u_{t-1}=\zeta_{t}-\zeta_{t-1}/0.5$, where $\zeta_{t}$ i.i.d. and
have the Student distribution with $9$ degrees of freedom. Since the root of
the $MA$ part is the inverse of the $AR$ root, the resulting process is
uncorrelated but the $u_{t}$ are dependent due to non-Gaussian $\zeta_{t}$.
Finally, experiment `ARRes' examines residuals from the $AR\left(  1\right)  $
$y_{t}=0.8y_{t-1}+u_{t}$, $\widehat{u}_{t}=y_{t}-\widehat{\theta}y_{t-1}$,
$\widehat{\theta}=\sum_{t=0}^{n-1}y_{t}y_{t+1}/\sum_{t=0}^{n-1}y_{t}^{2}$. The
$BP$, $Parz$ and $EL$ tests are all adapted to the estimation effect thanks to
the use of the critical values $\widehat{z}^{\ast}\left(  \alpha\right)  $ of
(\ref{Zstar}). The critical values of the $IMSE$ and $CvM\,$do not account for
estimation of residuals and the corresponding tests should be not be expected
to have a correct level under `ARRes'.%
\[
\text{\textbf{[INSERT\ TABLE 3 HERE]}}%
\]
The performance of the $BP$ and $Parz$ tests is very good with levels that are
not oversized in general. However the $BP$ and $Parz$ tests can be undersized,
see the case of `ARCH(1)'. But even in this case the value of $\%\left\{
\widehat{p}^{\ast}\neq1\right\}  $ remains very small suggesting that the size
distortion is due to the critical values of Lee (2007).\footnote{This is
confirmed by a not reported simulation experiment which shows that using
standard chi-squared critical values give good results.} The behavior of the
$EL$ test is more erratic, with levels that can be either oversized, as in the
case of `GARCH(1,1)', `All Pass' and `ARRes', or undersized. The $IMSE$ test
can also be severely oversized. The $CvM$ behaves well for `GARCH(1,1)' and
`ARCH(1)' but, as expected, is severely size distorted in the other cases.

We now consider $\mathcal{H}_{1}$. In what follows, the critical values of the
$EL\,$and $IMSE$ tests are adjusted to achieve the desired level under
normality. A first set of fixed alternatives is considered, $MA1$:
$u_{t}=\varepsilon_{t}+0.05\varepsilon_{t-1}$, $AR1$: $u_{t}=0.05u_{t-1}%
+\varepsilon_{t}$, $MA4$: $\varepsilon_{t}+0.2\varepsilon_{t-4}$ and $AR6$:
$u_{t}=0.3u_{t-6}+\varepsilon_{t}$ with i.i.d. standard normal innovations
$\varepsilon_{t}$ and $n=200$, $1,000$ is considered. The $CvM$ test is
expected to perform better for these alternatives, especially `$AR1$' and
`$MA1$'. In Tables 4 and 5, $\overline{\widehat{p}^{\ast}}$ and $s_{\widehat
{p}^{\ast}}$ are the simulation mean and standard deviation of $\widehat
{p}^{\ast}$. These statistics are useful for assessing the impact of
$\overline{p}_{n}$ on the power since large $\overline{\widehat{p}^{\ast}}$ or
$s_{\widehat{p}^{\ast}}$ suggests that decreasing $\overline{p}_{n}$ can
decrease the power.
\[
\text{\textbf{[INSERT TABLE 4 HERE]}}%
\]
The low-lag `$AR1$' and `$MA1$' experiments have very similar characteristics
with powers of the tests for $\alpha=10\%$ increasing from $17\%-18\%$ for
$n=200$ to $43\%-47\%$ for $n=1,000$. The data-driven tests all exhibit a
surprisingly high $\overline{\widehat{p}^{\ast}}$ or $s_{\widehat{p}^{\ast}}$.
The $BP$, $Parz$ and $EL$ seem to be outperformed by the $IMSE$ and
$CvM\,$tests. For the higher-order experiments `$MA4$' and `$AR6$' and
$n=1,000$, the $BP$, $Parz$ and $EL$ tests clearly outperform their
competitors with power close or equal to $100\%$. For $n=200$, the $EL$ test
outperforms its competitors with $BP$ as a second-best. The high values of
$\overline{\widehat{p}^{\ast}}$ and $s_{\widehat{p}^{\ast}}$ for the $BP$ and
$Parz$ tests illustrate the fact that $\widehat{p}^{\ast}$ is suitable for
testing but not as an estimator of the order of an $AR$ or $MA$ process.

The second experiment under $\mathcal{H}_{1}$ examines, for $n=200$, the power
of the $5\%$ level $BP$ and $Parz$ tests against $H_{\rho}:u_{t}=v_{t}-\rho
v_{t-1}$, $\rho\in\left[  0,1/2\right]  $, under the nine scenarios of Tables
2 and 3. For example, under `GARCH(1,1)' $v_{t}=s_{t}\zeta_{t}$ and $s_{t}%
^{2}=0.001+0.90s_{t-1}^{2}+0.05v_{t-1}^{2}$ where $\zeta_{t}$ are i.i.d.
standard normal innovations while, under `ARRes', the $v_{t}$ are i.i.d.
$N\left(  0,1\right)  $ and $u_{t}=v_{t}-\rho v_{t-1}$ is estimated from the
$AR(1)$ model $X_{t}=0.8X_{t-1}+u_{t}$. We do not consider the other tests to
avoid undesirable size correction effects, but we compare $BP$ and $Parz$ with
$\widetilde{\mathcal{M}}_{n}^{EP32}$ test of Lee (2007) which rejects the null
when $n\widehat{R}_{1}^{2}\geq\widehat{z}\left(  \alpha\right)  $ where
$\widehat{z}\left(  \alpha\right)  $ is defined in (\ref{Test}), and an
$\alpha$ level test which rejects the null when $n\widehat{R}_{1}^{2}\geq
c\left(  \alpha\right)  $, where the infeasible $c\left(  \alpha\right)  $,
dependent of the white noise process under consideration, is computed from
$10,000$ preliminary replications. Since the latter is locally optimal under
Gaussianity, it is labelled $LOT$. Figure 1 reports the nine power graphs
corresponding to each white noise experiments.%
\[
\text{\textbf{[INSERT FIGURE 1 HERE]}}%
\]
Except for white noise processes such as `NoMDS' for which the new tests are
undersized, the power of the four tests are quite similar in the vicinity of
$\rho=0$, suggesting that our data-driven tests are, for processes close to
Gaussianity, not far from being locally optimal as $LOT$. The global
performance of all tests deteriorate for nonlinear white noise processes as
`ARCH(1)', for which $LOT$ has a very low power compared to its competitors
$BP$, $Parz$ and $\widetilde{\mathcal{M}}_{n}^{EP32}$. $Parz$ dominates its
competitors for such white noise processes. As expected from (\ref{Improv}),
$Parz$ and $BP$ perform as well as or better than $\widetilde{\mathcal{M}}%
_{n}^{EP32}$ which is less powerful than $Parz$ for heteroskedastic noises the
`Bilinear', `ARCH(1)', `GARCH(1,1)' or `NoMDS'.

The third experiment under $\mathcal{H}_{1}$ considers a second set of
alternatives given by randomized \textquotedblleft small
correlation\textquotedblright\ processes defined in (\ref{Wildalt}),
\begin{equation}
u_{t}=\varepsilon_{t}+\frac{\left(  2.5\times\gamma_{n}\right)  ^{1/2}%
}{n^{1/2}P^{1/4}}\sum_{k=1}^{P}\psi_{k,b}\varepsilon_{t-k}\text{,\quad
\quad\quad\quad}\psi_{k,b}\overset{\text{i.i.d.}}{\sim}N\left(  0,1\right)
\text{.} \label{Bayeswild}%
\end{equation}
In this setting $b$ is the simulation index, $b=1,...,10,000$. New moving
average coefficients $\left\{  \psi_{k,b}\right\}  $ are drawn for each
simulation. Randomizing the moving average coefficients allows us to explore
various shapes of the correlation function. The noise $\left\{  \varepsilon
_{t}\right\}  $ is independent of the moving average coefficients $\left\{
\psi_{k,b}\right\}  $ and is drawn randomly from the standard normal
distribution. Since $\sum_{k=1}^{P}\psi_{k,b}^{2}=P\left(  1+o_{\mathbb{P}%
}\left(  1\right)  \right)  $ when $P$ tends to infinity, the covariances of
(\ref{Bayeswild}) can be $o\left(  n^{-1/2}\right)  $ as shown in Lemma
\ref{Wildlem}. We consider two scenarios. In the experiment `LOW', $P$ is set
to $15$ for $n=200$ and to $75$ when $n=1,000$. The experiment `HIGH' doubles
the order $P$, so $P=30$ for $n=200$ and $P=150$ for $n=1,000$. The next table
reports simulation results.%
\[
\text{\textbf{[INSERT TABLE 5 HERE]}}%
\]
The $BP$ test outperforms its competitors and $Parz$ comes as a second-best.
The $EL$ test achieves power similar to that of the $BP$ test only in the LOW
experiment when $P=15$ and $n=200$. The power of the $IMSE$ and $CvM$ tests
decreases with the sample size while the power of the other tests increases,
showing the importance of a proper data-driven choice of the order. The high
values of $\overline{\widehat{p}^{\ast}}_{Parz}$ may suggest that the $Parz$
test would be negatively affected by choosing a lower value of $\overline
{p}_{n}$. However setting $\overline{p}_{n}=3\left[  \left(  n/2\right)
^{1/2}\right]  $ instead of $\overline{p}_{n}=n-1$ as done in an experiment
not reported does not really affect the power of the $BP$ test.

\section{Concluding remarks\label{Concluding remarks}}

\setcounter{equation}
{0}The paper proposes an automatic test for the weak white noise null
hypothesis for observed variables or residuals from a parametric model. The
test is based on a new data-driven order selection procedure applied to the
Box-Pierce (1970) test statistic. The critical region uses robust critical
values of Lee (2007) which can account for estimation of residuals. An
important theoretical finding is that the new test can detect alternatives
with small autocorrelation coefficients of order $\rho_{n}=o\left(
n^{-1/2}\right)  $ where $n$ is the sample size, provided that the number of
autocorrelation coefficients at moderate lags is large enough. The proposed
test is shown to be adaptive rate-optimal against this class of alternatives.
The paper gives examples of moving average alternatives with small
autocorrelation coefficients of order $o\left(  n^{-1/2}\right)  $ which are
detected by the new test but not by tests previously proposed by Deo (2000),
Escanciano and Lobato (2009) or Xiao and Wu (2011). These alternatives
correspond to a plausible macroeconomic scenario where a temporary shock has
no significant impact whereas permanent shocks may cause significant changes.
They can also be of interest in finance where arbitrage should rule out strong
deviations from the difference of martingale hypothesis, since these
alternatives generate conditional expectation given the past of these
alternatives with order $o_{\mathbb{P}}\left(  n^{-1/2}\right)  $. A
simulation experiment has shown that the new test can cope with various weak
types of white noise processes including the ARCH\ or GARCH processes popular
in empirical finance. The simulation experiment has also confirmed good power
properties of the test regarding detection of standard $AR(1)$ and $MA(1)$
alternatives when the noise is highly nonlinear, for instance in the case of
the $ARCH(1)$ process considered in the experiment.

\section{References}

\textsc{Anderson, T.W.} (1993). Goodness of Fit Tests for Spectral
Distributions. \textit{The Annals of Statistics} \textbf{21}, 830--847.

\textsc{Anderson, T.W.} and \textsc{D.A. Darling} (1952). Asymptotic Theory of
Certain \textquotedblleft Goodness of Fit\textquotedblright\ Criteria Based on
Stochastic Processes. \textit{Annals of Mathematical Statistics} \textbf{23}, 193--212.

\textsc{Box, G.} and \textsc{D. Pierce} (1970). Distribution of Residual
Autocorrelations in Autoregressive-Integrated Moving Average Time Series
Models. \textit{Journal of American Statistical Association} \textbf{65}, 1509--1526.

\textsc{Campbell, J.Y., A.W. Lo} and \textsc{A.C. Craig MacKinlay} (1997).
\textit{The Econometrics of Financial Markets}. Second Edition, Princeton
University Press.

\textsc{Chen, S.X.} and \textsc{J. Gao} (2007). An adaptive Empirical
Likelihood Test for Parametric Time Series Regression Models. \textit{Journal
of Econometrics }\textbf{141}, 950--972.

\textsc{Delgado, M.A.} and \textsc{C. Velasco} (2012). An Asymptotically
Pivotal Transform of the Residuals Sample Autocorrelations with Application to
Model Checking. \textit{Journal of American Statistical Association
}\textbf{106}, 646---958.

\textsc{Delgado, M.A., J. Hidalgo} and \textsc{C. Velasco} (2005).
Distribution Free Goodness-of-Fit Tests for Linear Processes. \textit{The
Annals of Statistics }\textbf{33}, 2568-2609.

\textsc{Deo, R.S.} (2000). Spectral Tests of the Martingale Hypothesis under
Conditional Heteroscedasticity. \textit{Journal of Econometrics} \textbf{99}, 291-315.

\textsc{Donoho, D.} and \textsc{J. Jin} (2004). Higher Criticism for Detecting
Sparse Heterogeneous Mixtures. \textit{The Annals of Statistics} \textbf{32}, 962--994.

\textsc{Durlauf, S.N.} (1991). Spectral Based Testing of the Martingale
Hypothesis. \textit{Journal of Econometrics} \textbf{50}, 355-376.

\textsc{Escanciano, J.C.} and \textsc{I.N. Lobato} (2009). An Automatic
Portmanteau Test for Serial Correlation. \textit{Journal of Econometrics
}\textbf{151}, 140--149.

\textsc{Fan, J.} (1996). Test of Significance Based on Wavelet Thresholding
and Neyman's Truncation. \textit{Journal of the American Statistical
Association} \textbf{91}, 674--688.

\textsc{Fan, J.} and \textsc{Q. Yao} (2005). \textit{Nonlinear Time Series:
Nonparametric and Parametric Methods}. Springer.

\textsc{Francq, C.} , \textsc{R. Roy} and \textsc{J.M. Zakoian} (2005).
Diagnostic Checking in ARMA Models With Uncorrelated Errors. \textit{Journal
of the American Statistical Association }\textbf{100}, 532--544.

\textsc{Golubev, G.K.}, \textsc{M. Nussbaum} and \textsc{H.H. Zhou} (2010).
Asymptotic Equivalence of Spectrum Density Estimation and Gaussian White
Noise. \textit{The Annals of Statistics }\textbf{38}, 181--214.

\textsc{Grenander, U.} and \textsc{M. Rosenblatt} (1952). On Spectral Analysis
of Stationary Time-series. \textit{Proceedings of the National Academy of
Sciences U.S.A. }\textbf{38, }519-521.

\textsc{Guay, A.} and \textsc{E. Guerre} (2006). A Data-Driven Nonparametric
Specification Test for Dynamic Regression Models. \textit{Econometric Theory}
\textbf{22}, 543--586.

\textsc{Guerre, E.} and \textsc{P. Lavergne} (2002). Optimal Minimax Rates for
Nonparametric Specification Testing in Regression Models. \textit{Econometric
Theory} \textbf{18}, 1139--1171.

\textsc{Guerre, E.} and \textsc{P. Lavergne} (2005). Rate-Optimal Data-Driven
Specification Testing for Regression Models. \textit{The Annals of Statistics}
\textbf{33}, 840--870.

\textsc{Hong, Y}. (1996). Consistent Testing for Serial Correlation of Unknown
Form. \textit{Econometrica} \textbf{64}, 837--864.

\textsc{Hong, Y}. and \textsc{Y.J. Lee}. (2005). Generalized Spectral Tests
for Conditional Mean Models in Time Series with Conditional Heteroscedasticity
of Unknown Form. \textit{Review of Economic Studies} \textbf{72}, 499--541.

\textsc{Horowitz, J.L.} and \textsc{V.G. Spokoiny} (2001). An Adaptive,
Rate-Optimal Test of a Parametric Mean-Regression Model Against a
Nonparametric Alternative. \textit{Econometrica} \textbf{69}, 599--631.

\textsc{Kuan, C.M. }and \textsc{W.M. Lee} (2006). Robust \textit{M} Tests
without Consistent Estimation of the Asymptotic Covariance Matrix.
\textit{Journal of the American Statistical Association }\textbf{101}, 1264--1275.

\textsc{Lee, W.M.} (2007). Robust M Tests Using Kernel-based Estimators with
Bandwidth Equal to Sample Size. \textit{Economics Letters }\textbf{96}, 295--300.

\textsc{Lobato, I.N.} (2001). Testing That a Dependent Process Is
Uncorrelated. \textit{Journal of the American Statistical Association
}\textbf{96}, 1066--1076.

\textsc{Lobato, I.N.}, \textsc{J.C. Nankervis} and \textsc{N.E. Savin} (2002).
Testing for Zero Autocorrelation in the Presence of Statistical Dependence.
\textit{Econometric Theory }\textbf{18}, 730--743.

\textsc{Newey, W.K. }and\textsc{\ K. West} (1994). Automatic Lag Selection in
Covariance Matrix Estimation. \textit{Review of Economic Studies} \textbf{61}, 631--653.

\textsc{Phillips, P.C.B, Y. Sun \& S. Jin} (2006). Spectral Density Estimation
and Robust Hypothesis Testing Using Steep Origin Kernels Without Truncation.
\textit{International Economic Review}\textbf{ 47}, 837--894.

\textsc{Pollard, D.} (2002). \textit{A User's Guide to Measure Theoretic
Probability}. Cambridge University Press.

\textsc{Shao, X.} (2011a). A Bootstrap-assisted Spectral Test of White Noise
under Unknown Dependence. \textit{Journal of Econometrics} \textbf{162}, 213--224.

\textsc{Shao, X.} (2011b). Testing for White Noise under Unknown Dependence
and its Applications to Goodness-of-Fit for Time Series Models.
\textit{Econometric Theory }\textbf{27}, 312--343.

\textsc{Spokoiny, V.G.} (1996). Adaptive Hypothesis Testing Using Wavelets.
\textit{The Annals of Statistics} \textbf{24}, 2477--2498.

\textsc{Xiao, H.} and \textsc{W.B. Wu} (2011). Asymptotic Inference of
Autocovariances of Stationary Processes. University of Chicago, arXiv:11053423v1.

\textsc{Wu, W.B.} (2005). Nonlinear System Theory: Another Look at Dependence.
\textit{Proceedings of the National Academy of Sciences of the United States
of America }\textbf{102}, 14150--14154.

\textsc{Wu, W.B.} (2007). Strong Invariance Principles for Dependent Random
Variables. \textit{The Annals of Probability }\textbf{35}, 2294--2320.

\pagebreak%

\newpage
\begin{landscape}%
\begin{figure} [h]\centering
\begin{center}
\includegraphics[
natheight=4.122600in,
natwidth=7.386400in,
height=4.1226in,
width=7.3864in
]%
{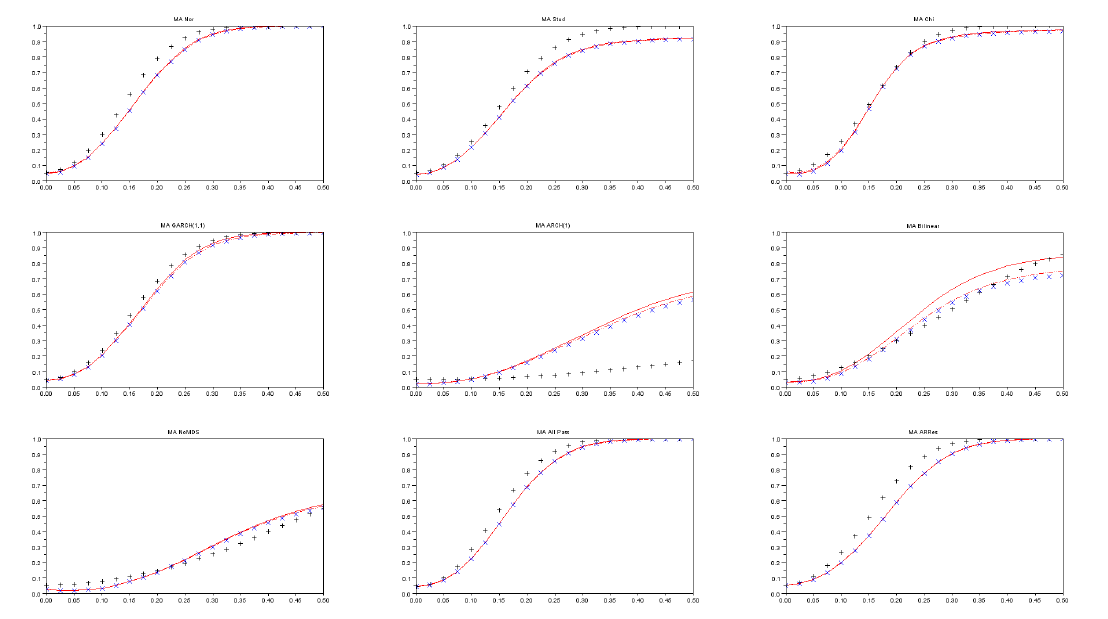}%
\\
{}%
\end{center}
\caption{\textbf{Empirical rejection probabilities} of $LOT$ (black `+' line), $Parz$ (red
solid line), $BP$ (red dotted line) and Lee (2007) $\protect\widetilde{\mathcal{M}}^{EP32}$
test (blue `x' line). The level of these tests is $5\%$ . The alternative is an
$MA(1)$ with a moving average coefficient ranging from $0$ to $1/2$ and
disturbances as in Tables 2-3. The sample size is $n=200$ and the number of
replications is $10,000$.
}%
\end{figure}%
\end{landscape}%

\pagebreak\includegraphics[width=1\textwidth]{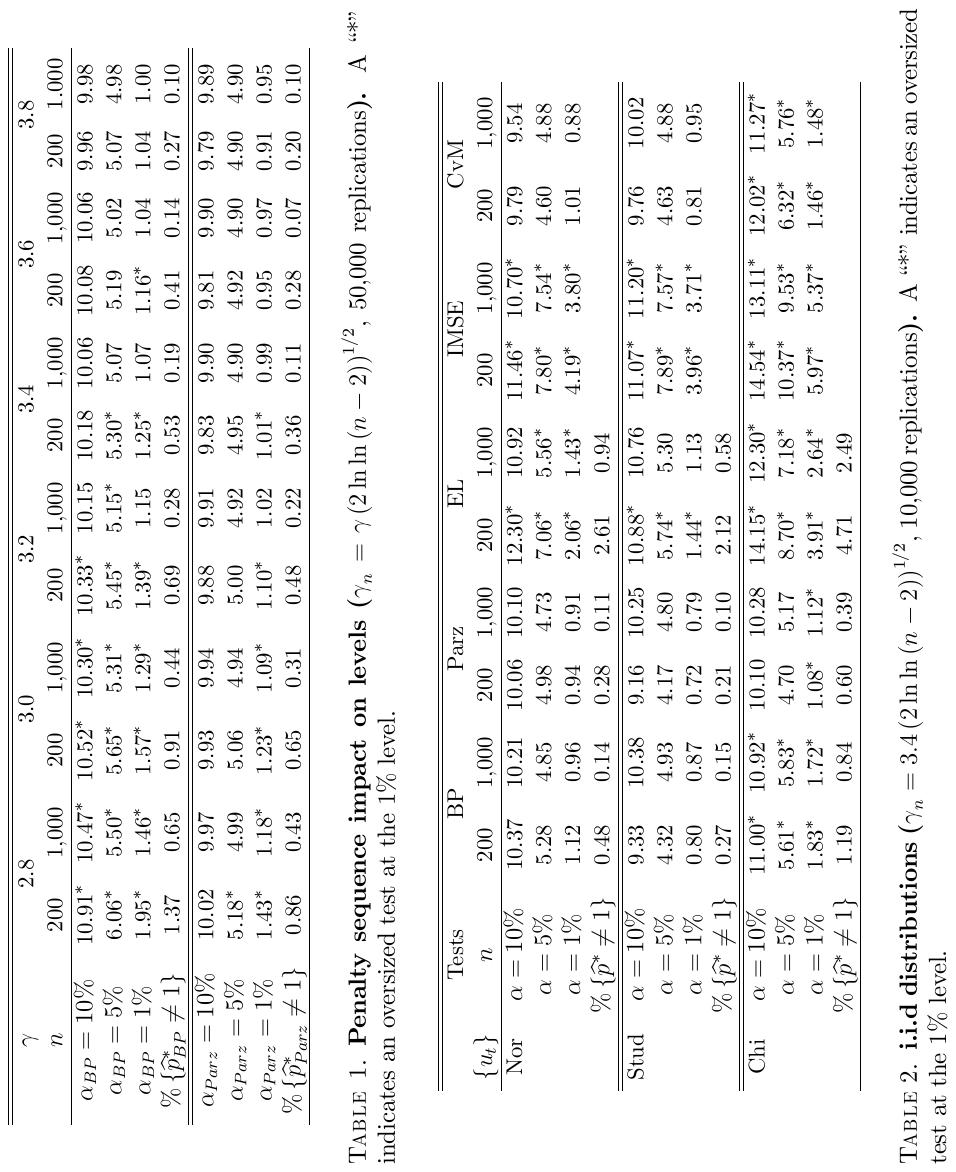}

\pagebreak

\includegraphics[width=1.00\textwidth]{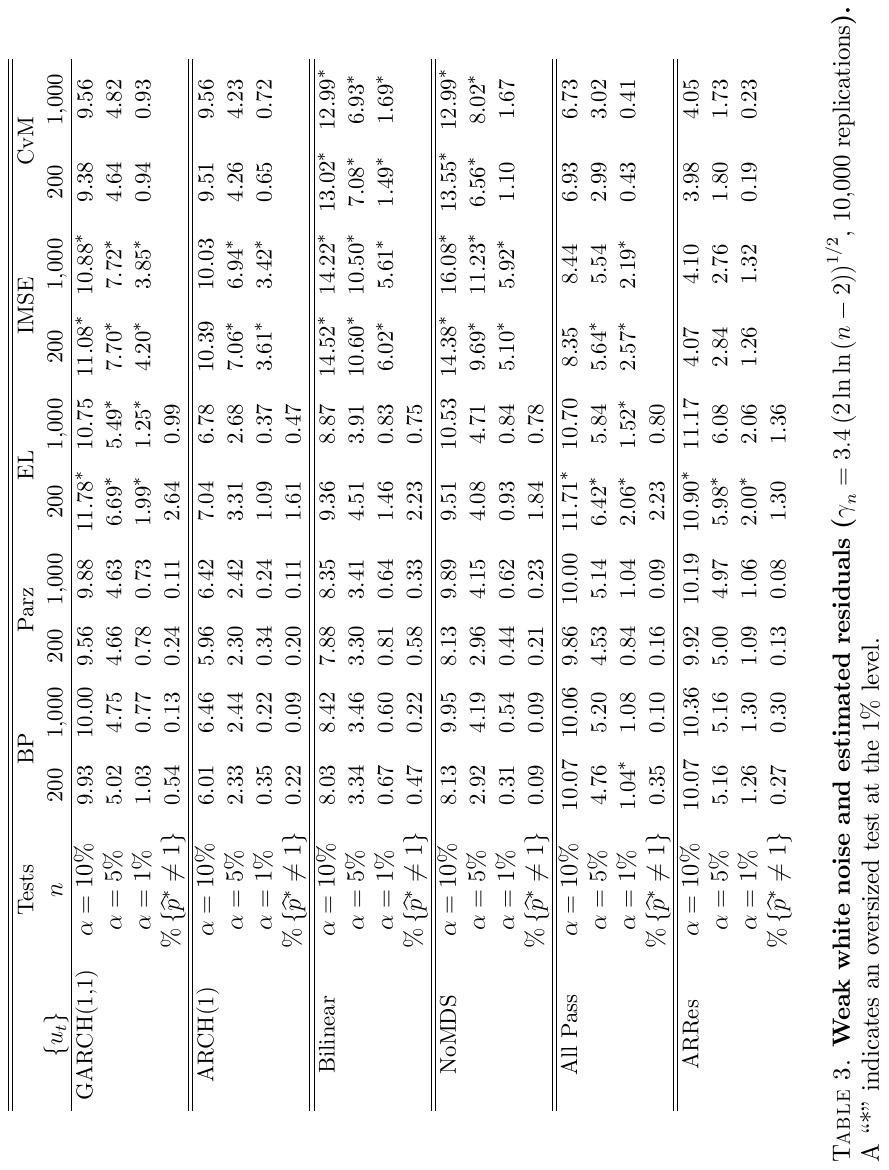}

\pagebreak

\includegraphics[width=1.00\textwidth]{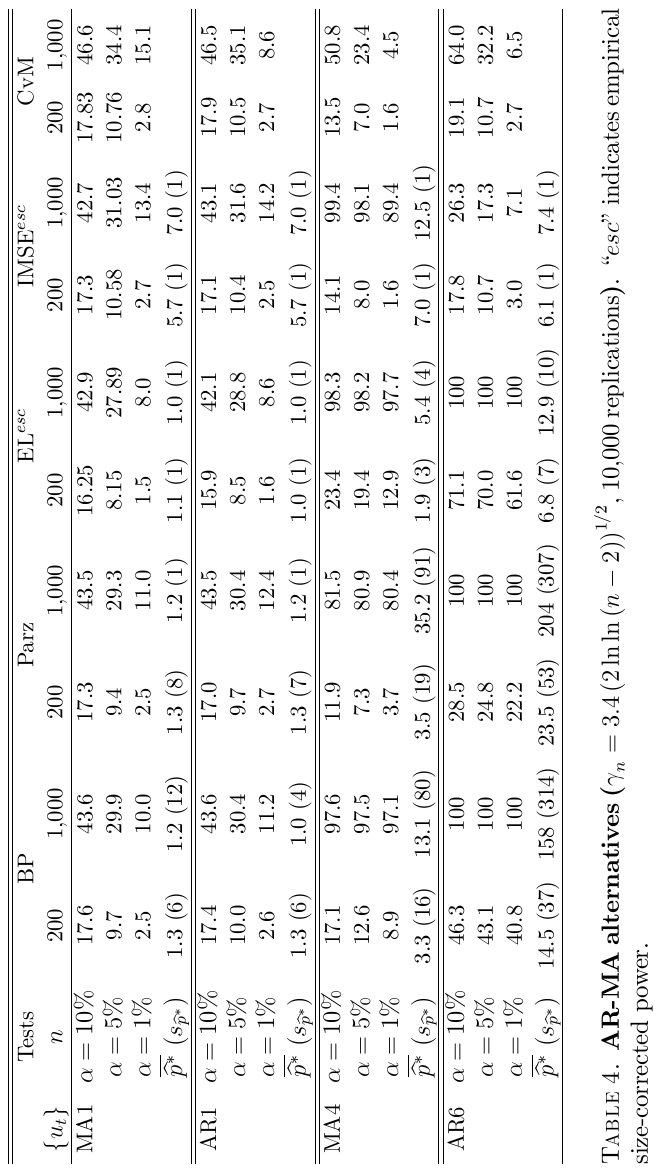}

\pagebreak

\includegraphics[width=1.00\textwidth]{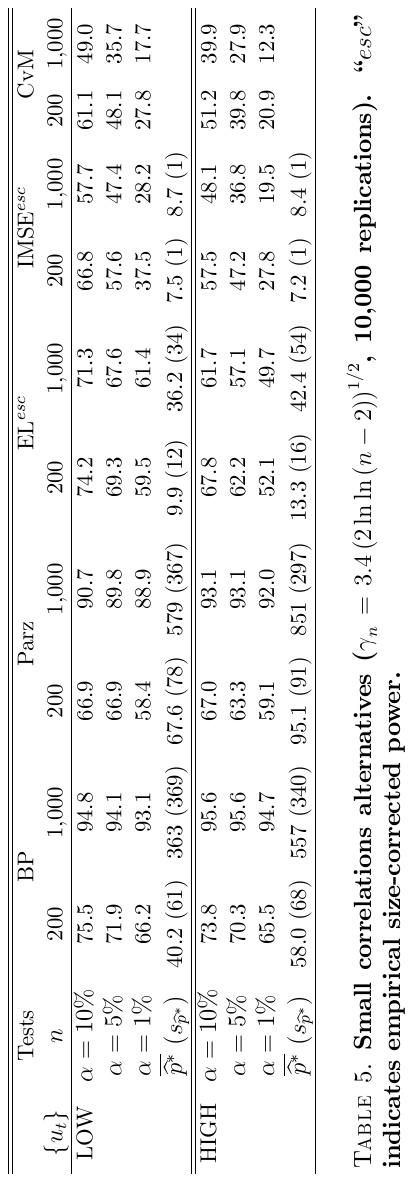}

\pagebreak

\newpage

\thispagestyle{empty} \vspace*{1cm}

\begin{center}
Robust Adaptive Rate-Optimal Testing for the White Noise Hypothesis:
Supplementary Material

\vspace*{0.5 cm}

Alain Guay\footnote{CIRP\'EE and CIREQ, Universit\'e du Qu\'ebec \`a
Montr\'eal, e-mail: \texttt{guay.alain@uqam.ca}} \\[0pt]

\medskip

Emmanuel Guerre\footnote{School of Economics and Finance, Queen Mary,
University of London, e-mail: \texttt{e.guerre@qmul.ac.uk }}\\[0pt]

\medskip

\v{S}t\v{e}p\'{a}na Lazarov\'{a}\footnote{School of Economics and Finance,
Queen Mary, University of London, e-mail: \texttt{s.lazarova@qmul.ac.uk }%
}\\[0pt]\vspace{0.5cm}

This version: 2nd November 2012

\vspace{1cm}
\end{center}

\renewcommand{\thefootnote}{\arabic{footnote}} \setcounter{page}{0} \newpage

\section*{Supplementary Material A: proofs of main results}

\renewcommand{\baselinestretch}{1.5} \setcounter{lem}{1} \setcounter{prop}{0}
\setcounter{thm}{0} \renewcommand{\thelem}{A.\arabic{lem}}
\renewcommand{\theprop}{A.\arabic{prop}}
\renewcommand{\thethm}{A.\arabic{thm}} \setcounter{equation}{0}
\setcounter{subsection}{0} \renewcommand{\theequation}{A.\arabic{equation}} \renewcommand{\thesubsection}{A.\arabic{subsection}}

This section contains the proofs of the results of Section 3. In what follows,
a tilde superscript, as in%
\begin{equation}
\widetilde{S}_{p}=n\sum_{j=1}^{p}K^{2}\left(  \frac{j}{p}\right)
\widetilde{R}_{j}^{2}\text{\quad where\quad}\widetilde{R}_{j}=\frac{1}{n}%
\sum_{t=1}^{n-|j|}u_{t}u_{t+|j|}. \label{TildeR}%
\end{equation}
indicates that the variables $u_{t}$ are observed. This also leads to define%
\[
\widetilde{\tau}_{j}=\frac{1}{n}\sum_{t=1}^{n-|j|}u_{t}^{2}u_{t+|j|}^{2}%
,\quad\widetilde{z}_{L}\left(  \alpha\right)  =\widehat{z}_{L}\left(
\alpha\right)  ,\quad\widetilde{z}_{L}^{\ast}\left(  \alpha\right)
=\widehat{z}_{L}^{\ast}\left(  \alpha\right)  ,
\]
but we keep the notation $\widehat{p}$. $C$ and $C^{\prime}$ are constants
that may vary from line to line but only depend on the constants of the
assumptions. Notation $\left[  \cdot\right]  $ is used for the integer part of
a real number and $a\vee b=\max\left(  a,b\right)  $, $a\wedge b=\min\left(
a,b\right)  $. Let $\overline{u}_{t}^{t-j}=\overline{u}_{t,n}^{t-j}$ be a copy
of $u_{t}=F_{n}\left(  \ldots,e_{t-1},e_{t}\right)  $ obtained by changing
$e_{t-j}$, $e_{t-j-1}$, $\ldots$ into $e_{t-j}^{\prime}$, $e_{t-j-1}^{\prime}%
$, $\ldots$. Then the condition $\left\Vert u_{t}-u_{t}^{t-j}\right\Vert
_{a}\leq\delta_{a}\left(  j\right)  $ ensures that%
\begin{equation}
\left\Vert u_{t}-\overline{u}_{t}^{t-j}\right\Vert _{a}\leq\Theta_{a}\left(
j\right)  \text{ where }\Theta_{a}\left(  j\right)  =\sum_{i=j}^{\infty}%
\delta_{a}\left(  j\right)  . \label{Thet}%
\end{equation}
We first state some intermediary results that are used in the proofs of our
main results. These intermediary results are proven in a section called
\textquotedblleft Supplementary Material B\textquotedblright. Lemma
\ref{Ordersums} gives the order of standardization terms $E(p)$, $E_{\Delta
}(p)$ and $V_{\Delta}(p)$. Propositions \ref{Covesti} and \ref{Esti} deal with
the impact of the estimation of $\theta$. Proposition \ref{Selection} is used
to study the asymptotic null behavior of the test and to show that
$\mathbb{P}\left(  \widehat{p}=1\right)  \rightarrow1$ in Theorem \ref{Level}.
Proposition \ref{Selection} deals with observed variables or residuals thanks
to Propositions \ref{Covesti} and \ref{Esti}. Propositions \ref{MeanH1} and
\ref{VarH1} are the key tools for our consistency result, Theorem
\ref{Sparse}. They dealt with observed variables but are combined with
Propositions \ref{Covesti} and \ref{Esti} to deal with estimation errors in
the proof of Theorem \ref{Sparse}.

\begin{lem}
\label{Ordersums}Suppose Assumption \ref{Kernel} holds and that $\overline
{p}_{n}/n\leq1/2$. (i) There exists a constant $C>1$ such that, for $q=1,2$
and for any $1\leq p\leq\overline{p}_{n}$, $\frac{p}{C}\leq\sum_{j=1}%
^{n-1}\left(  1-\frac{j}{n}\right)  ^{q}K^{2q}\left(  \frac{j}{p}\right)  \leq
Cp$, $\frac{p}{C}\leq\sum_{j=1}^{n-1}K^{2q}\left(  \frac{j}{p}\right)  \leq
Cp$, $V_{\Delta}^{2}(p)\leq Cp$, and $E_{\Delta}(p)\leq\sum_{j=1}^{n-1}\left(
K^{2}\left(  \frac{j}{p}\right)  -K^{2}\left(  j\right)  \right)  \leq
Cp^{1/2}V_{\Delta}(p)$; (ii) Under Assumption \ref{P}, for all $n$ and all
$p\in\left[  1,\overline{p}_{n}\right]  $, $V_{\Delta}(p)\geq C(p-1)^{1/2}$
and $E_{\Delta}(p)\geq0$.
\end{lem}

\begin{lem}
\label{L} Suppose Assumptions \ref{Kernel}, \ref{M} and \ref{Reg} hold. Then
the rejection regions $\widetilde{S}_{1}\geq\widetilde{z}_{L}\left(
\alpha\right)  $, $\widetilde{S}_{1}^{\ast}\geq\widetilde{z}_{L}^{\ast}\left(
\alpha\right)  $, $\widehat{S}_{1}\geq\widehat{z}_{KL}\left(  \alpha\right)  $
and $\widehat{S}_{1}^{\ast}\geq\widehat{z}_{KL}^{\ast}\left(  \alpha\right)  $
are asymptotically of level $\alpha$. Moreover, under $\mathcal{H}_{1}$,
$\widehat{z}_{L}\left(  \alpha\right)  $, $\widetilde{z}_{L}^{\ast}\left(
\alpha\right)  $, $\widehat{z}_{KL}\left(  \alpha\right)  $ and $\widehat
{z}_{KL}^{\ast}\left(  \alpha\right)  $ are all $O_{\mathbb{P}}\left(
1\right)  $.
\end{lem}

\begin{lem}
Under Assumption \ref{Reg}, $\sup_{0\leq j\leq n-1}\operatorname*{Var}\left(
\widetilde{R}_{j}\right)  \leq\frac{C}{n}$. \label{Varcov}
\end{lem}

\begin{prop}
\label{Covesti}Suppose Assumptions \ref{M}, \ref{P} and \ref{Reg} hold. Then
$\max_{j\in\left[  0,\overline{p}_{n}\right]  }\left\vert \widehat{R}%
_{j}-\widetilde{R}_{j}\right\vert =O_{\mathbb{P}}\left(  n^{-1/2}\right)  $,
$\max_{p\in\left[  0,n-1\right]  }n\sum_{j=1}^{p}\left(  \widehat{R}%
_{j}-\widetilde{R}_{j}\right)  ^{2}=O_{\mathbb{P}}\left(  1\right)  $, and%
\begin{align*}
\max_{j\in\left[  0,n-1\right]  }\left\vert \widetilde{R}_{j}-\left(
1-\frac{j}{n}\right)  R_{j,n}\right\vert  &  =O_{\mathbb{P}}\left(  \left(
\frac{\log n}{n}\right)  ^{1/2}\right)  ,\\
\max_{j\in\left[  0,\overline{p}_{n}\right]  }\left\vert \widehat{R}%
_{j}-R_{j,n}\right\vert  &  =O_{\mathbb{P}}\left(  \left(  \frac{\log n}%
{n}\right)  ^{1/2}\right)  ,\\
\max_{j\in\left[  0,n-1\right]  }\left(  1-\frac{j}{n}\right)  \left\vert
\widetilde{\tau}_{j}^{2}-\tau_{j,n}^{2}\right\vert  &  =O_{\mathbb{P}}\left(
\left(  \frac{\log n}{n}\right)  ^{1/2}\right)  ,\\
\max_{j\in\left[  0,\overline{p}_{n}\right]  }\left\vert \widehat{\tau}%
_{j}^{2}-\tau_{j,n}^{2}\right\vert  &  =O_{\mathbb{P}}\left(  \left(
\frac{\log n}{n}\right)  ^{1/2}\right)  .
\end{align*}

\end{prop}

\begin{prop}
\label{Esti}Let Assumptions \ref{Kernel}, \ref{M}, \ref{P} and \ref{Reg} hold.
Let $\widetilde{S}_{p}$ be as in (\ref{TildeR}). Then
\[
\max_{p\in\left[  2,\overline{p}_{n}\right]  }\frac{|\left(  \widehat{S}%
_{p}-\widehat{S}_{1}\right)  -\left(  \widetilde{S}_{p}-\widetilde{S}%
_{1}\right)  |}{1+\left(  n\sum_{j=1}^{p}R_{j,n}^{2}\right)  ^{1/2}%
}=O_{\mathbb{P}}\left(  1\right)
\]
and for any $p_{n}=O(n^{1/2})$, $\widehat{S}_{p_{n}}-\widetilde{S}_{p_{n}%
}=O_{\mathbb{P}}\left(  1+\left(  n\sum_{j=1}^{p_{n}}R_{j,n}^{2}\right)
^{1/2}\right)  .$
\end{prop}

\begin{prop}
\label{Selection}Suppose Assumptions \ref{Kernel}, \ref{M}, \ref{P} and
\ref{Reg} hold and that $\mathcal{H}_{0}$ is true. Then (\ref{Gam}) ensures
that
\[
\lim_{n\rightarrow\infty}\mathbb{P}\left(  \max_{p\in\left[  2,\overline
{p}_{n}\right]  }\frac{(\widehat{S}_{p}-\widehat{S}_{1})/\widehat{R}_{0}%
^{2}-E_{\Delta}(p)}{V_{\Delta}(p)}\geq\gamma_{n}\right)  =0.
\]

\end{prop}

\begin{prop}
Under Assumptions \ref{Kernel}, \ref{P} and \ref{Reg}, there are some
$C,C^{\prime}>0$ such that for $n$ large enough and uniformly in $p\in\left[
1,\overline{p}_{n}\right]  $,
\begin{align*}
\mathbb{E}\left[  \widetilde{S}_{p}\right]  -R_{0,n}^{2}E\left(  p\right)   &
\geq Cn\sum_{j=1}^{p/2}R_{j,n}^{2}-C^{\prime}R_{0,n}^{2},\\
\mathbb{E}\left[  \sum_{j=1}^{n-1}K\left(  \frac{j}{p}\right)  \frac
{\widetilde{R}_{j}^{2}}{\tau_{j,n}^{2}}\right]  -E\left(  p\right)   &  \geq
Cn\sum_{j=1}^{p/2}\left(  \frac{R_{j,n}}{R_{0,n}}\right)  ^{2}-C^{\prime}.
\end{align*}
\label{MeanH1}
\end{prop}

\begin{prop}
Under Assumptions \ref{Kernel}, \ref{P} and \ref{Reg}, there is a constant
$C>0$ such that for $n$ large enough and uniformly in $p\in\left[
1,\overline{p}_{n}\right]  $,%
\begin{align*}
\operatorname*{Var}\left(  \widetilde{S}_{p}\right)   &  \leq C\left(
n\sum_{j=1}^{p}R_{j,n}^{2}+p\right)  ,\\
\operatorname*{Var}\left(  \sum_{j=1}^{n-1}K\left(  \frac{j}{p}\right)
\frac{\widetilde{R}_{j}^{2}}{\tau_{j,n}^{2}}\right)   &  \leq C\left(
n\sum_{j=1}^{p}\frac{R_{j,n}^{2}}{R_{0,n}^{2}}+p\right)  .
\end{align*}
\label{VarH1}
\end{prop}

\subsection{Proof of Theorem \ref{Level}}

(\ref{Hatpnot}), (\ref{Gam}) and Proposition \ref{Selection} give that
$\lim_{n\rightarrow\infty}\mathbb{P}(\widehat{p}\neq1)=0$. Hence $\widehat
{S}_{\widehat{p}}=\widehat{S}_{1}+o_{\mathbb{P}}\left(  1\right)  $ and Lemma
\ref{L}, which ensures that the retained critical value satisfies
$\mathbb{P}\left(  \widehat{S}_{1}\geq\widehat{z}\left(  \alpha\right)
\right)  \rightarrow\alpha$, yield that the test (\ref{Test}) is
asymptotically of level $\alpha$.\hspace*{\fill}$\Box$

\subsection{Proof of Theorem \ref{Sparse}}

The definition (\ref{Hatp}) of $\widehat{p}$ gives, for any $p\in\left[
1,\overline{p}_{n}\right]  $,%
\begin{align*}
\widehat{S}_{\widehat{p}}  &  =\arg\max_{p\in\left[  1,\overline{p}%
_{n}\right]  }\left\{  \widehat{S}_{p}-\widehat{R}_{0}^{2}E\left(  p\right)
-\gamma_{n}\widehat{R}_{0}^{2}V_{\Delta}\left(  p\right)  \right\}
+\widehat{R}_{0}^{2}E\left(  \widehat{p}\right)  +\gamma_{n}\widehat{R}%
_{0}^{2}V_{\Delta}\left(  \widehat{p}\right) \\
&  \geq\widehat{S}_{p}-\widehat{R}_{0}^{2}E\left(  p\right)  -\gamma
_{n}\widehat{R}_{0}^{2}V_{\Delta}\left(  p\right)  .
\end{align*}
Note that this bound implies (\ref{Detect}). Since the critical value
$\widehat{z}\left(  \alpha\right)  $ in (\ref{Test}) is bounded under
$\mathcal{H}_{1}$ by Lemma \ref{L}, it is sufficient to find a $p_{n}%
\in\left[  1,\overline{p}_{n}\right]  $ such that $\widehat{S}_{p_{n}%
}-\widehat{R}_{0}^{2}E\left(  p_{n}\right)  -\gamma_{n}\widehat{R}_{0}%
^{2}V_{\Delta}\left(  p_{n}\right)  \overset{\mathbb{P}}{\rightarrow}+\infty$.
Let $p_{n}=2P_{n}$ where $P_{n}$ is as in (\ref{Sparse3}). Set
\[
\mathcal{R}_{n}^{2}=\sum_{j=1}^{P_{n}}\left(  \frac{R_{j,n}}{R_{0,n}}\right)
^{2}.
\]
The detection condition (\ref{Sparse3}) gives%
\begin{equation}
n\mathcal{R}_{n}^{2}\geq n\rho_{n}^{2}\sum_{j=1}^{P_{n}}\mathbb{I}\left\{
\left(  \frac{R_{j,n}}{R_{0,n}}\right)  ^{2}\geq\rho_{n}^{2}\right\}
=nN_{n}\rho_{n}^{2}\geq\frac{\kappa_{\ast}^{2}\gamma_{n}p_{n}^{1/2}}{2^{1/2}%
}\rightarrow\infty, \label{Sparse4}%
\end{equation}
with a constant $\kappa_{\ast}$ which can be chosen as large as needed. Lemmas
\ref{Ordersums}, \ref{Varcov}, Assumption \ref{P} which ensures $P_{n}%
=o\left(  n^{1/2}\right)  $ and $\gamma_{n}=o\left(  n^{1/4}\right)  $, and
Proposition \ref{Covesti} for the case of residuals yield that%
\begin{align*}
&  \widehat{S}_{p_{n}}-\widehat{R}_{0}^{2}E\left(  p_{n}\right)  -\gamma
_{n}\widehat{R}_{0}^{2}V_{\Delta}\left(  p_{n}\right) \\
&  =\widetilde{S}_{p_{n}}+O_{\mathbb{P}}\left(  1+n^{1/2}R_{0,n}%
\mathcal{R}_{n}\right)  -R_{0,n}^{2}E\left(  p_{n}\right)  -\gamma_{n}%
R_{0,n}^{2}V_{\Delta}\left(  p_{n}\right)  +O_{\mathbb{P}}\left(  \frac
{p_{n}+\gamma_{n}p_{n}^{1/2}}{n^{1/2}}\right) \\
&  \geq\widetilde{S}_{p_{n}}+O_{\mathbb{P}}\left(  1+n^{1/2}R_{0,n}%
\mathcal{R}_{n}\right)  -R_{0,n}^{2}E\left(  p_{n}\right)  -C\gamma_{n}%
R_{0,n}^{2}p_{n}^{1/2}.
\end{align*}
Now the Chebyshev inequality, Propositions \ref{MeanH1} and \ref{VarH1}, give%
\[
\widetilde{S}_{p_{n}}=\mathbb{E}\left[  \widetilde{S}_{p_{n}}\right]
+O_{\mathbb{P}}\left(  \operatorname*{Var}\nolimits^{1/2}\left(  \widetilde
{S}_{p_{n}}\right)  \right)  \geq R_{0,n}^{2}E\left(  p_{n}\right)
+C^{\prime}R_{0,n}^{2}n\mathcal{R}_{n}^{2}+O_{\mathbb{P}}\left(  p_{n}%
^{1/2}+n^{1/2}\mathcal{R}_{n}\right)  .
\]
Hence substituting gives, since $n\mathcal{R}_{n}^{2}\rightarrow\infty$ by
(\ref{Sparse4}),%
\[
\widehat{S}_{p_{n}}-\widehat{R}_{0}^{2}E\left(  p_{n}\right)  -\gamma
_{n}\widehat{R}_{0}^{2}V_{\Delta}\left(  p_{n}\right)  \geq C^{\prime}%
R_{0,n}^{2}n\mathcal{R}_{n}^{2}\left(  1+o_{\mathbb{P}}\left(  1\right)
\right)  -C\gamma_{n}R_{0,n}^{2}p_{n}^{1/2}\left(  1+o_{\mathbb{P}}\left(
1\right)  \right)  .
\]
Since Assumption \ref{Reg} ensures that $R_{0,n}^{2}$ stays bounded away from
$0$, (\ref{Sparse4}) gives that $\widehat{S}_{p_{n}}-\widehat{R}_{0}%
^{2}E\left(  p_{n}\right)  -\gamma_{n}\widehat{R}_{0}^{2}V_{\Delta}\left(
p_{n}\right)  \overset{\mathbb{P}}{\rightarrow}+\infty$ as requested provided
$\kappa_{\ast}^{2}>C^{\prime}/C$. $\hfill\square$

\subsection{Proof of Theorem \ref{Extension}}

Consider first the null hypothesis. As seen from the proof of Theorem
\ref{Level}, it suffices to show that%
\[
\lim_{n\rightarrow\infty}\mathbb{P}\left(  \max_{p\in\left[  2,\overline
{p}_{n}\right]  }\frac{(\widehat{S}_{p}^{\ast}-\widehat{S}_{1}^{\ast
})-E_{\Delta}(p)}{V_{\Delta}(p)}\geq\gamma_{n}\right)  =0,
\]
a statement which implies that $\widehat{p}^{\ast}=1+o_{\mathbb{P}}\left(
1\right)  $ so that Lemma \ref{L} implies that the conclusion of Theorem
\ref{Level} holds for the test based upon $\widehat{S}_{\widehat{p}^{\ast}%
}^{\ast}$. Since $\left\vert R_{j,n}\right\vert \leq\left\Vert u_{t,n}%
\right\Vert _{2}\left\Vert u_{t,n}-\overline{u}_{t,n}^{t-j}\right\Vert _{2}$
and%
\begin{align*}
\mathbb{E}\left[  u_{t-j,n}^{2}u_{t-j,n}^{2}\right]   &  =\mathbb{E}\left[
\left(  \overline{u}_{t,n}^{t-j}\right)  ^{2}u_{t-j,n}^{2}\right]
+\mathbb{E}\left[  \left(  u_{t,n}^{2}-\left(  \overline{u}_{t,n}%
^{t-j}\right)  ^{2}\right)  u_{t-j,n}^{2}\right] \\
&  =R_{0,n}^{2}+\mathbb{E}\left[  \left(  u_{t,n}-\overline{u}_{t,n}%
^{t-j}\right)  \left(  u_{t,n}+\overline{u}_{t,n}^{t-j}\right)  u_{t-j,n}%
^{2}\right]  ,
\end{align*}
(\ref{Thet})\ shows
\begin{equation}
\left\vert \tau_{j,n}^{2}-R_{0,n}^{2}\right\vert \leq C\left\Vert
u_{t,n}\right\Vert _{8}^{3}\Theta_{2}\left(  j\right)  \leq Cj^{-6}
\label{Tau2sig}%
\end{equation}
for all $j\geq1$. Now Lemmas \ref{Ordersums} and \ref{Varcov}, Assumptions
\ref{Kernel}, \ref{P} and \ref{Reg}, and Proposition \ref{Covesti} give%
\begin{align*}
&  \max_{p\in\left[  2,\overline{p}_{n}\right]  }\frac{\left\vert (\widehat
{S}_{p}^{\ast}-\widehat{S}_{1}^{\ast})-(\widehat{S}_{p}-\widehat{S}%
_{1})/\widehat{R}_{0}^{2}\right\vert }{V_{\Delta}(p)}\leq C\max_{p\in\left[
1,\overline{p}_{n}\right]  }\frac{\left\vert \widehat{S}_{p}^{\ast}%
-\widehat{S}_{p}/\widehat{R}_{0}^{2}\right\vert }{p^{1/2}}\\
&  \leq C\max_{p\in\left[  1,\overline{p}_{n}\right]  }\frac{n}{p^{1/2}}%
\sum_{j=1}^{p}\left(  \frac{\widehat{R}_{j}}{\widehat{R}_{0}}\right)
^{2}\left\{  \left\vert \frac{\widehat{\tau}_{j}^{2}}{\widehat{R}_{0}^{2}%
}-\frac{\tau_{j,n}^{2}}{R_{0,n}^{2}}\right\vert +\left\vert \frac{\tau
_{j,n}^{2}}{R_{0,n}^{2}}-1\right\vert \right\} \\
&  \leq Cn\overline{p}_{n}^{1/2}O_{\mathbb{P}}\left(  \left(  \frac{\log n}%
{n}\right)  ^{3/2}\right)  +O_{\mathbb{P}}\left(  1\right)  n\sum
_{j=1}^{\overline{p}_{n}}\frac{\widehat{R}_{j}^{2}}{j^{6}}\\
&  =o_{\mathbb{P}}\left(  1\right)  +O_{\mathbb{P}}\left(  \sum_{j=1}%
^{\overline{p}_{n}}\frac{\operatorname*{Var}\left(  n^{1/2}\widehat{R}%
_{j}\right)  }{j^{6}}\right)  =O_{\mathbb{P}}\left(  1\right)  .
\end{align*}
Hence (\ref{Gam}) and Proposition \ref{Selection}%
\begin{align*}
&  \mathbb{P}\left(  \max_{p\in\left[  2,\overline{p}_{n}\right]  }%
\frac{(\widehat{S}_{p}^{\ast}-\widehat{S}_{1}^{\ast})-E_{\Delta}(p)}%
{V_{\Delta}(p)}\geq\gamma_{n}\right) \\
&  \text{ }=\mathbb{P}\left(  \max_{p\in\left[  2,\overline{p}_{n}\right]
}\frac{(\widehat{S}_{p}-\widehat{S}_{1})/\widehat{R}_{0}^{2}-E_{\Delta}%
(p)}{V_{\Delta}(p)}+O_{\mathbb{P}}\left(  1\right)  \geq\gamma_{n}\right) \\
&  \text{ }\leq\mathbb{P}\left(  \max_{p\in\left[  2,\overline{p}_{n}\right]
}\frac{(\widehat{S}_{p}-\widehat{S}_{1})/\widehat{R}_{0}^{2}-E_{\Delta}%
(p)}{V_{\Delta}(p)}\geq\left(  1+\frac{\epsilon}{2}\right)  \left(  2\ln\ln
n\right)  ^{1/2}\right)  +o\left(  1\right) \\
&  \text{ }\mathbb{=}o\left(  1\right)  ,
\end{align*}
which gives the desired result under $\mathcal{H}_{0}$.

Consider now Theorem \ref{Sparse} and $\mathcal{H}_{1}$. Define%
\[
\widehat{S}_{p}^{\bigstar}=n\sum_{j=1}^{p}K^{2}\left(  \frac{j}{p}\right)
\frac{\widehat{R}_{j}^{2}}{\tau_{j,n}^{2}},\quad\widetilde{S}_{p}^{\bigstar
}=n\sum_{j=1}^{p}K^{2}\left(  \frac{j}{p}\right)  \frac{\widetilde{R}_{j}^{2}%
}{\tau_{j,n}^{2}}.
\]
Let $P_{n}$ be as in (\ref{Sparse3}) and define $p_{n}=2P_{n}$ and
$\mathcal{R}_{n}$ as in the proof of Theorem \ref{Sparse}. Then Assumptions
\ref{Kernel} and \ref{Reg}, Propositions \ref{Covesti} and \ref{Esti}
\begin{align*}
\left\vert \widehat{S}_{p_{n}}^{\ast}-\widehat{S}_{p_{n}}^{\bigstar
}\right\vert  &  \leq Cn\sum_{j=1}^{p_{n}}\frac{\widehat{R}_{j}^{2}}%
{\tau_{j,n}^{2}}\left\vert \frac{\tau_{j,n}^{2}}{\widehat{\tau}_{j}^{2}%
}-1\right\vert =O_{\mathbb{P}}\left(  \left(  \frac{\log n}{n}\right)
^{1/2}\right)  \check{S}_{p_{n}}^{\bigstar},\\
\left\vert \widehat{S}_{p_{n}}^{\bigstar}-\widetilde{S}_{p_{n}}^{\bigstar
}\right\vert  &  \leq C\left\vert \widehat{S}_{p_{n}}-\widetilde{S}_{p_{n}%
}\right\vert =O_{\mathbb{P}}\left(  n^{1/2}\mathcal{R}_{n}\right)  .
\end{align*}
Hence, for observed variables or residuals,%
\[
\widehat{S}_{p_{n}}^{\ast}=\left(  1+O_{\mathbb{P}}\left(  \left(  \frac{\log
n}{n}\right)  ^{1/2}\right)  \right)  \widetilde{S}_{p_{n}}^{\bigstar
}+O_{\mathbb{P}}\left(  n^{1/2}\mathcal{R}_{n}\right)
\]
The proof now follows the steps of the one of Theorem \ref{Sparse} based on
the order above, Proposition \ref{MeanH1} and \ref{VarH1}, and Lemma
\ref{Varcov} which gives $\mathbb{E}\left[  \widetilde{S}_{p_{n}}^{\bigstar
}\right]  \leq C\left(  p_{n}+n\mathcal{R}_{n}^{2}\right)  $. Hence, since
$p_{n}=o\left(  \left(  \log n/n\right)  ^{1/2}\right)  $,%
\begin{align*}
\widehat{S}_{\widehat{p}^{\ast}}^{\ast}  &  =\arg\max_{p\in\left[
1,\overline{p}_{n}\right]  }\left\{  \widehat{S}_{p}^{\ast}-E\left(  p\right)
-\gamma_{n}V_{\Delta}\left(  p\right)  \right\}  +E\left(  \widehat{p}^{\ast
}\right)  +\gamma_{n}V_{\Delta}\left(  \widehat{p}^{\ast}\right) \\
&  \geq\widehat{S}_{p_{n}}^{\ast}-E\left(  p_{n}\right)  -C\gamma_{n}%
p_{n}^{1/2}\\
&  =\left(  1+O_{\mathbb{P}}\left(  \left(  \frac{\log n}{n}\right)
^{1/2}\right)  \right)  \left(  \mathbb{E}\left[  \widetilde{S}_{p_{n}%
}^{\bigstar}\right]  +\operatorname*{Var}\nolimits^{1/2}\left(  \widetilde
{S}_{p_{n}}^{\bigstar}\right)  \right)  -E\left(  p_{n}\right)  -C\gamma
_{n}p_{n}^{1/2}\\
&  =C^{\prime}R_{0,n}^{2}n\mathcal{R}_{n}^{2}-C\gamma_{n}R_{0,n}^{2}%
p_{n}^{1/2}+O_{\mathbb{P}}\left(  p_{n}^{1/2}+n^{1/2}\mathcal{R}_{n}+\left(
\frac{\log n}{n}\right)  ^{1/2}\left(  p_{n}+n\mathcal{R}_{n}^{2}\right)
\right) \\
&  =C^{\prime}R_{0,n}^{2}n\mathcal{R}_{n}^{2}\left(  1+o_{\mathbb{P}}\left(
1\right)  \right)  -C\gamma_{n}R_{0,n}^{2}p_{n}^{1/2}\left(  1+o_{\mathbb{P}%
}\left(  1\right)  \right)  \overset{\mathbb{P}}{\rightarrow}+\infty
\end{align*}
provided $\kappa_{\ast}$ is large enough.$\hfill\square$

\subsection{Proof of Theorem \ref{AltCV}}

Since $\mathbb{P}\left(  \widehat{p}^{\ast}=1\right)  \rightarrow1$ under
$\mathcal{H}_{0}$, condition (A0) and (\ref{CVc}) give%
\begin{align*}
\lim_{n\rightarrow\infty}\mathbb{P}\left(  \widehat{S}_{\widehat{p}^{\ast}%
}^{\ast}\geq\widehat{c}_{n}^{\ast}\left(  \alpha\right)  \right)   &
=\lim_{n\rightarrow\infty}\mathbb{P}\left(  \widehat{S}_{1}^{\ast}\geq
\widehat{c}_{n}^{\ast}\left(  \alpha\right)  \right)  =\lim_{n\rightarrow
\infty}\mathbb{P}\left(  \widehat{S}_{1}^{\ast}\geq\widehat{S}_{1}^{\ast
}-\widehat{T}_{n}+\widehat{t}_{n}\left(  \alpha\right)  \right) \\
&  =\lim_{n\rightarrow\infty}\mathbb{P}\left(  \widehat{T}_{n}\geq\widehat
{t}_{n}\left(  \alpha\right)  \right)  =\alpha,
\end{align*}
so that the test of interest is asymptotically of level $\alpha$. Let us now
consider the alternative. Arguing as in the proof of Theorems \ref{Sparse} and
\ref{Extension} under condition (A1) shows that the test with critical value
$\widehat{c}_{n}\left(  \alpha\right)  $ detects the alternatives
(\ref{Sparse3}) provided $\kappa_{\ast}$ is taken large enough. Consider now
(\ref{Improv}). The definition of (\ref{Hatpstar}) gives, since $E_{\Delta
}\left(  \widehat{p}^{\ast}\right)  +\gamma_{n}V_{\Delta}(\widehat{p}^{\ast
})\geq0$ under Assumption \ref{Kernel},
\begin{align*}
\widehat{S}_{\widehat{p}^{\ast}}^{\ast}  &  =\max_{p\in\left[  1,\overline
{p}_{n}\right]  }\left(  \widehat{S}_{p}^{\ast}-E_{\Delta}\left(  p\right)
-\gamma_{n}V_{\Delta}(p)\right)  +E_{\Delta}\left(  \widehat{p}^{\ast}\right)
+\gamma_{n}V_{\Delta}(\widehat{p}^{\ast})\\
&  \geq\widehat{S}_{1}^{\ast}-E_{\Delta}\left(  1\right)  -\gamma_{n}%
V_{\Delta}(1)=\widehat{S}_{1}^{\ast}.
\end{align*}
Hence, by (\ref{CVc})%
\[
\mathbb{P}\left(  \widehat{S}_{\widehat{p}^{\ast}}^{\ast}\geq\widehat{c}%
_{n}\left(  \alpha\right)  \right)  \geq\mathbb{P}\left(  \widehat{S}%
_{1}^{\ast}\geq\widehat{c}_{n}\left(  \alpha\right)  \right)  =\mathbb{P}%
\left(  \widehat{S}_{1}^{\ast}\geq\widehat{S}_{1}^{\ast}-\widehat{T}%
_{n}+\widehat{t}_{n}\left(  \alpha\right)  \right)  =\mathbb{P}\left(
\widehat{T}_{n}\geq\widehat{t}_{n}\left(  \alpha\right)  \right)  ,
\]
which is (\ref{Improv}).$\hfill\square$

\subsection{Proof of Theorem \ref{Sparseopt}}

We first introduce a set of alternatives. Let $f\left(  \cdot\right)  $ denote
the spectral density of a centered Gaussian stationary process $\left\{
u_{t}\right\}  .$with covariance coefficients $R_{j}$. Define a H\"{o}lder
class of \ processes as
\[
\text{H\"{o}lder}\left(  L\right)  =\left\{  \left\{  u_{t}\right\}  \text{:
}1/3\leq\inf_{\lambda\in\left[  -\pi,\pi\right]  }f\left(  \lambda\right)
\leq\sup_{\lambda\in\left[  -\pi,\pi\right]  }f\left(  \lambda\right)
\leq3\text{, }\sup_{\lambda\in\left[  -\pi,\pi\right]  }\left\vert f^{\prime
}\left(  \lambda\right)  \right\vert \leq L,\text{ }\sum_{j=0}^{\infty
}\left\vert R_{j}\right\vert \leq L\right\}  .
\]
The next Lemma describes a family of alternatives which satisfies Assumption
\ref{Reg} uniformly for prescribed constants and a given $\delta_{a}\left(
j\right)  .$

\begin{lem}
\label{Altexp} Consider a centered stationary Gaussian process $\left\{
u_{t}\right\}  $ with spectral density function $f\left(  \lambda\right)
=\exp\left(  g\left(  \lambda\right)  \right)  /\left(  2\pi\right)  $, where
\begin{equation}
g\left(  \lambda\right)  =2\rho\sum_{k=1}^{p}b_{k}\cos\left(  k\lambda\right)
,\text{\quad\quad\quad\quad}b_{k}=-1,0,1. \label{Logsd}%
\end{equation}
If $p\geq1$ and $\rho\geq0$ are such that $p^{2}\rho\leq\epsilon\leq1/6$ then
there is some constant $L>0$, independent of $\epsilon$, $p$, $\rho$ and
$b=\left(  b_{k},k\in\left[  1,p\right]  \right)  $, such that (i) $\left\vert
R_{0}-1\right\vert \leq6\rho\epsilon$ and $\left\vert R_{j}-\rho
b_{j}\right\vert \leq6\rho\epsilon$ for $j\in\left[  1,p\right]  $; (ii)
$\left\vert R_{j}\right\vert \leq3\rho\left(  2\epsilon\right)  ^{\ell}$ for
all $j$ in $\left[  \ell p+1,\left(  \ell+1\right)  p\right)  $ and all
$\ell\geq1$; (iii) $\left\{  u_{t}\right\}  $ is in H\"{o}lder$\left(
L\right)  $; (iv) Suppose that $\rho_{n}^{2}=\rho_{n}^{2}(p)=2\kappa_{n}%
^{2}\left(  2\log\log n\right)  ^{1/2}/\left(  np^{1/2}\right)  $ for some
$\kappa_{n}>0$ and bounded away from infinity, and that $p\in\left[
1,P_{n}\right]  $ with $P_{n}=o\left(  \left(  n/\left(  \kappa_{n}^{2}%
\log\log n\right)  ^{1/2}\right)  ^{1/14}\right)  $. Then the associated
family of processes $\left\{  u_{t}\left(  b,p\right)  ;b\in\left\{
-1,0,1\right\}  ^{p},p\in\left[  1,P_{n}\right]  \right\}  $ satisfies
Assumption \ref{Reg} for any $a>0$ and a $\delta_{a}\left(  j\right)
=O\left(  j^{-7-1/4}\right)  $.
\end{lem}

\noindent\textbf{Proof of Lemma \ref{Altexp}. }Rewrite $g$ as $g\left(
\lambda\right)  =\rho\sum_{k=-p}^{p}b_{k}\exp\left(  ik\lambda\right)  $,
$b_{0}=0$, $b_{k}=b_{-k}=b_{\left\vert k\right\vert }$. Since $\exp\left(
x\right)  =\sum_{m=0}^{\infty}x^{m}/m!$ uniformly over any compact set and
$\max_{\lambda}\left\vert g\left(  \lambda\right)  \right\vert \leq2p\rho
\leq2\epsilon\leq1/3$, we have
\begin{equation}
R_{j}=\int_{-\pi}^{\pi}\exp\left(  -ij\lambda\right)  f\left(  \lambda\right)
d\lambda=\frac{1}{2\pi}\sum_{m=0}^{\infty}\frac{1}{m!}\int_{-\pi}^{\pi}%
\exp\left(  -ij\lambda\right)  \left(  g\left(  \lambda\right)  \right)
^{m}d\lambda. \label{Ralt}%
\end{equation}
For $m>0$, since $\int_{-\pi}^{\pi}\exp\left(  -ij\lambda\right)
d\lambda=2\pi$ if $j=0$ and $0$ if $j\neq0$,
\begin{align}
&  \frac{1}{2\pi}\int_{-\pi}^{\pi}\exp\left(  -ij\lambda\right)  \left(
g\left(  \lambda\right)  \right)  ^{m}d\lambda\nonumber\\
&  =\frac{\rho^{m}}{2\pi}\sum_{\left(  k_{1},...,k_{m}\right)  \in K_{m}%
}b_{k_{1}}\times\cdots\times b_{k_{m}}\int_{-\pi}^{\pi}\exp\left(  i\left(
k_{1}+\ldots+k_{m}-j\right)  \lambda\right)  d\lambda\nonumber\\
&  =\rho^{m}\sum_{\left(  k_{1},...,k_{m}\right)  \in K_{m}\left(  j\right)
}b_{k_{1}}\times\cdots\times b_{k_{m}}, \label{Fourierg}%
\end{align}
where $K_{m}$ is the set of $m$-tuples with entries in $\left[  -p,p\right]
\setminus\left\{  0\right\}  $ so that $\#K_{m}=\left(  2p\right)  ^{m}$ and
$K_{m}\left(  j\right)  $ contains $m$-tuples in $K_{m}$ for which
$k_{1}+\cdots+k_{m}=j$ so that $\#K_{m}(j)\leq\left(  2p\right)  ^{m-1}$.

\textit{Proof of (i).} Part (i) is a consequence of (\ref{Ralt}),
(\ref{Fourierg}) and inequality $2p\rho\leq2\epsilon<1$ which together imply
that for $j\in\left[  0,p\right]  $, $\left\vert R_{j}-\mathbb{I}\left(
j=0\right)  -\rho b_{j}\right\vert \leq\rho\sum_{m=2}^{\infty}\frac{\left(
2p\rho\right)  ^{m-1}}{m!}\leq2p\rho^{2}\sum_{m=0}^{\infty}\frac{1}{m!}%
\leq2e\rho\epsilon<6\rho\epsilon$.

\textit{Proof of (ii).} Let $\ell p+1\leq j>\left(  \ell+1\right)  p$. Observe
that $K_{m}\left(  j\right)  $ is an empty set when $m\leq\ell$. Hence it
follows from (\ref{Ralt}) and (\ref{Fourierg}) that $\left\vert R_{j}%
\right\vert \leq\left\vert \frac{1}{2\pi}\sum_{m=\ell+1}^{\infty}\frac{1}%
{m!}\int_{-\pi}^{\pi}\exp\left(  -ij\lambda\right)  \left(  g\left(
\lambda\right)  \right)  ^{m}d\lambda\right\vert \leq\rho\sum_{m=\ell
+1}^{\infty}\frac{\left(  2p\rho\right)  ^{m-1}}{m!}\leq\rho\left(
2\epsilon\right)  ^{\ell}e$.

\textit{Proof of (iii). }Observe that $\left\vert g\left(  \lambda\right)
\right\vert \leq2\rho p\leq2\epsilon\leq1/3$ and that therefore%
\[
1/3<1-1/3<\exp\left(  -1/3\right)  \leq f\left(  \lambda\right)  \leq
\exp\left(  1/3\right)  \leq e\leq3\quad\quad\quad\quad\text{for all }%
\lambda\in\left[  -\pi,\pi\right]  .
\]
Parts (i), (ii) and $0\leq\rho\leq\epsilon<1/6$, $p\rho\leq1/6$ yield that,
for $L$ large enough,%
\begin{align*}
\sum_{j=0}^{\infty}\left\vert R_{j}\right\vert  &  \leq R_{0}+\sum_{j=1}%
^{p}\left\vert R_{j}\right\vert +\sum_{\ell=1}^{\infty}\sum_{j=\ell
p+1}^{\left(  \ell+1\right)  p}\left\vert R_{j}\right\vert \leq1+6\rho
\epsilon+\left(  1+6\epsilon\right)  p\rho+3\sum_{\ell=1}^{\infty}\left(
\ell+1\right)  p\rho\left(  2\epsilon\right)  ^{\ell}\\
&  \leq1+1+1+1+\sum_{\ell=1}^{\infty}\left(  \ell+1\right)  \left(
2\epsilon\right)  ^{\ell}\leq L.
\end{align*}
Since $f^{\prime}\left(  \lambda\right)  =g^{\prime}\left(  \lambda\right)
f\left(  \lambda\right)  $ with $g^{\prime}\left(  \lambda\right)  =-2\rho
\sum_{k=1}^{p}b_{k}k\sin\left(  k\lambda\right)  $, we have $\sup_{\lambda
\in\left[  -\pi,\pi\right]  }\left\vert f^{\prime}\left(  \lambda\right)
\right\vert \leq3\times2p^{2}\rho\leq1$.

\textit{Proof of (iv). }Let $u_{t}=\varepsilon_{t}+\sum_{j=1}^{\infty}\psi
_{j}\varepsilon_{t-j}$ be the Wold decomposition of the process. Brillinger
(2001) and $\int_{-\pi}^{\pi}\log f\left(  \lambda\right)  \exp\left(
ij\lambda\right)  d\lambda/2\pi=\rho b_{j}$ gives%
\begin{align*}
\psi_{j}  &  =\frac{\int_{-\pi}^{\pi}\exp\left(  \rho\sum_{k=1}^{p}b_{k}%
\exp\left(  -ik\lambda\right)  \right)  \exp\left(  ij\lambda\right)
d\lambda}{\int_{-\pi}^{\pi}\exp\left(  \rho\sum_{k=1}^{p}b_{k}\exp\left(
-ik\lambda\right)  \right)  d\lambda},\\
\operatorname*{Var}\left(  \varepsilon_{t}\right)   &  =\left\vert \frac
{1}{2\pi}\int_{-\pi}^{\pi}\exp\left(  \rho\sum_{k=1}^{p}b_{k}\exp\left(
-ik\lambda\right)  \right)  d\lambda\right\vert ^{2}.
\end{align*}
Arguing as in (i) and (ii) with an expansion as in (\ref{Ralt}) give
$\operatorname*{Var}\left(  \varepsilon_{t}\right)  =1$, $\left\vert \psi
_{j}-\rho b_{j}\right\vert \leq C\rho\epsilon$ for $j\in\left[  1,p\right]  $
and $\left\vert \psi_{j}\right\vert \leq C\rho\left(  2\epsilon\right)
^{\ell}$ for all $j\in\left[  \ell p+1,\left(  \ell+1\right)  p\right)  $ and
all $\ell\geq1$. Gaussianity, the choice of $\rho$ in (iv) with the
restriction on $P_{n}$ and Wu (2005) give, for any $a>1$, $\delta_{12a}\left(
j\right)  \leq C_{a}\left\vert \psi_{j}\right\vert \leq C_{a}j^{-7-1/4}$. That
the other conditions of Assumption \ref{Reg} hold uniformly in $p\in\left[
1,P_{n}\right]  $ follows from (i) and (ii).$\hfill\square$

We will now define a family $\mathcal{F}_{n}$\ of correlated Gaussian
alternatives. We first introduce some notation. Consider $\widetilde{\gamma
}_{n}=\left(  2\ln\ln n\right)  ^{1/2}$ and $\mathcal{P}^{\prime}=\left\{
2^{j},j=1,\ldots,J_{n}\right\}  $, $2^{J_{n}}=P_{n}=o\left(  \overline{p}%
_{n}\wedge\left(  n/\widetilde{\gamma}_{n}\right)  ^{1/14}\right)  $ so that
$\mathcal{P}^{\prime}\subset\left[  1,\overline{p}_{n}\right]  $ for $n$ large
enough. Define also
\begin{equation}
\rho_{n}^{2}(p)=2\frac{\kappa_{n}^{2}\widetilde{\gamma}_{n}}{np^{1/2}}%
,\quad\widetilde{\rho}_{n}(p)=2\rho_{n}^{2}(p)\quad\epsilon_{n}=P_{n}^{2}%
\rho_{n}(P_{n})=\frac{\left(  \widetilde{\gamma}_{n}\right)  ^{1/2}\kappa
_{n}P_{n}^{7/4}}{n^{1/2}}=o\left(  1\right)  . \label{Rhop}%
\end{equation}
Since $p^{2}\rho_{n}(p)\leq\epsilon_{n}$ for all $p\in\mathcal{P}^{\prime}$,
$\epsilon_{n}$ plays the role of the real number $\epsilon$ of Lemma
\ref{Altexp} and we assume from now on that $n$ is so large that $\epsilon
_{n}\leq1/6$. Consider the following log-spectral density functions:%
\[
g\left(  \lambda;b,p\right)  =2\widetilde{\rho}_{n}(p)\sum_{k\in\left[
p,2p\right)  }b_{k}\cos\left(  k\lambda\right)  ,\quad b=\left(  b_{1}%
,\ldots,b_{P_{n}}\right)  \in\left\{  -1,1\right\}  ^{P_{n}},\quad
p\in\mathcal{P}^{\prime}.
\]
Functions $g$ are of the form specified in (\ref{Logsd}). Let $W$ be a
symmetric standard Brownian motion process. Consider a centered stationary
Gaussian processes%
\[
u_{t,n}\left(  b,p\right)  =\frac{1}{\left(  2\pi\right)  ^{1/2}}\int_{-\pi
}^{\pi}\exp\left(  \frac{g\left(  \lambda;b,p\right)  }{2}\right)  \exp\left(
it\lambda\right)  dW\left(  \lambda\right)  .
\]
Observe that $u_{t,n}\left(  0,p\right)  $ does not depend on $p$ and is a
Gaussian white noise process with variance 1. Let $\left\{  R_{j,n}\left(
b,p\right)  \right\}  $ denote the covariance function of $u_{t,n}\left(
b,p\right)  $. The family $\mathcal{F}_{n}$\ of Gaussian processes can now be
defined as\textbf{\ }%
\[
\mathcal{F}_{n}=\left\{  \left\{  u_{t,n}\left(  b,p\right)  \right\}
,b\in\left\{  -1,1\right\}  ^{P_{n}},p\in\mathcal{P}^{\prime}\right\}  .
\]
Lemma \ref{Altexp} implies that all sequences $\left\{  u_{t,n}\right\}  $ in
$\mathcal{F}_{n}$ satisfies Assumption \ref{Reg}\texttt{\ }and that
$\mathcal{F}_{n}\subset$H\"{o}lder$\left(  L\right)  $. We now study the
asymptotic behavior of the stochastic covariance sequence $\left\{
R_{j,n}\left(  B,P\right)  \right\}  $. Let $N_{n}\left(  b,p\right)  $ be as
in (\ref{Sparse2}), that is
\[
N_{n}\left(  b,p\right)  =N_{n}\left(  \left\{  u_{t,n}\left(  b,p\right)
\right\}  ,p,\rho_{n}\left(  p\right)  \right)  =\#\left\{  \left\vert
\frac{R_{j,n}\left(  b,p\right)  }{R_{0,n}\left(  b,p\right)  }\right\vert
\geq\rho_{n}\left(  p\right)  ,\text{ }j\in\left[  1,p\right]  \right\}  .
\]
Lemma \ref{Altexp}-(i,ii) and (\ref{Rhop}) gives that $N_{n}\left(
b,p\right)  =p/2$ for $n$ large enough and uniformly in $p=2^{j}\in
\mathcal{P}^{\prime}$, so that $\rho_{n}^{2}(p)=2\kappa_{n}^{2}\widetilde
{\gamma}_{n}/\left(  np^{1/2}\right)  =\kappa_{n}^{2}\widetilde{\gamma}%
_{n}p^{1/2}/\left(  nN_{n}\left(  b,p\right)  \right)  $. Hence the sequences
$\left\{  u_{t,n}\right\}  $ in $\mathcal{F}_{n}$ satisfies condition (i) in
Theorem \ref{Sparseopt}. Therefore the Theorem will be proved if we show that
$\sup_{T_{n}}\min_{\left\{  u_{t,n}\right\}  \in\mathcal{F}_{n}}%
\mathbb{P}\left(  T_{n}=0\right)  \leq\alpha+o\left(  1\right)  $, where
$\sup_{T_{n}}$ is a supremum over asymptotically $\alpha$-level tests. Since
the equivalence result of Golubev et al. (2010) holds over $\mathcal{F}%
_{n}\subset$H\"{o}lder$\left(  L\right)  $ this is equivalent to show that
$\sup_{T_{n}}\min_{\left\{  U_{n}\right\}  \in\mathcal{F}_{n}}\mathbb{Q}%
\left(  T_{n}=0\right)  \leq\alpha+o\left(  1\right)  $, $\mathbb{Q}$ being
the distribution of the continuous time regression model
\[
dU_{n}\left(  \lambda;b,p\right)  =g\left(  \lambda;b,p\right)  d\lambda
+2\pi^{1/2}\frac{dW\left(  \lambda\right)  }{n^{1/2}},\quad\quad\lambda
\in\left[  -\pi,\pi\right]  ,
\]
where $W\left(  \cdot\right)  $ is a Brownian motion over $\lambda\in\left[
-\pi,\pi\right]  $. This can be done as in Spokoiny (1996, Proof of Theorem
2.3) by bounding $\sup_{T_{n}}\min_{\left\{  U_{n}\right\}  \in\mathcal{F}%
_{n}}\mathbb{Q}\left(  T_{n}=0\right)  $ with a Bayes risk, based on the
choice of a uniform distribution for $p$ and a Bernoulli one for $b$%
.$\hfill\square$

\subsection{Proof of Lemma \ref{Wildlem}}

The first approximation $R_{0,n}=\sigma^{2}\left(  1+O\left(  \gamma_{n}%
P_{n}^{1/2}/n\right)  \right)  $ follows easily from the definition
(\ref{Wildalt}) of the alternative. To show that the second approximation is
valid, note that for $j=1,...,P_{n}$,%
\[
R_{j,n}=\frac{\nu\gamma_{n}^{1/2}}{n^{1/2}P_{n}^{1/4}}\psi_{j}\sigma
^{2}+\left(  \frac{\nu\gamma_{n}^{1/2}}{n^{1/2}P_{n}^{1/4}}\right)
^{2}\left(  \psi_{j+1}\psi_{1}+\cdot\cdot\cdot+\psi_{P_{n}}\psi_{P_{n}%
-j}\right)  \sigma^{2}.
\]
By the Cauchy-Schwarz inequality, $\left\vert \psi_{j+1}\psi_{1}+\cdot
\cdot\cdot+\psi_{P_{n}}\psi_{P_{n}-j}\right\vert \leq\sum_{k=1}^{P_{n}}%
\psi_{k}^{2}=O(P_{n})$\ for all $j=1,...,P_{n}$, hence, uniformly in
$j=1,...,P_{n}$,%
\[
R_{j,n}=\frac{\nu\gamma_{n}^{1/2}}{n^{1/2}P_{n}^{1/4}}\psi_{j}\sigma
^{2}+O\left(  \frac{\gamma_{n}P_{n}^{1/2}}{n}\right)  =\frac{\nu\gamma
_{n}^{1/2}}{n^{1/2}P_{n}^{1/4}}\psi_{j}\sigma^{2}+o\left(  \frac{\gamma
_{n}^{1/2}}{n^{1/2}P_{n}^{1/4}}\right)
\]
since $P_{n}=o((n/\gamma_{n})^{2/3})$.

For the expression of $\mathbb{E}\left[  u_{t}|u_{t-k},k\geq1\right]  $,
observe that (\ref{Wildalt}) gives, for $n$ large enough,%
\begin{align*}
&  \mathbb{E}\left[  u_{t}|u_{t-k},k\geq1\right]  =\frac{\nu\gamma_{n}^{1/2}%
}{n^{1/2}P_{n}^{1/4}}\sum_{k=1}^{P_{n}}\psi_{k}\varepsilon_{t-k}\\
&  \text{ }=\frac{\nu\gamma_{n}^{1/2}}{n^{1/2}P_{n}^{1/4}}\sum_{k=1}^{P_{n}%
}\psi_{k}\left(  u_{t-k}-\frac{\nu\gamma_{n}^{1/2}}{n^{1/2}P_{n}^{1/4}}%
\sum_{j=1}^{P_{n}}\psi_{j}\varepsilon_{t-k-j}\right) \\
&  \text{ }=\frac{\nu\gamma_{n}^{1/2}}{n^{1/2}P_{n}^{1/4}}\sum_{k=1}^{P_{n}%
}\psi_{k}u_{t-k}-\frac{\nu^{2}\gamma_{n}}{nP_{n}^{1/2}}\sum_{k=1}^{P_{n}}%
\psi_{k}\sum_{j=1}^{P_{n}}\psi_{j}\varepsilon_{t-k-j}.
\end{align*}
Now, since $\left\{  \varepsilon_{t}\right\}  $ is a strong white noise and
$\sum_{k=1}^{P_{n}}\psi_{k}^{2}=O\left(  P_{n}\right)  $,%
\begin{align*}
&  \frac{\nu^{2}\gamma_{n}}{nP_{n}^{1/2}}\sum_{k=1}^{P_{n}}\psi_{k}\sum
_{j=1}^{P_{n}}\psi_{j}\varepsilon_{t-k-j}=\frac{\nu^{2}\gamma_{n}}%
{nP_{n}^{1/2}}\sum_{\ell=2}^{2P_{n}}\left(  \sum_{k=1}^{\max\left(  P_{n}%
,\ell-1\right)  }\psi_{k}\psi_{\ell-k}\right)  \varepsilon_{t-\ell}\\
&  \text{ }=O_{\mathbb{P}}\left(  \left(  \frac{\gamma_{n}^{2}}{n^{2}P_{n}%
}\sum_{\ell=2}^{2P_{n}}\left(  \sum_{k=1}^{\max\left(  P_{n},\ell-1\right)
}\psi_{k}\psi_{\ell-k}\right)  ^{2}\right)  ^{1/2}\right) \\
&  \text{ }=O_{\mathbb{P}}\left(  \left(  \frac{\gamma_{n}^{2}\left(
\sum_{k=1}^{P_{n}}\psi_{k}^{2}\right)  ^{2}}{n^{2}}\right)  ^{1/2}\right)
=O_{\mathbb{P}}\left(  \frac{\gamma_{n}P_{n}}{n}\right)  ,
\end{align*}
which ends the proof of the Lemma.$\hspace{0in}\hfill\Box$

\subsection{Proof of Proposition \ref{WildCvM}}

Let us now check consistency of the test (\ref{Test}) under the assumption
that $\min_{k\in\left[  1,P_{n}\right]  }\left\vert \psi_{k}\sigma
^{2}\right\vert \geq1$. Define $\rho_{n}=\left(  \nu/2\right)  \gamma
_{n}^{1/2}/\left(  n^{1/2}P_{n}^{1/4}\right)  $. Lemma \ref{Wildlem} implies
that $N_{n}=P_{n}\left(  1+o(1)\right)  $ for such a $\rho_{n}$, which
therefore satisfies
\[
\rho_{n}=\left(  1+o\left(  1\right)  \right)  \left(  \nu/2\right)  \left(
\gamma_{n}P_{n}^{1/2}/N_{n}\right)  ^{1/2}/n^{1/2},
\]
so that (\ref{Sparse3}) asymptotically holds provided $\nu\geq3\kappa^{\ast}$
and the test is consistent if $1\leq P_{n}\leq\overline{p}_{n}/2$ by Theorem
\ref{Sparse} provided the considered alternatives satisfies Assumption
\ref{Reg}. Wu (2005) gives that the alternative (\ref{Wildalt}) satisfies for
any $a>0$,%
\[
\delta_{12a}\left(  j\right)  \leq C_{a}\frac{\nu\gamma_{n}^{1/2}}%
{n^{1/2}P_{n}^{1/4}}\left\vert \sigma\psi_{j}\right\vert \text{ for all }%
j\in\left[  1,P_{n}\right]  \text{, }\delta_{12a}\left(  j\right)  =0\text{
for all }j>P_{n}.
\]
Hence the condition $P_{n}=O\left(  \left(  n/\gamma_{n}\right)
^{1/14}\right)  $ gives that $\delta_{12a}\left(  j\right)  \leq Cj^{-7-1/4}$
since the $\left\vert \sigma\psi_{j}\right\vert $ are bounded away from
infinity. Moreover Gaussianity ensures that%
\[
\left\Vert u_{t,n}-\varepsilon_{t}\right\Vert _{12a}\leq C_{a}\sigma\left(
\frac{\nu^{2}\gamma_{n}}{nP_{n}^{1/2}}\sum_{k=1}^{P_{n}}\psi_{k}^{2}\right)
^{1/2}=O\left(  \frac{\nu\gamma_{n}^{1/2}P_{n}^{1/4}}{n^{1/2}}\right)
=o\left(  1\right)  ,
\]
which gives $\operatorname*{Var}\left(  u_{t,n}\right)  =\sigma^{2}+o\left(
1\right)  $ and $\max_{j\in\left[  1,n\right]  }\operatorname*{Var}^{2}\left(
u_{t,n}\right)  /\operatorname*{Var}\left(  u_{t,n}u_{t+j,n}\right)
=1+o\left(  1\right)  $ so that Assumption \ref{Reg} holds. This ends the
proof of Proposition \ref{WildCvM}-(i).

Consider now the other tests in Proposition \ref{WildCvM}-(ii). Define
$\widetilde{R}_{1,j}=\sum_{t=1}^{n-j}u_{t,n}u_{t+j,n}/n$, $\widetilde{R}%
_{0,j}=\sum_{t=1}^{n-j}\varepsilon_{t}\varepsilon_{t+j}/n$, $\widetilde{\tau
}_{1,j}^{2}=\sum_{t=1}^{n-j}u_{t,n}^{2}u_{t+j,n}^{2}/\left(  n-j\right)
-n\widetilde{R}_{1,j}^{2}/\left(  n-j\right)  $ and $\widetilde{\tau}%
_{0,j}^{2}=\sum_{t=1}^{n-j}\varepsilon_{t}^{2}\varepsilon_{t+j}^{2}/\left(
n-j\right)  -n\widetilde{R}_{0,j}^{2}/\left(  n-j\right)  $. Define also
$\eta_{t}=\eta_{t,n}=\nu\sum_{k=1}^{\infty}\psi_{k}\varepsilon_{t-k}$, setting
$\psi_{k}=0$ for $k>P_{n}$, so that $u_{t,n}=\varepsilon_{t}+\gamma_{n}%
^{1/2}\eta_{t}/\left(  n^{1/2}P_{n}^{1/4}\right)  $. We have%
\[
\left\vert \widetilde{R}_{j}-\widetilde{R}_{0,j}\right\vert \leq\frac
{\gamma_{n}^{1/2}}{n^{3/2}P_{n}^{1/4}}\left\vert \sum_{t=1}^{n-j}\eta
_{t}\varepsilon_{t+j}\right\vert +\frac{\gamma_{n}^{1/2}}{n^{3/2}P_{n}^{1/4}%
}\left\vert \sum_{t=1}^{n-j}\varepsilon_{t}\eta_{t+j}\right\vert +\frac
{\gamma_{n}}{n^{2}P_{n}^{1/2}}\left\vert \sum_{t=1}^{n-j}\eta_{t}\eta
_{t+j}\right\vert .
\]
The Burkholder inequality gives, for any $a>1$,%
\begin{align*}
&  \left\Vert \frac{\gamma_{n}^{1/2}}{n^{3/2}P_{n}^{1/4}}\sum_{t=1}^{n-j}%
\eta_{t}\varepsilon_{t+j}\right\Vert _{a}\leq C\frac{\gamma_{n}^{1/2}\left(
n-j\right)  ^{1/2}}{n^{3/2}P_{n}^{1/4}}\left\Vert \eta_{t}\right\Vert _{a}\leq
C\frac{\gamma_{n}^{1/2}P_{n}^{1/4}}{n},\\
&  \left\Vert \frac{\gamma_{n}^{1/2}}{n^{3/2}P_{n}^{1/4}}\sum_{t=1}%
^{n-j}\left(  \varepsilon_{t}\eta_{t+j}-\psi_{j}\varepsilon_{t}^{2}\right)
\right\Vert _{a}\leq\left\Vert \frac{\gamma_{n}^{1/2}}{n^{3/2}P_{n}^{1/4}}%
\sum_{t=1}^{n-j}\varepsilon_{t}\left(  \sum_{k=0}^{j-1}\psi_{j}\varepsilon
_{t+j-k}\right)  \right\Vert _{a}\\
&  +\left\Vert \frac{\gamma_{n}^{1/2}}{n^{3/2}P_{n}^{1/4}}\sum_{t=1}%
^{n-j}\left(  \sum_{k=j+1}^{\infty}\psi_{j}\varepsilon_{t+j-k}\right)
\varepsilon_{t}\right\Vert _{a}\leq C\frac{\gamma_{n}^{1/2}P_{n}^{1/4}}{n},\\
&  \left\Vert \frac{\gamma_{n}^{1/2}}{n^{3/2}P_{n}^{1/4}}\sum_{t=1}%
^{n-j}\left(  \varepsilon_{t}^{2}-\sigma^{2}\right)  \right\Vert _{a}\leq
C\frac{\gamma_{n}^{1/2}}{nP_{n}^{1/4}},\quad\left\Vert \frac{\gamma_{n}}%
{n^{2}P_{n}^{1/2}}\sum_{t=1}^{n}\eta_{t}^{2}\right\Vert _{a}\leq\frac
{\gamma_{n}}{nP_{n}^{1/2}}\leq C\frac{\gamma_{n}P_{n}^{1/2}}{n},
\end{align*}
for all $j$. Note also that $\left\vert \sum_{t=1}^{n-j}\eta_{t}\eta
_{t+j}\right\vert \leq\sum_{t=1}^{n}\eta_{t}^{2}$ and the Markov inequality
give for $a$ large enough, since $\gamma_{n}P_{n}^{1/2}=o(n^{1/4})$
\begin{align*}
&  \max_{j\in\left[  1,n\right]  }\left\vert \widetilde{R}_{1,j}-\widetilde
{R}_{0,j}\right\vert ^{a}=O_{\mathbb{P}}\left(  \max_{j\in\left[  1,n\right]
}\left\vert \widetilde{R}_{1,j}-\widetilde{R}_{0,j}\right\vert ^{a}\right) \\
&  \text{ }=O_{\mathbb{P}}\left(  \sum_{j=1}^{n}\left\Vert \frac{\gamma
_{n}^{1/2}}{n^{3/2}P_{n}^{1/4}}\sum_{t=1}^{n-j}\eta_{t}\varepsilon_{t+j}%
+\sum_{t=1}^{n-j}\varepsilon_{t}\eta_{t+j}\right\Vert _{a}^{a}+\left\Vert
\frac{\gamma_{n}}{n^{2}P_{n}^{1/2}}\sum_{t=1}^{n}\eta_{t}^{2}\right\Vert
_{a}^{a}\right) \\
&  \text{ }=O_{\mathbb{P}}\left(  n\left(  \frac{\gamma_{n}^{1/2}P_{n}^{1/4}%
}{n}\right)  ^{a}+\left(  \frac{\gamma_{n}P_{n}^{1/2}}{n}\right)  ^{a}\right)
=o_{\mathbb{P}}\left(  \frac{1}{n^{7a/8-1}}+\frac{1}{n^{3a/4}}\right) \\
&  \text{ }=o_{\mathbb{P}}\left(  \frac{1}{\left(  n\log n\right)  ^{a/2}%
}\right)  .
\end{align*}
Hence%
\begin{equation}
\max_{j\in\left[  1,n\right]  }\left\vert \widetilde{R}_{1,j}-\widetilde
{R}_{0,j}\right\vert =o_{\mathbb{P}}\left(  \frac{1}{\left(  n\log n\right)
^{1/2}}\right)  . \label{MaxtildR}%
\end{equation}
Arguing similarly for the $\widetilde{\tau}_{k,j}^{2}$ give, since
$J_{n}=O\left(  n^{1/2}\right)  $
\begin{equation}
\max_{j\in\left[  1,J_{n}\right]  }\left\vert \widetilde{\tau}_{1,j}%
^{2}-\widetilde{\tau}_{0,j}^{2}\right\vert =o_{\mathbb{P}}\left(  \frac
{1}{\left(  n\log n\right)  ^{1/2}}\right)  \quad,\max_{j\in\left[
1,J_{n}\right]  }\left\vert \widetilde{\tau}_{0,j}^{2}-\sigma^{4}\right\vert
=O_{\mathbb{P}}\left(  \frac{\log^{1/2}n}{n^{1/2}}\right)  ,
\label{Maxtildsig}%
\end{equation}
where the latter is from Proposition \ref{Covesti}. Note that (\ref{MaxtildR})
and (\ref{Maxtildsig}) gives (\ref{G01cov}). Let $W_{k,n}$, $CvM_{k,n}$,
$EL_{k,n}$ be the statistic computed under $G_{k}$, $k=0,1$, i.e. with
$\widetilde{R}_{0,j}/\widetilde{\tau}_{0,j}$ and $\widetilde{R}_{1,j}%
/\widetilde{\tau}_{1,j}$. Note that (\ref{MaxtildR}) and (\ref{Maxtildsig})
gives $W_{1,n}=W_{0,n}+o_{\mathbb{P}}\left(  1\right)  $. (\ref{G01cov}) and
Proposition \ref{Covesti} give%
\begin{align*}
&  \left\vert CvM_{1,n}-CvM_{0,n}\right\vert \leq\frac{2}{\pi^{2}}\sum
_{j=1}^{J_{n}}\frac{n\left\vert \left(  \widetilde{R}_{1,j}/\widetilde{\tau
}_{1,j}+\widetilde{R}_{0,j}/\widetilde{\tau}_{0,j}\right)  \left(
\widetilde{R}_{1,j}/\widetilde{\tau}_{1,j}-\widetilde{R}_{0,j}/\widetilde
{\tau}_{0,j}\right)  \right\vert }{j^{2}}\\
&  \text{ }\leq2\max_{j\in\left[  1,J_{n}\right]  }\frac{\left\vert
n^{1/2}\widetilde{R}_{0,j}\right\vert }{\widetilde{\tau}_{0,j}}\times
\max_{j\in\left[  1,J_{n}\right]  }\left\vert n^{1/2}\left(  \frac
{\widetilde{R}_{1,j}}{\widetilde{\tau}_{1,j}}-\frac{\widetilde{R}_{0,j}%
}{\widetilde{\tau}_{0,j}}\right)  \right\vert \frac{2}{\pi^{2}}\sum
_{j=1}^{J_{n}}\frac{1}{j^{2}}\\
&  \text{ }+\max_{j\in\left[  1,J_{n}\right]  }n\left(  \frac{\widetilde
{R}_{1,j}}{\widetilde{\tau}_{1,j}}-\frac{\widetilde{R}_{0,j}}{\widetilde{\tau
}_{0,j}}\right)  ^{2}\frac{2}{\pi^{2}}\sum_{j=1}^{J_{n}}\frac{1}{j^{2}}\\
&  \text{ }=n^{1/2}O_{\mathbb{P}}\left(  \left(  \frac{\log n}{n}\right)
^{1/2}\right)  n^{1/2}o_{\mathbb{P}}\left(  \frac{1}{\left(  n\log n\right)
^{1/2}}\right)  +no_{\mathbb{P}}\left(  \frac{1}{n\log n}\right)
=o_{\mathbb{P}}\left(  1\right)  ,
\end{align*}
Hence $CvM_{1,n}=CvM_{0,n}+o_{\mathbb{P}}\left(  1\right)  $. For $EL_{n}$,
$W_{1,n}=W_{0,n}+o_{\mathbb{P}}\left(  1\right)  $ and Xiao and Wu (2011)
gives that $\max_{j\in\left[  1,J_{n}\right]  }\left\vert \widetilde{R}%
_{k,j}/\widetilde{\tau}_{k,j}\right\vert \leq\left(  2\ln n\right)
^{1/2}\left(  1+o_{\mathbb{P}}\left(  1\right)  \right)  $ for $k=0,1$ so that
$\mathbb{P}\left(  \widehat{\gamma}_{EL}^{\ast}=\ln n\right)  \rightarrow1$
under $G_{0}$ and $G_{1}$.We now show that $\mathbb{P}\left(  \widehat{p}%
_{EL}^{\ast}=1\right)  \rightarrow1$ under $G_{0}$. Propositions \ref{MeanH1}
and \ref{VarH1}, (\ref{Maxtildsig}) give
\begin{align*}
&  \mathbb{P}\left(  \widetilde{p}_{0,EL}^{\ast}\neq1\right)  =\mathbb{P}%
\left(  \max_{p\in\left[  2,J_{n}\right]  }\frac{\widetilde{BP}_{0,p}^{\ast
}-\widetilde{BP}_{0,1}^{\ast}}{p-1}>\ln n\right)  +o\left(  1\right) \\
&  \text{ }=\mathbb{P}\left(  \left(  1+o_{\mathbb{P}}\left(  1\right)
\right)  \max_{p\in\left[  2,J_{n}\right]  }\frac{n\sum_{j=2}^{p}\widetilde
{R}_{0,j}^{2}/\sigma^{4}}{p-1}>\ln n\right)  +o\left(  1\right) \\
&  \text{ }=\mathbb{P}\left(  \frac{n\sum_{j=2}^{p}\widetilde{R}_{0,j}%
^{2}/\sigma^{4}}{p-1}>\frac{1}{2}\ln n\text{ for some }p\in\left[
2,J_{n}\right]  \right)  +o\left(  1\right) \\
&  \text{ }\leq\sum_{p=2}^{J_{n}}\mathbb{P}\left(  \frac{n\sum_{j=2}%
^{p}\left(  \widetilde{R}_{0,j}^{2}/\sigma^{4}-\mathbb{E}\left[  \widetilde
{R}_{0,j}^{2}/\sigma^{4}\right]  \right)  }{p-1}>\frac{1}{2}\ln n-\frac
{n\sum_{j=2}^{p}\mathbb{E}\left[  \widetilde{R}_{0,j}^{2}/\sigma^{4}\right]
}{p-1}\right)  +o\left(  1\right) \\
&  \text{ }\leq\sum_{p=2}^{J_{n}}\frac{\operatorname*{Var}\left(  \frac
{n\sum_{j=2}^{p}\left(  \widetilde{R}_{0,j}^{2}/\sigma^{4}-\mathbb{E}\left[
\widetilde{R}_{0,j}^{2}/\sigma^{4}\right]  \right)  }{p-1}\right)  }{\left(
\frac{1}{2}\ln n-\frac{1}{p-1}\sum_{j=2}^{p}\left(  1-j/n\right)  \right)
^{2}}+o\left(  1\right) \\
&  \text{ }\leq\frac{C}{\log^{2}n}\sum_{p=2}^{J_{n}}\frac{1}{p-1}+o\left(
1\right)  =O\left(  \frac{1}{\log n}\right)  +o\left(  1\right)  =o\left(
1\right)  .
\end{align*}
Now, observe that Proposition \ref{Covesti} and (\ref{G01cov}) give%
\begin{align*}
&  \max_{p\in\left[  2,J_{n}\right]  }\left\vert \frac{\widetilde{BP}%
_{0,p}^{\ast}-\widetilde{BP}_{0,1}^{\ast}}{p-1}-\frac{\widetilde{BP}%
_{1,p}^{\ast}-\widetilde{BP}_{1,1}^{\ast}}{p-1}\right\vert \leq\max
_{p\in\left[  2,J_{n}\right]  }\left\vert \frac{n\sum_{j=2}^{p}\left(
\widetilde{R}_{0,j}^{2}/\widetilde{\tau}_{0,j}^{2}-\widetilde{R}_{1,j}%
^{2}/\widetilde{\tau}_{1,j}^{2}\right)  }{p-1}\right\vert \\
&  \text{ }\leq2\max_{p\in\left[  2,J_{n}\right]  }\left\vert n^{1/2}%
\frac{\widetilde{R}_{0,j}}{\widetilde{\tau}_{0,j}}\right\vert \times\max
_{p\in\left[  2,J_{n}\right]  }\left\vert n^{1/2}\left(  \frac{\widetilde
{R}_{0,j}}{\widetilde{\tau}_{0,j}}-\frac{\widetilde{R}_{1,j}}{\widetilde{\tau
}_{1,j}}\right)  \right\vert +\left(  \max_{p\in\left[  2,J_{n}\right]
}\left\vert n^{1/2}\left(  \frac{\widetilde{R}_{0,j}}{\widetilde{\tau}_{0,j}%
}-\frac{\widetilde{R}_{1,j}}{\widetilde{\tau}_{1,j}}\right)  \right\vert
\right)  ^{2}\\
&  \text{ }=n^{1/2}O_{\mathbb{P}}\left(  \left(  \frac{\log n}{n}\right)
^{1/2}\right)  n^{1/2}o_{\mathbb{P}}\left(  \frac{1}{\left(  n\log n\right)
^{1/2}}\right)  +no_{\mathbb{P}}\left(  \frac{1}{n\log n}\right)
=o_{\mathbb{P}}\left(  1\right)  .
\end{align*}
This, since arguing as in the bound above gives $\max_{p\in\left[
2,J_{n}\right]  }\left\vert \left(  \widetilde{BP}_{0,p}^{\ast}-\widetilde
{BP}_{0,1}^{\ast}\right)  /\left(  p-1\right)  \right\vert =O_{\mathbb{P}%
}\left(  \log^{1/2}n\right)  $, implies that $\max_{p\in\left[  2,J_{n}%
\right]  }\left\vert \left(  \widetilde{BP}_{1,p}^{\ast}-\widetilde{BP}%
_{1,1}^{\ast}\right)  /\left(  p-1\right)  \right\vert \leq\log n$ with a
probability tending to $1$ and then $\mathbb{P}\left(  \widehat{p}_{EL}^{\ast
}=1\right)  \rightarrow1$ under $G_{1}$. Hence (\ref{G01cov}) gives that
$EL_{1,n}=\widetilde{BP}_{1,1}^{\ast}+o_{\mathbb{P}}\left(  1\right)
=\widetilde{BP}_{0,1}^{\ast}+o_{\mathbb{P}}\left(  1\right)  =EL_{0,n}%
+o_{\mathbb{P}}\left(  1\right)  $, so that $EL_{n}$ converges in distribution
to a Chi square one with one degree of freedom under $G_{0}$ and $G_{1}%
$.\hspace*{\fill}$\Box$

\section*{Supplementary Material B: Proofs of intermediary
results\label{Proofs of intermediary results}}

\setcounter{lem}{0} \setcounter{thm}{0} \renewcommand{\thelem}{B.\arabic{lem}}
\renewcommand{\thethm}{B.\arabic{thm}} \setcounter{equation}{0}
\renewcommand{\theequation}{B.\arabic{equation}} \setcounter{subsection}{0}
\renewcommand{\thesubsection}{B.\arabic{subsection}} \renewcommand{\thesubsubsection}{B.\arabic{subsection}.\arabic{subsubsection}}

The proofs also use the notion of cumulants, see for example Brillinger (2001,
p. 19) or Xiao and Wu (2011)\ for a definition. Let
\[
\mathrm{Cum}\left(  u_{t_{1,n}},\ldots,u_{t_{q,n}}\right)  =\Gamma_{n}%
(t_{1},\ldots,t_{q})
\]
stands for the $q$th cumulants of $\left\{  u_{t,n}\right\}  $. The next
theorem on cumulant summability is Theorem 21 in Xiao and Wu (2011). These
authors do not formally consider sequences $\left\{  u_{t,n}\right\}  $ but
the following result is a straightforward extension of Xiao and Wu (2011).

\begin{thm}
[Xiao and Wu (2011)]\label{XW11} Suppose $\left\{  u_{t,n}\right\}  $ is
stationary for each $n$, with
\[
\sup_{n}\left\Vert u_{t,n}\right\Vert _{q+1}<\infty\text{ and }\sup
_{n}\left\Vert u_{t,n}-u_{t,n}^{t-j}\right\Vert _{q}\leq\delta_{q}\left(
j\right)  \text{ where }\sum_{j=0}^{\infty}j^{q-2}\delta_{q}\left(  j\right)
<\infty.
\]
Then there is a $\mathcal{C}$ which only depends on $\sup_{n}\left\Vert
u_{t,n}\right\Vert _{q+1}$ and $\sum_{j=0}^{\infty}j^{q-2}\delta_{q}\left(
j\right)  $ such that%
\[
\sum_{t_{2},\ldots,t_{q}=-\infty}^{\infty}\left\vert \Gamma_{n}(0,t_{2}%
,\ldots,t_{q})\right\vert \leq\mathcal{C}\text{.}%
\]

\end{thm}

In what follows, we drop subscript $n$ in expressions like $u_{t,n}$,
$R_{j,n}$, $\Gamma_{n}\left(  \cdot\right)  $ and $\theta_{n}$ when there is
no ambiguity. We denote%
\begin{equation}
K_{jp}=K^{2}\left(  \frac{j}{p}\right)  -K^{2}\left(  j\right)  \quad
\quad\text{and}\quad\quad K_{1n}(p)=\sum_{j=1}^{n-1}K_{jp}. \label{kjdifcov}%
\end{equation}

\subsection{\textbf{Proof of Lemma \ref{Ordersums}}}

(i) The first three bounds of the lemma follow directly from Assumption
\ref{Kernel} which implies that $K^{2}\left(  j/p\right)  \geq K^{2}\left(
j\right)  $ for all $j$ and $\mathbb{I}(x\in\lbrack0,1/2])/C\leq K^{2q}(x)\leq
C\mathbb{I}(x\in\lbrack0,1])$ for some $C>0$. The Cauchy-Schwarz inequality
implies that for any $p\in\lbrack1,n/2]$, $E_{\Delta}(p)=\sum_{j=1}%
^{n-1}\left(  1-\frac{j}{n}\right)  K_{jp}\leq K_{1n}(p)\leq p^{1/2}\left(
\sum_{j=1}^{n-1}k_{j}^{2}(p)\right)  ^{1/2}\leq Cp^{1/2}V_{\Delta}(p)$, which
is the last bound in (i). (ii) Write $p=1+\nu$. Since $p\leq\overline{p}%
_{n}\leq n/2$, the support of $K\left(  \cdot\right)  $ is $\left[
0,1\right]  $ and $K\left(  \cdot\right)  $ is a decreasing function, we have
\begin{align*}
V_{\Delta}^{2}(p)  &  \geq\frac{1}{2}\times2\sum_{j=2}^{p}K^{2}\left(
\frac{j}{p}\right)  \geq\sum_{j=1}^{\nu}K^{2}\left(  \frac{1+j}{1+\nu}\right)
\geq\sum_{j=1}^{\nu}\int_{j}^{j+1}K^{2}\left(  \frac{1+x}{1+\nu}\right)  dx\\
&  =\int_{1}^{\nu+1}K^{2}\left(  \frac{1+x}{1+\nu}\right)  dx=\nu\int_{0}%
^{1}K^{2}\left(  \frac{2+z\nu}{1+\nu}\right)  dz.
\end{align*}
The map $\nu\longmapsto$ $\left(  2+z\nu\right)  /\left(  1+\nu\right)  $,
$z\in\left[  0,1\right)  $, is decreasing. Hence, for $\nu\geq2$, $V_{\Delta
}^{2}(p)\geq\nu\int_{0}^{1/2}K^{2}\left(  \frac{2+2z}{3}\right)  dz\geq
C\left(  p-1\right)  $. Now $V_{\Delta}^{2}(2)\geq2\left(  K^{2}\left(
\frac{1}{2}\right)  -K^{2}\left(  1\right)  \right)  ^{2}>0$ gives the desired
result for $V_{\Delta}(p)$. Since $K$ is nonincreasing, $p\longmapsto
E_{\Delta}(p)$ is non decreasing and $E_{\Delta}(p)\geq0$ for all
$p\in\mathcal{P}$.\hspace*{\fill}$\Box$

\subsection{Proof of Lemma \ref{L}}

Under $\mathcal{H}_{0}$, The proof repeats the steps of Lee (2007), Lobato
(2001) and Kuan and Lee (2006) using the joint FCLT of Assumption \ref{M}. The
joint FCLT of Assumption \ref{M} gives that the critical values are
$O_{\mathbb{P}}\left(  1\right)  $ under $\mathcal{H}_{1}$.\hspace*{\fill
}$\Box$

\subsection{\textbf{Proof of Lemma \ref{Varcov}.}}

Equation (5.3.21) in Priestley (1981) and Theorem \ref{XW11} gives uniformly
in $j$,
\begin{align*}
\operatorname*{Var}\left(  \widetilde{R}_{j}\right)   &  =\frac{1}{n}%
\sum_{j_{1}=-n+j+1}^{n-j-1}\left(  1-\frac{|j_{1}|+j}{n}\right)  \left(
R_{j_{1}}^{2}+R_{j_{1}+j}R_{j_{1}-j}+\Gamma\left(  0,j_{1},j,j_{1}+j\right)
\right) \\
&  \leq\frac{2}{n}\sum_{j_{1}=-2n}^{2n}R_{j_{1}}^{2}+\frac{1}{n}\sum
_{j_{2},j_{3},j_{4}=-\infty}^{+\infty}\left\vert \Gamma\left(  0,j_{2}%
,j_{3},j_{4}\right)  \right\vert \\
&  \leq\frac{4}{n}\sum_{j=0}^{\infty}R_{j}^{2}+\frac{1}{n}\sum_{j_{2}%
,j_{3},j_{4}=-\infty}^{+\infty}\left\vert \Gamma\left(  0,j_{2},j_{3}%
,j_{4}\right)  \right\vert <C.\hfill\square
\end{align*}

\subsection{\textbf{Proof of Proposition \ref{Covesti}}}

For the sake of brevity we assume that $\theta$ is unidimensional. That
\begin{align*}
\max_{j\in\left[  0,n-1\right]  }\left\vert \widetilde{R}_{j}-\left(
1-\frac{j}{n}\right)  R_{j,n}\right\vert  &  =O_{\mathbb{P}}\left(  \left(
\frac{\log n}{n}\right)  ^{1/2}\right)  ,\\
\max_{j\in\left[  0,n-1\right]  }\left(  1-\frac{j}{n}\right)  \left\vert
\widetilde{\tau}_{j}^{2}-\tau_{j,n}^{2}\right\vert  &  =O_{\mathbb{P}}\left(
\left(  \frac{\log n}{n}\right)  ^{1/2}\right)  ,
\end{align*}
follow from Xiao and Wu (2011, Theorem 2). Note that these authors do not
consider stationary sequences $\left\{  u_{t,n}\right\}  $ but their arguments
carry over under Assumption \ref{Reg}. Hence it suffices to study $\max
_{j\in\left[  0,\overline{p}_{n}\right]  }\left\vert \widehat{R}%
_{j}-\widetilde{R}_{j}\right\vert $ and $\max_{j\in\left[  0,\overline{p}%
_{n}\right]  }\left\vert \widehat{\tau}_{j}^{2}-\widetilde{\tau}_{j}%
^{2}\right\vert $ since $\overline{p}_{n}/n=o\left(  n^{-1/2}\right)  $ under
Assumption \ref{P}. We then now show that $\max_{j\in\left[  0,\overline
{p}_{n}\right]  }\left\vert \widehat{R}_{j}-\widetilde{R}_{j}\right\vert
=O_{\mathbb{P}}\left(  n^{-1/2}\right)  $. Let $e_{t}=\widehat{u}_{t}-u_{t}$,
so that%
\[
\widehat{R}_{j}=\frac{1}{n}\sum_{t=1}^{n-{j}}\left(  u_{t}+e_{t}\right)
\left(  u_{t+j}+e_{t+j}\right)  =\widetilde{R}_{j}+\frac{1}{n}\sum
_{t=1}^{n-{j}}\left(  u_{t}e_{t+j}+e_{t}u_{t+j}\right)  +\frac{1}{n}\sum
_{t=1}^{n-{j}}e_{t}e_{t+j}%
\]
with, by the Cauchy-Schwarz inequality, $\left\vert \sum_{t=1}^{n-{j}}%
e_{t}e_{t+j}\right\vert /n\leq\sum_{t=1}^{n}e_{t}^{2}/n$ and, under Assumption
\ref{M}, for $\widehat{\mathfrak{r}}_{t}=\mathfrak{r}_{t}\left(
\widehat{\theta}\right)  $,%
\[
\frac{1}{n}\sum_{t=1}^{n-{j}}u_{t}e_{t+j}=\left(  \widehat{\theta}%
-\theta\right)  \frac{1}{n}\sum_{t=1}^{n-{j}}u_{t}u_{t+j}^{\left(  1\right)
}+\frac{1}{2}\left(  \widehat{\theta}-\theta\right)  ^{2}\frac{1}{n}\sum
_{t=1}^{n-{j}}u_{t}u_{t+j}^{\left(  2\right)  }+\frac{1}{n}\sum_{t=1}^{n-{j}%
}u_{t}\widehat{\mathfrak{r}}_{t+j}.
\]
Now, observe that Assumption \ref{M} gives $\widehat{\theta}-\theta
=O_{\mathbb{P}}\left(  n^{-1/2}\right)  $, $\max_{t\in\left[  1,n\right]
}\left\vert \widehat{\mathfrak{r}}_{t}\right\vert =o_{\mathbb{P}}\left(
1/n\right)  $ and%
\[
\frac{1}{n}\sum_{t=1}^{n}e_{t}^{2}\leq3\left(  \widehat{\theta}-\theta\right)
^{2}\frac{1}{n}\sum_{t=1}^{n}\left(  u_{t}^{\left(  1\right)  }\right)
^{2}+\frac{3}{4}\left(  \widehat{\theta}-\theta\right)  ^{4}\frac{1}{n}%
\sum_{t=1}^{n}\left(  u_{t}^{\left(  1\right)  }\right)  ^{2}+\frac{3}{n}%
\sum_{t=1}^{n}\left\vert \widehat{\mathfrak{r}}_{t}\right\vert =O_{\mathbb{P}%
}\left(  \frac{1}{n}\right)  ,
\]%
\[
\max_{j\in\left[  1,n\right]  }\left\vert \frac{1}{n}\sum_{t=1}^{n-{j}}\left(
u_{t}\widehat{\mathfrak{r}}_{t+j}+u_{t+j}\widehat{\mathfrak{r}}_{t}\right)
\right\vert \leq\frac{2\max_{t\in\left[  1,n\right]  }\left\vert
\widehat{\mathfrak{r}}_{t}\right\vert }{n}\sum_{t=1}^{n-{j}}\left\vert
u_{t}\right\vert =o_{\mathbb{P}}\left(  \frac{1}{n}\right)  .
\]
This gives, uniformly in $j\in\left[  1,n\right]  $%
\begin{align}
&  \left\vert \widehat{R}_{j}-\widetilde{R}_{j}\right\vert \leq\left\vert
\widehat{\theta}-\theta\right\vert \left\vert \mathbb{E}\left[  u_{t}%
u_{t+j}^{\left(  1\right)  }+u_{t+j}u_{t}^{\left(  1\right)  }\right]
\right\vert \nonumber\\
&  \text{ }+\left\vert \widehat{\theta}-\theta\right\vert \left\vert \frac
{1}{n}\sum_{t=1}^{n-{j}}\left(  u_{t}u_{t+j}^{\left(  1\right)  }+u_{t+j}%
u_{t}^{\left(  1\right)  }-\mathbb{E}\left[  u_{t}u_{t+j}^{\left(  1\right)
}+u_{t+j}u_{t}^{\left(  1\right)  }\right]  \right)  \right\vert
+O_{\mathbb{P}}\left(  \frac{1}{n}\right)  . \label{HatR2tildeR}%
\end{align}
It also follows from Assumption \ref{M} and $\overline{p}_{n}=o\left(
n^{1/2}\right)  $ that $\left\vert \widehat{\theta}-\theta\right\vert
\max_{j\in\left[  1,n\right]  }\left\vert \mathbb{E}\left[  u_{t}%
u_{t+j}^{\left(  1\right)  }+u_{t+j}u_{t}^{\left(  1\right)  }\right]
\right\vert =O_{\mathbb{P}}\left(  1/n^{1/2}\right)  $, $n\left(
\widehat{\theta}-\theta\right)  ^{2}\sum_{j=0}^{\infty}\mathbb{E}^{2}\left[
u_{t}u_{t+j}^{\left(  1\right)  }+u_{t+j}u_{t}^{\left(  1\right)  }\right]
=O_{\mathbb{P}}\left(  1\right)  $, and for $A_{t}\left(  j\right)
=u_{t}u_{t+j}^{\left(  1\right)  }+u_{t+j}u_{t}^{\left(  1\right)
}-\mathbb{E}\left[  u_{t}u_{t+j}^{\left(  1\right)  }+u_{t+j}u_{t}^{\left(
1\right)  }\right]  $%
\begin{align*}
&  \left\vert \widehat{\theta}-\theta\right\vert \max_{j\in\left[
0,\overline{p}_{n}\right]  }\left\vert \frac{1}{n}\sum_{t=1}^{n-{j}}%
A_{t}\left(  j\right)  \right\vert \leq O_{\mathbb{P}}\left(  \frac{1}%
{n^{1/2}}\right)  \sum_{j=0}^{\overline{p}_{n}}\left\vert \frac{1}{n}%
\sum_{t=1}^{n-{j}}A_{t}\left(  j\right)  \right\vert \\
&  \text{ }=O_{\mathbb{P}}\left(  \frac{1}{n}\right)  O_{\mathbb{P}}\left(
\sum_{j=0}^{\overline{p}_{n}}\mathbb{E}^{1/2}\left[  \left(  \frac{1}{n^{1/2}%
}\sum_{t=1}^{n-{j}}A_{t}\left(  j\right)  \right)  ^{2}\right]  \right) \\
&  \text{ }=O_{\mathbb{P}}\left(  \frac{1}{n}\right)  O_{\mathbb{P}}\left(
\overline{p}_{n}\max_{j\in\left[  0,\overline{p}_{n}\right]  }\left[  \left(
\frac{1}{n^{1/2}}\sum_{t=1}^{n-{j}}A_{t}\left(  j\right)  \right)
^{2}\right]  \right)  =O_{\mathbb{P}}\left(  \frac{1}{n^{1/2}}\right)  ,
\end{align*}%
\begin{align*}
&  n\sum_{j=0}^{n-1}\left(  \widehat{\theta}-\theta\right)  ^{2}\left(
\frac{1}{n}\sum_{t=1}^{n-{j}}A_{t}\left(  j\right)  \right)  ^{2}\\
&  \text{ }=O_{\mathbb{P}}\left(  1\right)  \frac{1}{n}O_{\mathbb{P}}\left(
\sum_{j=0}^{n-1}\mathbb{E}\left[  \left(  \frac{1}{n^{1/2}}\sum_{t=1}^{n-{j}%
}A_{t}\left(  j\right)  \right)  ^{2}\right]  \right) \\
&  \text{ }=O_{\mathbb{P}}\left(  1\right)  \frac{1}{n}O_{\mathbb{P}}\left(
n\max_{j\in\left[  0,n\right]  }\mathbb{E}\left[  \left(  \frac{1}{n^{1/2}%
}\sum_{t=1}^{n-{j}}A_{t}\left(  j\right)  \right)  ^{2}\right]  \right)
=O_{\mathbb{P}}\left(  1\right)  .
\end{align*}
This gives $\max_{j\in\left[  0,\overline{p}_{n}\right]  }\left\vert
\widehat{R}_{j}-\widetilde{R}_{j}\right\vert =O_{\mathbb{P}}\left(
n^{-1/2}\right)  $ and $\max_{p\in\left[  0,n-1\right]  }n\sum_{j=1}%
^{p}\left(  \widehat{R}_{j}-\widetilde{R}_{j}\right)  ^{2}=$ $O_{\mathbb{P}%
}\left(  1\right)  $. The study of $\max_{j\in\left[  0,\overline{p}%
_{n}\right]  }\left\vert \widehat{\tau}_{j}^{2}-\widetilde{\tau}_{j}%
^{2}\right\vert $ is similar.$\hfill\square$

\subsection{\textbf{Proof of Proposition \ref{Esti}}}

For the sake of brevity we assume that $\theta$ is unidimensional. Since
$\widehat{R}_{j}^{2}-\widetilde{R}_{j}^{2}=\left(  \widehat{R}_{j}%
-\widetilde{R}_{j}\right)  ^{2}+2\widetilde{R}_{j}\left(  \widehat{R}%
_{j}-\widetilde{R}_{j}\right)  $, Proposition \ref{Esti} is a direct
consequence of Proposition \ref{Covesti} and Lemma \ref{Crosscov} below.

\begin{lem}
\label{Crosscov}Assume that Assumptions \ref{Kernel}, \ref{M}, \ref{P} and
\ref{Reg} hold. Then
\[
\max_{p\in\left[  2,\overline{p}_{n}\right]  }\frac{\left\vert n\sum
_{j=1}^{n-1}\left(  K^{2}(j/p)-K^{2}(j)\right)  \widetilde{R}_{j}\left(
\widehat{R}_{j}-\widetilde{R}_{j}\right)  \right\vert }{\left(  1+n\sum
_{j=1}^{p}R_{j}^{2}\right)  ^{1/2}}=O_{\mathbb{P}}\left(  1\right)
\]
and $n\sum_{j=1}^{n-1}K^{2}(j/p_{n})\widetilde{R}_{j}\left(  \widehat{R}%
_{j}-\widetilde{R}_{j}\right)  =O_{\mathbb{P}}\left(  \left(  1+n\sum
_{j=1}^{p_{n}}R_{j}^{2}\right)  ^{1/2}\right)  $ for any $p_{n}=O(n^{1/2})$.
\end{lem}

{\small \noindent}\textbf{Proof of Lemma \ref{Crosscov}.}{\small \ }We just
prove the first equality since the proof of the second is very similar. Define
$\overline{R}_{j}=\mathbb{E}\left[  \widetilde{R}_{j}\right]  =(1-j/n)R_{j}$.
We have
\begin{align*}
&  \left\vert n\sum_{j=1}^{n-1}K_{jp}\widetilde{R}_{j}\left(  \widehat{R}%
_{j}-\widetilde{R}_{j}\right)  \right\vert \leq C_{n}(p)+D_{n}(p),\text{
where}\\
&  C_{n}(p)=\left\vert n\sum_{j=1}^{n-1}K_{jp}R_{j}\left(  \widehat{R}%
_{j}-\widetilde{R}_{j}\right)  \right\vert ,\\
&  D_{n}(p)=\left\vert n\sum_{j=1}^{n-1}K_{jp}\left(  \widetilde{R}%
_{j}-\overline{R}_{j}\right)  \left(  \widehat{R}_{j}-\widetilde{R}%
_{j}\right)  \right\vert .
\end{align*}
The Cauchy-Schwarz inequality and Assumption \ref{Kernel} gives
\[
C_{n}(p)\leq C\left(  n\sum_{j=1}^{p}R_{j}^{2}\right)  ^{1/2}\left(
n\sum_{j=1}^{p}\left(  \widehat{R}_{j}-\widetilde{R}_{j}\right)  ^{2}\right)
^{1/2}.
\]
Hence Proposition \ref{Covesti} yields that $\max_{p\in\left[  2,\overline
{p}_{n}\right]  }|C_{n}(p)/\left(  n\sum_{j=1}^{p}R_{j}^{2}\right)
^{1/2}|=O_{\mathbb{P}}\left(  1\right)  $. For $D_{n}(p)$, Assumptions
\ref{Kernel}, \ref{M}, (\ref{HatR2tildeR}) and $\widehat{\mathfrak{r}}%
_{t}=\mathfrak{r}_{t}\left(  \widehat{\theta}\right)  $ give
\begin{align*}
\max_{p\in\left[  2,\overline{p}_{n}\right]  }D_{n}(p)  &  \leq O_{\mathbb{P}%
}(n^{-1/2})\left(  \max_{p\in\left[  2,\overline{p}_{n}\right]  }%
D_{1n}(p)+\max_{p\in\left[  2,\overline{p}_{n}\right]  }D_{2n}(p)\right)
+O_{\mathbb{P}}(n^{-1})\max_{p\in\left[  2,\overline{p}_{n}\right]  }%
D_{3n}(p)\\
&  +\left(  \frac{1}{n}\sum_{t=1}^{n}e_{t}^{2}+2\frac{\max_{t\in\left[
1,n\right]  }\left\vert \mathfrak{r}_{t}\right\vert }{n}\sum_{t=1}%
^{n}\left\vert u_{t}\right\vert \right)  \max_{p\in\left[  2,\overline{p}%
_{n}\right]  }D_{4n}(p),
\end{align*}
where $D_{1n}(p)=n\sum_{j=1}^{p}\left\vert \widetilde{R}_{j}-\overline{R}%
_{j}\right\vert \left\vert \mathbb{E}\left[  u_{t}u_{t+j}^{(1)}+u_{t+j}%
u_{t}^{(1)}\right]  \right\vert $,
\begin{align*}
D_{2n}(p)  &  =n\sum_{j=1}^{p}\left\vert \widetilde{R}_{j}-\overline{R}%
_{j}\right\vert \left\vert \frac{1}{n}\sum_{t=1}^{n-j}\left(  u_{t}%
u_{t+j}^{(1)}+u_{t+j}u_{t}^{(1)}-\mathbb{E}\left[  u_{t}u_{t+j}^{(1)}%
+u_{t+j}u_{t}^{(1)}\right]  \right)  \right\vert ,\\
D_{3n}(p)  &  =n\sum_{j=1}^{p}\left\vert \widetilde{R}_{j}-\overline{R}%
_{j}\right\vert \left\vert \frac{1}{n}\sum_{t=1}^{n-j}\left(  u_{t}%
u_{t+j}^{(2)}+u_{t+j}u_{t}^{(2)}\right)  \right\vert ,\\
D_{4n}(p)  &  =n\sum_{j=1}^{p}\left\vert \widetilde{R}_{j}-\overline{R}%
_{j}\right\vert .
\end{align*}
By Assumption \ref{Kernel} and \ref{M} and by Lemma \ref{Varcov}, we have%
\[
\mathbb{E}\left[  \max_{p\in\left[  2,\overline{p}_{n}\right]  }%
D_{1n}(p)\right]  \leq Cn\sum_{j=1}^{\overline{p}_{n}}\operatorname*{Var}%
\nolimits^{1/2}\left(  \widetilde{R}_{j}\right)  \left\vert \mathbb{E}\left[
u_{t}u_{t+j}^{(1)}+u_{t+j}u_{t}^{(1)}\right]  \right\vert \leq Cn^{1/2},
\]%
\begin{align*}
\mathbb{E}\left[  \max_{p\in\left[  2,\overline{p}_{n}\right]  }%
D_{2n}(p)\right]   &  \leq Cn^{1/2}\sum_{j=1}^{\overline{p}_{n}}%
\operatorname*{Var}\nolimits^{1/2}\left(  \widetilde{R}_{j}\right) \\
&  \times\mathbb{E}^{1/2}\left[  \left\vert \frac{1}{n^{1/2}}\sum_{t=1}%
^{n}\left(  u_{t}u_{t+j}^{(1)}+u_{t+j}u_{t}^{(1)}-\mathbb{E}\left[
u_{t}u_{t+j}^{(1)}+u_{t+j}u_{t}^{(1)}\right]  \right)  \right\vert ^{2}\right]
\\
&  \leq C\overline{p}_{n},
\end{align*}%
\[
\mathbb{E}\left[  \max_{p\in\left[  2,\overline{p}_{n}\right]  }%
D_{3n}(p)\right]  \leq Cn\sum_{j=1}^{\overline{p}_{n}}\operatorname*{Var}%
\nolimits^{1/2}\left(  \widetilde{R}_{j}\right)  \mathbb{E}^{1/2}\left[
\left\vert \frac{1}{n}\sum_{t=1}^{n}\left(  u_{t}u_{t+j}^{(2)}+u_{t+j}%
u_{t}^{(2)}\right)  \right\vert ^{2}\right]  \leq C\overline{p}_{n}n^{1/2},
\]%
\[
\mathbb{E}\left[  \max_{p\in\left[  2,\overline{p}_{n}\right]  }%
D_{4n}(p)\right]  \leq Cn\sum_{j=1}^{\overline{p}_{n}}\mathbb{E}\left[
\left\vert \widetilde{R}_{j}-\overline{R}_{j}\right\vert \right]  \leq
Cn\sum_{j=1}^{\overline{p}_{n}}\operatorname*{Var}\nolimits^{1/2}\left(
\widetilde{R}_{j}\right)  \leq Cn^{1/2}\overline{p}_{n}.
\]
The Markov inequality gives us the stochastic orders of magnitude of the four
maxima in the bound for $\max_{p\in\left[  2,\overline{p}_{n}\right]  }%
D_{n}(p)$. Since $\overline{p}_{n}=O\left(  n^{1/2}\right)  $ by Assumption
\ref{P}, $\max_{t\in\left[  1,n\right]  }\left\vert \widehat{\mathfrak{r}}%
_{t}\right\vert =o_{\mathbb{P}}\left(  1/n\right)  $ and $n^{-1}\sum_{t=1}%
^{n}e_{t}^{2}=O_{\mathbb{P}}(n^{-1})$\ by Assumption \ref{M}, we have
$\max_{p\in\left[  2,\overline{p}_{n}\right]  }\left\vert D_{n}(p)\right\vert
=O_{\mathbb{P}}\left(  1+\frac{\overline{p}_{n}}{n^{1/2}}\right)
=O_{\mathbb{P}}\left(  1\right)  $. This together with $\max_{p\in\left[
2,\overline{p}_{n}\right]  }|C_{n}(p)/\left(  n\sum_{j=1}^{p}R_{j}^{2}\right)
^{1/2}|=O_{\mathbb{P}}\left(  1\right)  $ shows that the Lemma is
proved.\hspace*{\fill}$\Box$

\subsection{Proof of Proposition \ref{Selection}}

The proof of Proposition \ref{Selection} is long and divided in three steps.
In the two first steps, we focus on observed variables. In the first step, we
approximate the sample covariance $\widetilde{R}_{j}$ by a martingale
counterpart $\sum_{t=1}^{n}D_{jt}/n$, $j\in\left[  1,\overline{p}_{n}\right]
$, as in Shao (2011b), see the notations below and Lemmas \ref{CrossDD},
\ref{Sumvareta}. and \ref{Momsig}. The second step deals with the deviation
probability of
\[
\frac{n\sum_{j=1}^{p}\left(  \frac{1}{n}\sum_{t=j+1}^{n}D_{jt}\right)
^{2}\left(  K^{2}\left(  j/p\right)  -K^{2}\left(  1\right)  \right)
-\sigma^{4}E_{\Delta}\left(  p\right)  }{\sigma^{4}V_{\Delta}\left(  p\right)
}%
\]
which is approximated with some Gaussian counterparts through the Lindeberg
technique, see Lemma \ref{Lindeberg}. The third step concludes and explicitly
deals with the case of residuals thanks to Propositions \ref{Covesti} and
\ref{Esti}.

Let us now introduce additional notations. Let $\mathcal{F}_{k}$ be the sigma
field generated by $e_{k},e_{k-1},\ldots$. Define $\mathbf{P}_{t}\left[
Z\right]  =\mathbb{E}\left[  Z\left\vert \mathcal{F}_{t}\right.  \right]
-\mathbb{E}\left[  Z\left\vert \mathcal{F}_{t-1}\right.  \right]  $. Wu (2007,
Proposition 3) establishes that $\left\Vert \mathbf{P}_{t}\left[
u_{t+k}\right]  \right\Vert _{a}\leq\delta_{a}\left(  k\right)  $ and Shao
(2011b) has shown that%
\begin{equation}
\left\Vert \mathbf{P}_{0}\left[  u_{k}u_{k-j}\right]  \right\Vert _{a}%
\leq2\left\Vert u_{k}\right\Vert _{2a}\left(  \delta_{2a}\left(  k\right)
+\delta_{2a}\left(  k-j\right)  \mathbb{I}\left(  j\leq k\right)  \right)  ,
\label{Dep}%
\end{equation}
which is smaller than $4\left\Vert u_{k}\right\Vert _{2a}\delta_{2a}\left(
k-j\right)  $ when $j\leq k$. Define now the vector of martingale difference
$D_{t}=\left[  D_{1t},\ldots,D_{\overline{p}_{n}t}\right]  ^{\prime}$ with%
\[
D_{jt}=\sum_{k=t}^{\infty}\mathbf{P}_{t}\left[  u_{k}u_{k-j}\right]
\]
which converges a.s. and satisfies $\mathbb{E}\left[  D_{jt}\left\vert
\mathcal{F}_{t-1}\right.  \right]  =0$, $\max_{j}\mathbb{E}\left[  \left\vert
D_{jt}\right\vert ^{a}\right]  <\infty$, provided $\left\Vert u_{t}\right\Vert
_{2a}<\infty$ and $\sum_{k=0}^{\infty}\delta_{2a}\left(  k\right)  <\infty$.
Consider the martingale $M_{j}=M_{jn}=\sum_{t=j+1}^{n}D_{jt}$ which is an
approximation of $\widetilde{R}_{j}$. Shao (Lemma A.1, 2011b) gives under
Assumption \ref{Reg} and for any $\mathfrak{a\in}\left[  1,6a\right]  $,%
\begin{equation}
\left(  \mathbb{E}^{\frac{1}{\mathfrak{a}}}\left[  \left\vert \sum_{t=j+1}%
^{n}u_{t}u_{t-j}-M_{j}\right\vert ^{\mathfrak{a}}\right]  \right)  ^{2}\leq C.
\label{Cov2M}%
\end{equation}
We shall also use a $\mathfrak{p}$-dependent version of $D_{t}$, denoted
$D_{t}^{_{t-\mathfrak{p}+1}}$, with entries%
\begin{align}
&  D_{jt}^{_{t-\mathfrak{p}+1}}=\mathbb{E}\left[  D_{jt}\left\vert
e_{t},\ldots,e_{t-\mathfrak{p}+1}\right.  \right]  =\sum_{k=t}^{\infty
}\mathbf{P}_{t}^{\prime}\left[  u_{k}u_{k-j}\right]  ,\text{ where}%
\label{Djtp}\\
\text{ }  &  \mathbf{P}_{t}^{\prime}\left[  Z\right]  =\mathbf{P}%
_{t}^{t-\mathfrak{p}+1}\left[  Z\right]  =\mathbb{E}\left[  Z\left\vert
e_{t},\ldots,e_{t-\mathfrak{p}+1}\right.  \right]  -\mathbb{E}\left[
Z\left\vert e_{t-1},\ldots,e_{t-\mathfrak{p}+1}\right.  \right]  .\nonumber
\end{align}
Arguing as in Shao (2011b, Lemma A.2-(iii)) gives%
\begin{equation}
\left\Vert D_{jt}-D_{jt}^{_{t-\mathfrak{p}+1}}\right\Vert _{\mathfrak{a}}\leq
C\left\Vert u_{t}\right\Vert _{2\mathfrak{a}}\Theta_{2\mathfrak{a}}\left(
\mathfrak{p-}j\right)  ,\quad\text{for all }j\in\left[  1,\mathfrak{p}\right]
. \label{Shaolem52}%
\end{equation}

\subsubsection{Martingale approximation and preliminary lemmas}

An important property of $D_{t}$ and $D_{t}^{_{t-\mathfrak{p}+1}}$ is as follows.

\begin{lem}
\label{CrossDD}Suppose Assumption \ref{Kernel} and \ref{Reg} hold. Let
$K_{jp}$ be as in (\ref{kjdifcov}). Then for any $p\leq\mathfrak{p}$, $t$, and
any $s\leq t-\mathfrak{p}$, $\left\Vert \sum_{j=1}^{p}K_{jp}D_{js}%
D_{jt}^{_{t-\mathfrak{p}+1}}\right\Vert _{3a}\leq Cp^{1/2}$.
\end{lem}

\noindent\textbf{Proof of Lemma \ref{CrossDD}}. We have%
\begin{align}
&  \left\Vert \sum_{j=1}^{p}K_{jp}D_{js}D_{jt}^{t-\mathfrak{p}+1}\right\Vert
_{3a}\nonumber\\
&  \quad=\left\Vert \sum_{j=1}^{p}K_{jp}\sum_{k_{1}=0}^{\infty}\mathbf{P}%
_{s}\left[  u_{s+k_{1}}u_{s+k_{1}-j}\right]  \sum_{k_{2}=0}^{\infty}%
\mathbf{P}_{t}^{\prime}\left[  u_{t+k_{2}}u_{t+k_{2}-j}\right]  \right\Vert
_{3a}\nonumber\\
&  \quad\leq\left\Vert \sum_{j=1}^{p}K_{jp}\sum_{k_{1}=0}^{j-1}\mathbf{P}%
_{s}\left[  u_{s+k_{1}}u_{s+k_{1}-j}\right]  \sum_{k_{2}=0}^{j-1}%
\mathbf{P}_{t}^{\prime}\left[  u_{t+k_{2}}u_{t+k_{2}-j}\right]  \right\Vert
_{3a}\label{TildeDD1}\\
&  \quad+\left\Vert \sum_{j=1}^{p}K_{jp}\sum_{k_{1}=0}^{j-1}\mathbf{P}%
_{s}\left[  u_{s+k_{1}}u_{s+k_{1}-j}\right]  \sum_{k_{2}=j}^{\infty}%
\mathbf{P}_{t}^{\prime}\left[  u_{t+k_{2}}u_{t+k_{2}-j}\right]  \right\Vert
_{3a}\label{TildeDD2}\\
&  \quad+\left\Vert \sum_{j=1}^{p}K_{jp}\sum_{k_{1}=j}^{\infty}\mathbf{P}%
_{s}\left[  u_{s+k_{1}}u_{s+k_{1}-j}\right]  \sum_{k_{2}=0}^{j-1}%
\mathbf{P}_{t}^{\prime}\left[  u_{t+k_{2}}u_{t+k_{2}-j}\right]  \right\Vert
_{3a}\label{TildeDD3}\\
&  \quad+\left\Vert \sum_{j=1}^{p}K_{jp}\sum_{k_{1}=j}^{\infty}\mathbf{P}%
_{s}\left[  u_{s+k_{1}}u_{s+k_{1}-j}\right]  \sum_{k_{2}=j}^{\infty}%
\mathbf{P}_{t}^{\prime}\left[  u_{t+k_{2}}u_{t+k_{2}-j}\right]  \right\Vert
_{3a}. \label{TildeDD4}%
\end{align}
We have for (\ref{TildeDD1})%
\begin{align*}
\text{(\ref{TildeDD1})}  &  =\left\Vert \sum_{j=1}^{p}K_{jp}\sum_{k_{1}%
=0}^{p-1}\mathbb{I}\left(  k_{1}<j\right)  u_{s+k_{1}-j}\mathbf{P}_{s}\left[
u_{s+k_{1}}\right]  \sum_{k_{2}=0}^{p-1}\mathbb{I}\left(  k_{2}<j\right)
u_{t+k_{2}-j}\mathbf{P}_{t}^{\prime}\left[  u_{t+k_{2}}\right]  \right\Vert
_{3a}\\
&  =\left\Vert \sum_{k_{1}=0}^{p-1}\sum_{k_{2}=0}^{p-1}\left(  \sum
_{j=k_{1}\vee k_{2}}^{p-1}K_{jp}u_{s+k_{1}-j}u_{t+k_{2}-j}\right)
\mathbf{P}_{s}\left[  u_{s+k_{1}}\right]  \mathbf{P}_{t}^{\prime}\left[
u_{t+k_{2}}\right]  \right\Vert _{3a}\\
&  \leq\sum_{k_{1}=0}^{p-1}\sum_{k_{2}=0}^{p-1}\left\Vert \sum_{j=k_{1}\vee
k_{2}}^{p-1}K_{jp}u_{s+k_{1}-j}u_{t+k_{2}-j}\right\Vert _{6a}\delta
_{12a}\left(  k_{1}\right)  \delta_{12a}\left(  k_{2}\right)  ,
\end{align*}
using $\left\Vert \mathbf{P}_{t}^{\prime}\left[  u_{t+k_{2}}\right]
\right\Vert _{12a}\leq\left\Vert \mathbf{P}_{t}\left[  u_{t+k_{2}}\right]
\right\Vert _{12a}=\delta_{12a}\left(  k_{2}\right)  $. Now (\ref{Cov2M}) and
the Burkholder inequality give%
\begin{align*}
&  \left\Vert \sum_{j=k_{1}\vee k_{2}}^{p-1}K_{jp}u_{s+k_{1}-j}u_{t+k_{2}%
-j}\right\Vert _{6a}\leq\left\Vert \sum_{j=k_{1}\vee k_{2}}^{p-1}%
K_{jp}D_{t+k_{2}-j,t-s+k_{2}-k_{1}}\right\Vert _{6a}\\
&  \quad+\left\Vert \sum_{j=k_{1}\vee k_{2}}^{p-1}K_{jp}\left(  u_{s+k_{1}%
-j}u_{t+k_{2}-j}-D_{t+k_{2}-j,t-s+k_{2}-k_{1}}\right)  \right\Vert _{6a}\leq
Cp^{1/2}.
\end{align*}
Hence (\ref{TildeDD1}) is smaller than $Cp^{1/2}$. For (\ref{TildeDD2}), we
have since $\left\{  u_{s+k_{1}-j},j\in\left[  1,k_{1}\right]  \right\}  $ and
$\left\{  \mathbf{P}_{t}^{\prime}\left[  u_{t+k_{2}}u_{t+k_{2}-j}\right]
,j\in\left[  1,k_{1}\right]  ,k_{2}\geq0\right\}  $ are independent,%
\begin{align*}
\text{(\ref{TildeDD2})}  &  =\left\Vert \sum_{k_{1}=0}^{p-1}\sum_{k_{2}%
=0}^{\infty}\left(  \sum_{j=k_{1}}^{p-1}K_{jp}u_{s+k_{1}-j}\mathbf{P}%
_{t}^{\prime}\left[  u_{t+k_{2}+j}u_{t+k_{2}}\right]  \right)  \mathbf{P}%
_{s}\left[  u_{s+k_{1}}\right]  \right\Vert _{3a}\\
&  \leq\sum_{k_{1}=0}^{p-1}\sum_{k_{2}=0}^{\infty}\left\Vert \sum_{j=k_{1}%
}^{p-1}K_{jp}u_{s+k_{1}-j}\mathbf{P}_{t}^{\prime}\left[  u_{t+k_{2}%
+j}u_{t+k_{2}}\right]  \right\Vert _{6a}\delta_{6a}\left(  k_{1}\right)  .
\end{align*}
Let $d_{t}=\sum_{k=t}^{\infty}\mathbf{P}_{t}\left[  u_{k}\right]  $ be the
martingale difference approximation of $u_{t}$, see Wu (2007). Now, since
$\left\{  u_{s+k_{1}-j},d_{s+k_{1}-j,}j\in\left[  1,k_{1}\right]  \right\}  $
and $\left\{  \mathbf{P}_{t}^{\prime}\left[  u_{t+k_{2}}u_{t+k_{2}-j}\right]
,j\in\left[  1,k_{1}\right]  ,k_{2}\geq0\right\}  $ are independent, arguing
as in the proof of Theorem 1 in Wu (2007), (\ref{Cov2M}) and the Burkholder
inequality give%
\begin{align*}
&  \left\Vert \sum_{j=k_{1}}^{p-1}K_{jp}u_{s+k_{1}-j}\mathbf{P}_{t}^{\prime
}\left[  u_{t+k_{2}+j}u_{t+k_{2}}\right]  \right\Vert _{6a}^{2}\\
&  \text{ }\leq2\left\Vert \sum_{j=k_{1}}^{p-1}K_{jp}d_{s+k_{1}-j}%
\mathbf{P}_{t}^{\prime}\left[  u_{t+k_{2}+j}u_{t+k_{2}}\right]  \right\Vert
_{6a}^{2}+2\left\Vert \sum_{j=k_{1}}^{p-1}K_{jp}\left(  u_{s+k_{1}-j}%
-d_{t}\right)  \mathbf{P}_{t}^{\prime}\left[  u_{t+k_{2}+j}u_{t+k_{2}}\right]
\right\Vert _{6a}^{2}\\
&  \text{ }\leq C\left\Vert \sum_{j=k_{1}}^{p-1}K_{jp}d_{s+k_{1}-j}^{2}\left(
\mathbf{P}_{t}^{\prime}\left[  u_{t+k_{2}+j}u_{t+k_{2}}\right]  \right)
^{2}\right\Vert _{3a}+C\left\Vert \mathbf{P}_{t}^{\prime}\left[  u_{t+k_{2}%
+j}u_{t+k_{2}}\right]  \right\Vert _{6a}^{2}\leq Ck_{1}\delta_{6a}^{2}\left(
k_{2}\right)  .
\end{align*}
Hence Assumption \ref{Reg} gives (\ref{TildeDD2})$\leq\sum_{k_{1}=0}^{p-1}%
\sum_{k_{2}=0}^{\infty}k_{1}\delta_{6a}^{2}\left(  k_{2}\right)  \delta
_{6a}\left(  k_{1}\right)  \leq C$.

For (\ref{TildeDD3}), observe first that (\ref{Cov2M}) gives%
\begin{align*}
\text{(\ref{TildeDD3})}  &  =\left\Vert \sum_{k_{1}=0}^{\infty}\sum_{k_{2}%
=0}^{p-1}\sum_{j=1}^{p}K_{jp}\mathbb{I}\left(  j\leq k_{1}\right)
\mathbf{P}_{s}\left[  u_{s+k_{1}}u_{s+k_{1}-j}\right]  \mathbb{I}\left(
k_{2}<j\right)  \mathbf{P}_{t}^{\prime}\left[  u_{t+k_{2}}u_{t+k_{2}%
-j}\right]  \right\Vert _{3a}\\
&  \leq\sum_{k_{1}=0}^{\infty}\sum_{k_{2}=0}^{p-1}\sum_{j=k_{2}}^{p}%
\mathbb{I}\left(  j\leq k_{1}\right)  \delta_{6a}\left(  k_{1}-j\right)
\left\Vert \mathbf{P}_{t}^{\prime}\left[  u_{t+k_{2}}u_{t+k_{2}-j}\right]
\right\Vert _{6a}\\
&  \leq\left(  \sum_{k_{1}=0}^{\infty}\delta_{6a}\left(  k_{1}\right)
\right)  \times\sum_{k_{2}=0}^{p-1}\sum_{j=k_{2}}^{p}\left\Vert \mathbf{P}%
_{t}^{\prime}\left[  u_{t+k_{2}}u_{t+k_{2}-j}\right]  \right\Vert _{6a}.
\end{align*}
Since $\overline{u}_{t+k_{2}-j}^{t}$ is independent of $e_{t},\ldots
,e_{t-\mathfrak{p}+1}$ and $\mathbf{P}_{t}\left[  u_{t+k_{2}}\right]  $,%
\begin{align}
&  \left\Vert \mathbf{P}_{t}^{\prime}\left[  u_{t+k_{2}}u_{t+k_{2}-j}\right]
\right\Vert _{6a}\leq\left\Vert \underset{0}{\underbrace{\mathbb{E}\left[
\overline{u}_{t+k_{2}-j}^{t}\mathbf{P}_{t}\left[  u_{t+k_{2}}\right]
\left\vert e_{t},\ldots,e_{t-\mathfrak{p}+1}\right.  \right]  }}\right\Vert
_{6a}\nonumber\\
&  \text{ }+\left\Vert \mathbb{E}\left[  \left(  u_{t+k_{2}-j}-\overline
{u}_{t+k_{2}-j}^{t}\right)  \mathbf{P}_{t}\left[  u_{t+k_{2}}\right]
\left\vert e_{t},\ldots,e_{t-\mathfrak{p}+1}\right.  \right]  \right\Vert
_{6a}\nonumber\\
&  \text{ }\leq\left\Vert u_{t+k_{2}-j}-\overline{u}_{t+k_{2}-j}%
^{t}\right\Vert _{12a}\left\Vert \mathbf{P}_{t}\left[  u_{t+k_{2}}\right]
\right\Vert _{12a}\leq\Theta_{12a}\left(  k_{2}-j\right)  \delta_{12a}\left(
k_{2}\right)  . \label{Pprime}%
\end{align}
Substituting gives that (\ref{TildeDD3})$\leq C\sum_{k_{2}=0}^{p-1}%
\sum_{j=k_{2}}^{p}\Theta_{12a}\left(  k_{2}-j\right)  \delta_{12a}\left(
k_{2}\right)  \leq C$.

For (\ref{TildeDD4}), (\ref{Dep}) and (\ref{Pprime}) give%
\begin{align*}
\text{(\ref{TildeDD4})}  &  \leq C\sum_{j=1}^{p}\left(  \sum_{k_{1}=j}%
^{\infty}\left\Vert \mathbf{P}_{s}\left[  u_{s+k_{1}}u_{s+k_{1}-j}\right]
\right\Vert _{6a}\right)  \sum_{k_{2}=j}^{\infty}\left\Vert \mathbf{P}%
_{t}^{\prime}\left[  u_{t+k_{2}}u_{t+k_{2}-j}\right]  \right\Vert _{6a}\\
&  \leq C\sum_{j=1}^{p}\left(  \sum_{k_{1}=j}^{\infty}\delta_{6a}\left(
k_{1}-j\right)  \right)  \sum_{k_{2}=j}^{\infty}\Theta_{12a}\left(
k_{2}-j\right)  \delta_{12a}\left(  k_{2}\right)  \leq C.
\end{align*}
Hence substituting gives $\left\Vert \sum_{j=1}^{p}K_{jp}D_{js}D_{jt}%
^{t-\mathfrak{p}+1}\right\Vert _{3a}\leq Cp^{1/2}$. \hfill$\square$

\bigskip

We now define a suitable sequence of Gaussian vector. Let $2\overline{p}%
_{n}\leq\ell\leq3\overline{p}_{n}$ be an integer number. Consider a sequence
of independent centered Gaussian vectors $\eta_{t}=\left[  \eta_{1t}%
,\ldots,\eta_{\overline{p}_{n}t}\right]  ^{\prime}$ with%
\begin{equation}
\mathbb{E}\left[  \eta_{j_{1}t}\eta_{j_{2}t}\right]  =\mathbb{E}\left[
D_{j_{1}t}^{_{t-\ell+1}}D_{j_{2}t}^{_{t-\ell+1}}\right]  . \label{Vareta}%
\end{equation}
We shall also assume that $\left\{  \eta_{t}\right\}  $ and $\left\{
e_{t}\right\}  $ are independent.

\begin{lem}
\label{Sumvareta}Let $\left\{  \eta_{t}\right\}  $ be as in (\ref{Vareta}) and
suppose Assumption \ref{Reg} holds. Then for all $p\in\left[  1,\overline
{p}_{n}\right]  $ and $t,s\in\left[  1,n\right]  $,%
\begin{align*}
&  \sum_{j_{1}\neq j_{2}\in\left[  1,\overline{p}_{n}\right]  }\left\vert
\operatorname*{Cov}\left(  \eta_{j_{1}t},\eta_{j_{2}t}\right)  \right\vert
\leq C\text{ and }\sum_{j=1}^{\overline{p}_{n}}\left\vert \operatorname*{Var}%
\left(  \eta_{jt}\right)  -\sigma^{4}\right\vert \leq C,\\
&  \left\vert \sum_{j=1}^{p}\left(  1-\frac{j}{n}\right)  K_{jp}\left(
\operatorname*{Var}\left(  \eta_{jt}\right)  -\sigma^{4}\right)  \right\vert
\leq C,\\
&  \left\vert \left(  2\sum_{j=1}^{p}\left(  1-\frac{j}{n}\right)  ^{2}%
K_{jp}^{2}\operatorname*{Var}\nolimits^{2}\left(  \eta_{jt}\right)  \right)
^{1/2}-\sigma^{4}V_{\Delta}\left(  p\right)  \right\vert \leq C,\\
&  \operatorname*{Var}\left(  \frac{1}{p^{1/2}}\sum_{j=1}^{p}K_{jp}D_{js}%
\eta_{jt}\left\vert D_{s}\right.  \right)  \leq\frac{C}{p}\sum_{j=1}^{p}%
K_{jp}^{2}D_{js}^{2}.
\end{align*}

\end{lem}

\noindent\textbf{Proof of Lemma \ref{Sumvareta}.} (\ref{Cov2M}) gives for all
$j_{1}$, $j_{2}$,%
\[
\operatorname*{Cov}\left(  D_{j_{1}t},D_{j_{2}t}\right)  =\lim_{n\rightarrow
\infty}\operatorname*{Cov}\left(  \frac{\sum_{t=j_{1}+1}^{n}u_{t}u_{t-j_{1}}%
}{\left(  n-j_{1}\right)  ^{1/2}},\frac{\sum_{t=j_{2}+1}^{n}u_{t}u_{t-j_{2}}%
}{\left(  n-j_{2}\right)  ^{1/2}}\right)  =\sum_{k=-\infty}^{\infty}%
\mathbb{E}\left[  u_{0}u_{j_{1}}u_{k}u_{k+j_{2}}\right]  ,
\]
see also Lemma A.2 in Shao (2011b), provided $\sum_{k=-\infty}^{\infty
}\left\vert \mathbb{E}\left[  u_{0}u_{j_{1}}u_{k}u_{k+j_{2}}\right]
\right\vert <\infty$ as shown below. (\ref{Shaolem52}) and (\ref{Vareta}) give%
\begin{equation}
\max_{j_{1},j_{2}\in\left[  0,\overline{p}_{n}\right]  }\left\vert
\operatorname*{Cov}\left(  \eta_{j_{1}t},\eta_{j_{2}t}\right)  -\sum
_{k=-\infty}^{\infty}\mathbb{E}\left[  u_{0}u_{j_{1}}u_{k}u_{k+j_{2}}\right]
\right\vert \leq C\Theta_{12a}\left(  \overline{p}_{n}\right)  .
\label{Coveta}%
\end{equation}
Now relation between cumulants and moments in Brillinger (2001) and Theorem
\ref{XW11} gives absolute summability of the $4$th moments. Hence
$\Theta_{12a}\left(  \overline{p}_{n}\right)  =O(\overline{p}_{n}^{-6})$ gives
the first bound of the Lemma. For the second and the third bound, observe that
under the null%
\[
\left\vert \sum_{k=-\infty}^{\infty}\mathbb{E}\left[  u_{0}u_{j}u_{k}%
u_{k+j}\right]  -\sigma^{4}\right\vert \leq\left\vert \mathbb{E}\left[
u_{0}^{2}u_{j}^{2}\right]  -\mathbb{E}\left[  u_{0}^{2}\right]  \mathbb{E}%
\left[  u_{j}^{2}\right]  \right\vert +2\left\vert \sum_{k=1}^{\infty
}\mathbb{E}\left[  u_{0}u_{j}u_{k}u_{k+j}\right]  \right\vert .
\]
$\left\vert \mathbb{E}\left[  u_{0}^{2}u_{j}^{2}\right]  -\mathbb{E}\left[
u_{0}^{2}\right]  \mathbb{E}\left[  u_{j}^{2}\right]  \right\vert \leq
C\Theta_{12a}\left(  j\right)  =O\left(  j^{-6}\right)  $ and absolute
summability of the $4$th moments gives the second bound. This also gives the
fourth one since%
\begin{align*}
&  \left\vert \left(  2\sum_{j=1}^{p}\left(  1-\frac{j}{n}\right)  ^{2}%
K_{jp}^{2}\operatorname*{Var}\nolimits^{2}\left(  \eta_{jt}\right)  \right)
^{1/2}-\sigma^{4}V_{\Delta}\left(  p\right)  \right\vert \\
&  \text{ }\leq\left(  2\sum_{j=1}^{p}\left(  1-\frac{j}{n}\right)  ^{2}%
K_{jp}^{2}\left(  \operatorname*{Var}\left(  \eta_{jt}\right)  -\sigma
^{4}\right)  ^{2}\right)  ^{1/2}\\
&  \text{ }\leq2^{1/2}\left\vert \sum_{j=1}^{p}\left(  1-\frac{j}{n}\right)
K_{jp}\left(  \operatorname*{Var}\left(  \eta_{jt}\right)  -\sigma^{4}\right)
\right\vert \leq C.
\end{align*}
For the last one, observe first that%
\[
\sum_{1\leq j_{1}<j_{2}\leq\overline{p}_{n}}\left\vert \operatorname*{Cov}%
\left(  \eta_{j_{1}t},\eta_{j_{2}t}\right)  \right\vert ^{2}\leq\left(
\sum_{1\leq j_{1}<j_{2}\leq\overline{p}_{n}}\left\vert \operatorname*{Cov}%
\left(  \eta_{j_{1}t},\eta_{j_{2}t}\right)  \right\vert \right)  ^{2}<\infty
\]
by Theorem \ref{XW11} since the $2$th cumulants are the covariance. This
gives, for any $z=\left[  z_{1},\ldots,z_{\overline{p}_{n}}\right]  ^{\prime}
$,
\begin{align*}
\operatorname*{Var}\left(  z^{\prime}\eta\right)   &  =z^{\prime}%
\mathbb{E}\left[  \eta\eta^{\prime}\right]  z\leq\sum_{j=1}^{\overline{p}_{n}%
}\operatorname*{Var}\left(  \eta_{jt}\right)  z_{j}^{2}+2\sum_{1\leq
j_{1}<j_{2}\leq\overline{p}_{n}}\left\vert \operatorname*{Cov}\left(
\eta_{j_{1}t},\eta_{j_{2}t}\right)  \right\vert \left\vert z_{j_{1}%
}\right\vert \left\vert z_{j_{2}}\right\vert \\
&  \leq Czz^{\prime}+2\left(  \sum_{1\leq j_{1}<j_{2}\leq\overline{p}_{n}%
}\left\vert \operatorname*{Cov}\left(  \eta_{j_{1}t},\eta_{j_{2}t}\right)
\right\vert ^{2}\right)  ^{1/2}\left(  \sum_{1\leq j_{1}<j_{2}\leq\overline
{p}_{n}}z_{j_{1}}^{2}z_{j_{2}}^{2}\right)  ^{1/2}\\
&  \leq Cz^{\prime}z.
\end{align*}
Hence $\operatorname*{Var}\left(  \sum_{j=1}^{p}K_{jp}D_{js}\eta
_{jt}\left\vert D_{s}\right.  \right)  \leq C\left(  \sum_{j=1}^{p}K_{jp}%
^{2}D_{js}^{2}\right)  ^{1/2}$ since $\left\{  D_{t}\right\}  $ and $\left\{
\eta_{t}\right\}  $ are independent.\hfill$\square$

\subsubsection{The deviation probability of the maximum of Proposition
\ref{Selection}}

The proof is based on a smooth approximation of the maximum of real numbers
$x_{1},\ldots,x_{\overline{p}_{n}}$. Consider an increasing and three times
continuously differentiable real function $f$ with%
\begin{equation}
\lim_{x\rightarrow-\infty}f\left(  x\right)  =1,\quad f\left(  x\right)
=x\text{ for }x\geq2,\quad\max_{i=1,2,3}\sup_{x}\left\vert f^{(i)}\left(
x\right)  \right\vert <\infty. \label{f}%
\end{equation}
Let $e=e_{n}\rightarrow\infty$ with $\ln\left(  \overline{p}_{n}\right)
/e=o(1)$. Then $\max_{p\in\left[  1,\overline{p}_{n}\right]  }\left\{
f\left(  x_{p}\right)  \right\}  \leq\left(  \sum_{p=1}^{\overline{p}_{n}%
}f^{e}\left(  x_{p}\right)  \right)  ^{1/e}\leq\overline{p}_{n}^{1/e}%
\max_{p\in\left[  1,\overline{p}_{n}\right]  }\left\{  f\left(  x_{p}\right)
\right\}  $ gives that%
\begin{equation}
\left(  \sum_{p=1}^{\overline{p}_{n}}f^{e}\left(  x_{p}\right)  \right)
^{1/e}=\left(  1+O\left(  \frac{\ln\overline{p}_{n}}{e}\right)  \right)
\max_{p\in\left[  1,\overline{p}_{n}\right]  }\left\{  f\left(  x_{p}\right)
\right\}  . \label{Smthmax}%
\end{equation}
We will first find a suitable approximation for the distribution of%
\begin{equation}
\mathcal{M}=\left(  \sum_{p=1}^{\overline{p}_{n}}f^{e}\left(  \check{s}%
_{p}\right)  \right)  ^{1/e}\text{ where }\check{S}_{p}=n\sum_{j=1}^{p}%
K_{jp}\left(  \frac{M_{jn}}{n}\right)  ^{2},\quad\check{s}_{p}=\frac{\check
{S}_{p}-\sigma^{4}E_{\Delta}(p)}{\sigma^{4}V_{\Delta}(p)}\text{.}
\label{Pseudomax}%
\end{equation}
Define, for $\eta=\left[  \eta_{1},\ldots,\eta_{\overline{p}_{n}}\right]
^{\prime}$ and $x\in\left[  0,1\right]  $,%
\begin{align}
M_{jt}\left(  x;\eta\right)   &  =\sum_{s=j+1}^{t-1}D_{js}+x\eta_{j}%
+\sum_{s=t+1}^{n}\eta_{js},\quad R_{jt}\left(  x;\eta\right)  =\frac
{M_{jt}\left(  x;\eta\right)  }{n}\nonumber\\
\check{s}_{pt}\left(  x;\eta\right)   &  =\frac{n\sum_{j=1}^{p}K_{jp}%
R_{jt}^{2}\left(  x;\eta\right)  -\sigma^{4}E_{\Delta}(p)}{\sigma^{4}%
V_{\Delta}\left(  p\right)  },\quad\Sigma_{t}\left(  x;\eta\right)  =f\left(
\check{s}_{pt}\left(  x;\eta\right)  \right)  ,\nonumber\\
\mathcal{M}_{t}\left(  x;\eta\right)   &  =\left(  \sum_{p=1}^{\overline
{p}_{n}}\Sigma_{t}^{e}\left(  x;\eta\right)  \right)  ^{\frac{1}{e}}%
,\quad\mathcal{M}_{t}\left(  \eta\right)  =\mathcal{M}_{t}\left(
1;\eta\right)  , \label{Meta}%
\end{align}
and%
\begin{align*}
\check{s}_{pt}^{(1)}\left(  x;\eta\right)   &  =\frac{d\check{s}_{pt}\left(
x;\eta\right)  }{dx}=\frac{2\sum_{j=1}^{p}K_{jp}\left(  \sum_{s=j+1}%
^{t-1}D_{js}+x\eta_{j}+\sum_{s=t+1}^{n}\eta_{js}\right)  \eta_{j}}{n\sigma
^{4}V_{\Delta}\left(  p\right)  },\\
\check{s}_{pt}^{(2)}\left(  x;\eta\right)   &  =\frac{d_{pt}^{2}\check
{s}\left(  x;\eta\right)  }{dx^{2}}=\frac{2\sum_{j=1}^{p}K_{jp}\eta_{j}^{2}%
}{n\sigma^{4}V_{\Delta}\left(  p\right)  },\\
\Sigma_{pt}^{(1)}\left(  x;\eta\right)   &  =f^{(1)}\left(  \check{s}%
_{pt}\left(  x;\eta\right)  \right)  \check{s}_{pt}^{(1)}\left(
x;\eta\right)  ,\\
\Sigma_{pt}^{(2)}\left(  x;\eta\right)   &  =f^{(2)}\left(  \check{s}%
_{pt}\left(  x;\eta\right)  \right)  \left(  \check{s}_{pt}^{(1)}\left(
x;\eta\right)  \right)  ^{2}+f^{(1)}\left(  \check{s}_{pt}\left(
x;\eta\right)  \right)  \check{s}_{pt}^{(2)}\left(  x;\eta\right)  ,\\
\Sigma_{pt}^{(3)}\left(  x;\eta\right)   &  =f^{(3)}\left(  \check{s}%
_{pt}\left(  x;\eta\right)  \right)  \left(  \check{s}_{pt}^{(1)}\left(
x;\eta\right)  \right)  ^{3}+3f^{(2)}\left(  \check{s}_{pt}\left(
x;\eta\right)  \right)  \check{s}_{pt}^{(1)}\left(  x;\eta\right)  \check
{s}_{pt}^{(2)}\left(  x;\eta\right)  .
\end{align*}
We first bound the moments of $\Sigma_{pt}^{(1)}\left(  x;\eta\right)  $,
$\Sigma_{pt}^{(2)}\left(  x;\eta\right)  $ and $\Sigma_{pt}^{(3)}\left(
x;\eta\right)  $ when $\eta$ is set to $D_{t}$ or $\eta_{t}$.

\begin{lem}
Under Assumption \ref{Reg} and if $\overline{p}_{n}=O\left(  n^{1/2}\right)
$, we have uniformly in $p\in\left[  1,\overline{p}_{n}\right]  $,
$x\in\left[  0,1\right]  $ and $t=1,\ldots,n$,\label{Momsig}%
\begin{align}
\max\left\{  \left\Vert \Sigma_{pt}^{(1)}\left(  x;D_{t}\right)  \right\Vert
_{3a},\left\Vert \Sigma_{pt}^{(1)}\left(  x;\eta_{t}\right)  \right\Vert
_{3a}\right\}   &  \leq\frac{C}{n^{1/2}},\label{Sig1}\\
\max\left\{  \left\Vert \Sigma_{pt}^{(2)}\left(  x;D_{t}\right)  \right\Vert
_{3a/2},\left\Vert \Sigma_{pt}^{(2)}\left(  x;\eta_{t}\right)  \right\Vert
_{3a/2}\right\}   &  \leq\frac{Cp^{1/2}}{n},\label{Sig2}\\
\max\left\{  \left\Vert \Sigma_{pt}^{(3)}\left(  x;D_{t}\right)  \right\Vert
_{a},\left\Vert \Sigma_{pt}^{(3)}\left(  x;\eta_{t}\right)  \right\Vert
_{a}\right\}   &  \leq\frac{Cp^{1/2}}{n^{3/2}}. \label{Sig3}%
\end{align}

\end{lem}

\noindent\textbf{Proof of Lemma \ref{Momsig}.} (\ref{f}) gives%
\begin{align}
\left\vert \Sigma_{pt}^{(1)}\left(  x;\eta\right)  \right\vert  &  \leq
C\left\vert \check{s}_{pt}^{(1)}\left(  x;\eta\right)  \right\vert
,\quad\left\vert \Sigma_{pt}^{(2)}\left(  x;\eta\right)  \right\vert \leq
C\left(  \left(  \check{s}_{pt}^{(1)}\left(  x;\eta\right)  \right)
^{2}+\left\vert \check{s}_{pt}^{(2)}\left(  x;\eta\right)  \right\vert
\right)  ,\nonumber\\
\left\vert \Sigma_{pt}^{(3)}\left(  x;\eta\right)  \right\vert  &  \leq
C\left\vert \check{s}_{pt}^{(1)}\left(  x;\eta\right)  \right\vert \left(
\left(  \check{s}_{pt}^{(1)}\left(  x;\eta\right)  \right)  ^{2}+\left\vert
\check{s}_{pt}^{(2)}\left(  x;\eta\right)  \right\vert \right)  .
\label{Boundss}%
\end{align}
(\ref{Boundss}) shows that the lemma directly follows from%
\begin{align}
\max\left\{  \left\Vert \check{s}_{pt}^{(1)}\left(  x;D_{t}\right)
\right\Vert _{3a},\left\Vert \check{s}_{pt}^{(1)}\left(  x;\eta_{t}\right)
\right\Vert _{3a}\right\}   &  \leq\frac{C}{n^{1/2}},\label{sig1}\\
\max\left\{  \left\Vert \check{s}_{pt}^{(2)}\left(  x;D_{t}\right)
\right\Vert _{3a/2},\left\Vert \check{s}_{pt}^{(2)}\left(  x;\eta_{t}\right)
\right\Vert _{3a/2}\right\}   &  \leq\frac{Cp^{1/2}}{n}. \label{sig2}%
\end{align}
(\ref{sig2}) directly follow from the triangular inequality. For (\ref{sig1}),
we first bound $\left\Vert \check{s}_{pt}^{(1)}\left(  x;D_{t}\right)
\right\Vert _{3a}$. We have%
\begin{align}
&  \left\Vert \check{s}_{pt}^{(1)}\left(  x;D_{t}\right)  \right\Vert
_{3a}\leq C\left\Vert \frac{\sum_{s=1}^{t-1}\left(  \sum_{j=1}^{p}K_{jp}%
D_{js}D_{jt}\right)  }{np^{1/2}}\right\Vert _{3a}\label{Momsig1}\\
&  +C\left\Vert \frac{\sum_{j=1}^{p}K_{jp}D_{jt}^{2}}{np^{1/2}}\right\Vert
_{3a}+C\left\Vert \frac{\sum_{s=t+1}^{n}\left(  \sum_{j=1}^{p}K_{jp}D_{jt}%
\eta_{js}\right)  }{np^{1/2}}\right\Vert _{3a}. \label{Momsig2}%
\end{align}
We have, for the first item (\ref{Momsig1})%
\begin{align*}
\text{(\ref{Momsig1})}  &  \leq\left\Vert \frac{\sum_{j=1}^{p}D_{jt}\sum
_{s=1}^{t-\mathfrak{p}}K_{jp}D_{js}}{np^{1/2}}\right\Vert _{3a}+\left\Vert
\frac{\sum_{s=t-\mathfrak{p+1}}^{t-1}D_{jt}\sum_{j=1}^{p}K_{jp}D_{js}%
}{np^{1/2}}\right\Vert _{3a}\\
&  \leq\left\Vert \frac{\sum_{j=1}^{p}D_{jt}\sum_{s=1}^{t-\mathfrak{p}}%
K_{jp}D_{js}}{np^{1/2}}\right\Vert _{3a}+\frac{1}{np^{1/2}}\sum_{j=1}%
^{p}\left\Vert K_{jp}D_{jt}\right\Vert _{6a}\left\Vert \sum
_{s=t-\mathfrak{p+1}}^{t-1}D_{js}\right\Vert _{6a}\\
&  \leq\left\Vert \frac{\sum_{s=1}^{t-\mathfrak{p}}K_{jp}\sum_{j=1}^{p}%
D_{jt}D_{js}}{np^{1/2}}\right\Vert _{3a}+\frac{Cp^{1/2}\mathfrak{p}^{1/2}}{n},
\end{align*}
where $\mathfrak{p}\geq p$ and by the Burkholder inequality. Now let
$\widetilde{D}_{jt}=D_{jt}^{t-\mathfrak{p}+1}$ be as in (\ref{Djtp}). Since
$\sum_{j=1}^{p}K_{jp}D_{js}\widetilde{D}_{jt}$ is a martingale difference
given $e_{t},\ldots,e_{t-\mathfrak{p}+1}$, (\ref{Shaolem52}), the Burkholder
and triangular inequalities, Lemma \ref{CrossDD} give%
\begin{align*}
&  \left\Vert \frac{\sum_{j=1}^{p}\sum_{s=1}^{t-\mathfrak{p}}K_{jp}%
D_{js}D_{jt}}{np^{1/2}}\right\Vert _{3a}\\
&  \text{ }\leq\left\Vert \frac{\sum_{s=1}^{t-\mathfrak{p}}\sum_{j=1}%
^{p}K_{jp}D_{js}\widetilde{D}_{jt}}{np^{1/2}}\right\Vert _{3a}+\frac
{1}{np^{1/2}}\sum_{j=1}^{p}\left\vert K_{jp}\right\vert \left\Vert \sum
_{s=1}^{t-\mathfrak{p}}D_{js}\right\Vert _{6a}\left\Vert D_{jt}-\widetilde
{D}_{jt}\right\Vert _{6a}\\
&  \text{ }\leq\frac{C}{np^{1/2}}\left(  \sum_{s=1}^{t-\mathfrak{p}}\left\Vert
\sum_{j=1}^{p}K_{jp}D_{js}\widetilde{D}_{jt}\right\Vert _{3a}^{2}\right)
^{1/2}+C\frac{\Theta_{6a}\left(  \mathfrak{p}-p\right)  }{p^{1/2}}\\
&  \text{ }\leq\frac{C}{np^{1/2}}\left(  \left\vert t-\mathfrak{p}\right\vert
p\right)  ^{1/2}+C\frac{\Theta_{6a}\left(  \mathfrak{p}-p\right)  }{p^{1/2}%
}\leq C\left(  \frac{1}{n^{1/2}}+\frac{\Theta_{6a}\left(  \mathfrak{p}%
-p\right)  }{p^{1/2}}\right)  .
\end{align*}
Hence substituting gives%
\begin{equation}
\left\Vert \frac{\sum_{s=1}^{t-1}\left(  \sum_{j=1}^{p}K_{jp}D_{js}%
D_{jt}\right)  }{np^{1/2}}\right\Vert _{3a}\leq C\left(  \frac{1}{n^{1/2}%
}+\frac{p^{1/2}\mathfrak{p}^{1/2}}{n}+\frac{\Theta_{6a}\left(  \mathfrak{p}%
-p\right)  }{p^{1/2}}\right)  . \label{Boundsig1}%
\end{equation}

For the first item in (\ref{Momsig2}), (\ref{sig2}) gives a bound $C/n^{1/2}$.
For the second item in (\ref{Momsig2}), conditional Gaussianity of the
$\left\{  \sum_{j=1}^{p}K_{jp}D_{jt}\eta_{js}\right\}  $ and Lemma
\ref{Sumvareta} give%
\begin{align*}
&  \left\Vert \frac{\sum_{s=t+1}^{n}\left(  \sum_{j=1}^{p}K_{jp}D_{jt}%
\eta_{js}\right)  }{np^{1/2}}\right\Vert _{3a}\\
&  \text{ }\leq\frac{C}{np^{1/2}}\left\Vert \left\{  \sum_{s=t+1}^{n}\left(
\sum_{j=1}^{p}K_{jp}^{2}D_{jt}^{2}\right)  \right\}  ^{1/2}\right\Vert
_{3a}\leq\frac{C}{np^{1/2}}\left\Vert \sum_{s=t+1}^{n}\left(  \sum_{j=1}%
^{p}K_{jp}^{2}D_{jt}^{2}\right)  \right\Vert _{3a/2}^{1/2}\\
&  \text{ }\leq\frac{C}{np^{1/2}}\left(  \sum_{s=t+1}^{n}\sum_{j=1}^{p}%
K_{jp}^{2}\left\Vert D_{jt}\right\Vert _{3a}^{2}\right)  ^{1/2}\leq\frac
{C}{np^{1/2}}\left(  \left(  n-t\right)  p\right)  ^{1/2}\leq\frac{C}{n^{1/2}%
}.
\end{align*}
Substituting the two last bounds in (\ref{Momsig2}) and (\ref{Boundsig1}) in
(\ref{Momsig1}) shows that
\begin{equation}
\max\left\{  \left\Vert \check{s}_{pt}^{(1)}\left(  x;D_{t}\right)
\right\Vert _{3a},\left\Vert \check{s}_{pt}^{(1)}\left(  x;\eta_{t}\right)
\right\Vert _{3a}\right\}  \leq C\left(  \frac{1}{n^{1/2}}+\frac
{p^{1/2}\mathfrak{p}^{1/2}}{n}+\frac{\Theta_{6a}\left(  \mathfrak{p}-p\right)
}{p^{1/2}}\right)  . \label{sig1a}%
\end{equation}
Observe that $\Theta_{6a}\left(  \mathfrak{p}-p\right)  \leq C\left(
\mathfrak{p}-p\right)  ^{-11/2}$ by Assumption \ref{Reg}. Consider now%
\[
\mathfrak{p}=\max\left(  2p,\left(  \frac{n}{p}\right)  ^{\frac{1}{6}}\right)
\geq2p,
\]
which is such that, since $p\in\left[  1,\overline{p}_{n}\right]  $ with
$\overline{p}_{n}=O\left(  n^{1/2}\right)  $,%
\begin{align*}
\text{If }\left(  \frac{n}{p}\right)  ^{\frac{1}{6}}  &  \geq2p\text{, }%
\frac{\left(  \mathfrak{p}-p\right)  ^{-11/2}}{p^{1/2}}\asymp\frac
{p^{1/2}\mathfrak{p}^{1/2}}{n}\leq\frac{\mathfrak{p}}{n}\leq\frac{1}{n^{5/6}%
}\leq\frac{1}{n^{1/2}},\\
\text{If }\left(  \frac{n}{p}\right)  ^{\frac{1}{6}}  &  <2p\Leftrightarrow
\left(  \frac{n}{2^{6}}\right)  ^{\frac{1}{7}}<p\text{, }\frac{\Theta
_{6a}\left(  \mathfrak{p}-p\right)  }{p^{1/2}}\leq Cp^{-6}\leq\frac{C}%
{n^{1/2}}\text{, }\frac{p^{1/2}\mathfrak{p}^{1/2}}{n}\leq\frac{\overline
{p}_{n}}{n}\leq\frac{C}{n^{1/2}}.
\end{align*}
Hence (\ref{sig1a}) gives (\ref{sig1}).\hfill$\square$

\bigskip

Let $I\left(  \cdot\right)  $ be a three times differentiable real function
and define for $\mathcal{M}_{t}\left(  \eta\right)  $ as in (\ref{Meta}),%
\[
\mathcal{I}_{t}\left(  \eta\right)  =\mathcal{I}_{tn}\left(  \eta\right)
=I\left(  \mathcal{M}_{t}\left(  \eta\right)  \right)  ,\quad\mathcal{I}%
_{t}\left(  x;\eta\right)  =\mathcal{I}\left(  x\eta\right)  ,\quad
\mathcal{I}_{t}^{(j)}\left(  x;\eta\right)  =\frac{d_{t}^{j}\mathcal{I}\left(
x;\eta\right)  }{d^{j}x},\quad j=1,2.
\]
Observe that $I\left(  \mathcal{M}\right)  =I\left(  \mathcal{M}_{n}\left(
D_{n}\right)  \right)  =\mathcal{I}_{n}\left(  D_{n}\right)  $, $\mathcal{I}%
_{t}\left(  D_{t}\right)  =\mathcal{I}_{t+1}\left(  \eta_{t+1}\right)  $, and
that $I\left(  \mathcal{M}_{1}\left(  \eta_{1}\right)  \right)  $
$=\mathcal{I}_{1}\left(  \eta_{1}\right)  $ is a function of the Gaussian
vectors $\eta_{1},\ldots,\eta_{n}$ only.

\begin{lem}
Let $\mathcal{M}$ and $\mathcal{M}_{1}\left(  \eta_{1}\right)  $ be as in
(\ref{Pseudomax}) and (\ref{Meta}). Consider a real function $I\left(
\cdot\right)  $ which may depend on $n$ and three times continuously
differentiable with $\max_{j=1,2,3}\sup_{x}\left\vert I^{(j)}\left(  x\right)
\right\vert \leq C$. Then under Assumptions \ref{P}, \ref{Reg} and if
$e=O\left(  \overline{p}_{n}^{1/(2a)}\right)  $,\label{Lindeberg}%
\[
\left\vert \mathbb{E}\left[  I\left(  \mathcal{M}\right)  -I\left(
\mathcal{M}_{1}\left(  \eta_{1}\right)  \right)  \right]  \right\vert \leq
C\left(  \frac{\overline{p}_{n}^{1+3/a}}{n^{1/2}}+\frac{1}{\overline{p}%
_{n}^{1-1/a}}\right)  .
\]

\end{lem}

\noindent\textbf{Proof of Lemma \ref{Lindeberg}. }The proof of the Lemma works
by changing $D_{n}$ into $\eta_{n}$, $D_{n-1}$ into $\eta_{n-1}$ and so on,
the so called Lindeberg technique described in Pollard (2002,
p.179).\textbf{\ }This amounts to decompose $I\left(  \mathcal{M}\right)
-I\left(  \mathcal{M}_{n}\left(  \eta_{n}\right)  \right)  $ into the
following sum of differences,
\begin{align*}
&  I\left(  \mathcal{M}\right)  -I\left(  \mathcal{M}_{n}\left(  \eta
_{n}\right)  \right) \\
&  \text{ }=\mathcal{I}_{n}\left(  D_{n}\right)  -\mathcal{I}_{n-1}\left(
D_{n-1}\right)  +\mathcal{I}_{n-1}\left(  D_{n-1}\right)  -\mathcal{I}%
_{n-2}\left(  D_{n-2}\right)  +\cdots+\mathcal{I}_{1}\left(  D_{1}\right)
-\mathcal{I}_{1}\left(  \eta_{1}\right) \\
&  \text{ }=\mathcal{I}_{n}\left(  D_{n}\right)  -\mathcal{I}_{n}\left(
\eta_{n}\right)  +\mathcal{I}_{n-1}\left(  D_{n-1}\right)  -\mathcal{I}%
_{n-1}\left(  \eta_{n-1}\right)  +\cdots+\mathcal{I}_{1}\left(  D_{1}\right)
-\mathcal{I}_{1}\left(  \eta_{1}\right)  .
\end{align*}
Since $\mathcal{I}_{t}(\eta)=\mathcal{I}_{t}(1;\eta)$ and $\mathcal{I}%
_{t}(0;\eta)=\mathcal{I}_{t}(0)$, a third-order Taylor expansion around
$\eta=0$ with integral remainder gives%
\begin{align*}
&  \left[  \mathcal{I}_{t}(D_{t})-\mathcal{I}_{t}(\eta_{t})\right]
=\mathbb{E}\left[  \mathcal{I}_{t}^{(1)}(0;D_{t})-\mathcal{I}_{t}^{(1)}%
(0;\eta_{t})\right] \\
&  +\frac{1}{2}\mathbb{E}\left[  \mathcal{I}_{t}^{(2)}(0;D_{t})-\mathcal{I}%
_{t}^{(2)}(0;\eta_{t})\right]  +\frac{1}{2}\int_{0}^{1}(1-x)^{2}%
\mathbb{E}\left[  \mathcal{I}_{t}^{(3)}(x;D_{t})-\mathcal{I}_{t}^{(3)}%
(x;\eta_{t})\right]  dx.
\end{align*}
Since $\left\{  D_{t}\right\}  $ is a sequence of martingale difference,
$\mathbb{E}\left[  \mathcal{I}_{t}^{(1)}(0;D_{t})-\mathcal{I}_{t}^{(1)}%
(0;\eta_{t})\right]  =0$ due to the expression of $\mathcal{I}_{t}%
^{(1)}\left(  0;\eta\right)  $ given above. Hence%
\begin{align}
&  \left\vert \mathbb{E}\left[  I\left(  \mathcal{M}\right)  \right]
-\mathbb{E}\left[  I\left(  \mathcal{M}_{1}\left(  \eta_{1}\right)  \right)
\right]  \right\vert \leq\frac{1}{2}\left\vert \sum_{t=1}^{n}\mathbb{E}\left[
\mathcal{I}_{t}^{(2)}(0;D_{t})-\mathcal{I}_{t}^{(2)}(0;\eta_{t})\right]
\right\vert \label{Lind2}\\
&  +\frac{1}{2}\int_{0}^{1}(1-x)^{2}\left\{  \sum_{t=1}^{n}\left\vert
\mathbb{E}\left[  \mathcal{I}_{t}^{(3)}(x;D_{t})-\mathcal{I}_{t}^{(3)}%
(x;\eta_{t})\right]  \right\vert \right\}  dx. \label{Lind3}%
\end{align}

We now compute the differentials $\mathcal{I}_{t}^{(j)}\left(  x;\eta\right)
$, $j=1,2,3$. We have%
\begin{align*}
\mathcal{I}_{t}^{(1)}\left(  x;\eta\right)   &  =I^{\prime}\left(
\mathcal{M}_{t}\left(  x;\eta\right)  \right)  \mathcal{M}_{t}^{(1)}\left(
x;\eta\right)  ,\\
\mathcal{I}_{t}^{(2)}\left(  x;\eta\right)   &  =I^{^{\prime\prime}}\left(
\mathcal{M}_{t}\left(  x;\eta\right)  \right)  \left(  \mathcal{M}_{t}%
^{(1)}\left(  x;\eta\right)  \right)  ^{2}+I^{\prime}\left(  \mathcal{M}%
_{t}\left(  x;\eta\right)  \right)  \mathcal{M}_{t}^{(2)}\left(
x;\eta\right)  ,\\
\mathcal{I}_{t}^{(3)}\left(  x;\eta\right)   &  =I^{^{\prime\prime\prime}%
}\left(  \mathcal{M}_{t}\left(  x;\eta\right)  \right)  \left(  \mathcal{M}%
_{t}^{(1)}\left(  x;\eta\right)  \right)  ^{3}+3I^{^{\prime\prime}}\left(
\mathcal{M}_{t}\left(  x;\eta\right)  \right)  \mathcal{M}_{t}^{(1)}\left(
x;\eta\right)  \mathcal{M}_{t}^{(2)}\left(  x;\eta\right) \\
&  +I^{\prime}\left(  \mathcal{M}_{t}\left(  x;\eta\right)  \right)
\mathcal{M}_{t}^{(3)}\left(  x;\eta\right)  .
\end{align*}
We compute the differentials of $\mathcal{M}_{t}$. We have%
\begin{align*}
\mathcal{M}_{t}^{(1)}\left(  x;\eta\right)   &  =\left(  \sum_{p=1}%
^{\overline{p}_{n}}\Sigma_{pt}^{e}\left(  x;\eta\right)  \right)  ^{1/e-1}%
\sum_{p=1}^{\overline{p}_{n}}\Sigma_{pt}^{e-1}\left(  x;\eta\right)
\Sigma_{pt}^{(1)}\left(  x;\eta\right) \\
&  =\mathcal{M}_{t}^{1-e}\left(  x;\eta\right)  \sum_{p=1}^{\overline{p}_{n}%
}\Sigma_{pt}^{e-1}\left(  x;\eta\right)  \Sigma_{pt}^{(1)}\left(
x;\eta\right)  ,\\
\mathcal{M}_{t}^{(2)}\left(  x;\eta\right)   &  =\mathcal{M}_{1t}^{(2)}\left(
x;\eta\right)  +\mathcal{M}_{2t}^{(2)}\left(  x;\eta\right)  +\mathcal{M}%
_{3t}^{(2)}\left(  x;\eta\right)  ,\\
\mathcal{M}^{(3)}\left(  x;\eta\right)   &  =\mathcal{M}_{1t}^{(3)}\left(
x;\eta\right)  +\cdots+\mathcal{M}_{6t}^{(3)}\left(  x;\eta\right)  ,
\end{align*}
where, dropping the variables $x$, $\eta$ for notational convenience%
\begin{align*}
\mathcal{M}_{1t}^{(2)}  &  =\left(  \frac{1}{e}-1\right)  \mathcal{M}%
_{t}^{1-2e}\left(  \sum_{p=1}^{\overline{p}_{n}}\Sigma_{pt}^{e-1}\Sigma
_{pt}^{(1)}\right)  ^{2},\\
\mathcal{M}_{2t}^{(2)}  &  =\mathcal{M}_{t}^{1-e}\sum_{p=1}^{\overline{p}_{n}%
}\Sigma_{pt}^{e-1}\Sigma_{pt}^{(2)},\\
\mathcal{M}_{3t}^{(2)}  &  =\left(  e-1\right)  \mathcal{M}_{t}^{1-e}%
\sum_{p=1}^{\overline{p}_{n}}\Sigma_{pt}^{e-2}\left(  \Sigma_{pt}%
^{(1)}\right)  ^{2},\\
\mathcal{M}_{1t}^{(3)}  &  =\left(  \frac{1}{e}-1\right)  \left(  \frac{1}%
{e}-2\right)  \mathcal{M}_{t}^{1-3e}\left(  \sum_{p=1}^{\overline{p}_{n}%
}\Sigma_{pt}^{e-1}\Sigma_{pt}^{(1)}\right)  ^{3},\\
\mathcal{M}_{2t}^{(3)}  &  =3\left(  \frac{1}{e}-1\right)  \mathcal{M}%
_{t}^{1-2e}\sum_{p=1}^{\overline{p}_{n}}\Sigma_{pt}^{e-1}\Sigma_{pt}^{(1)}%
\sum_{p=1}^{\overline{p}_{n}}\Sigma_{pt}^{e-1}\Sigma_{pt}^{(2)},\\
\mathcal{M}_{3t}^{(3)}  &  =3\left(  \frac{1}{e}-1\right)  \left(  e-1\right)
\mathcal{M}_{t}^{1-2e}\sum_{p=1}^{\overline{p}_{n}}\Sigma_{pt}^{e-1}%
\Sigma_{pt}^{(1)}\sum_{p=1}^{\overline{p}_{n}}\Sigma_{pt}^{e-2}\left(
\Sigma_{pt}^{(1)}\right)  ^{2},\\
\mathcal{M}_{4t}^{(3)}  &  =\left(  3e-1\right)  \mathcal{M}_{t}^{1-e}%
\sum_{p=1}^{\overline{p}_{n}}\Sigma_{pt}^{e-2}\Sigma_{pt}^{(2)}\Sigma
_{pt}^{(1)},\\
\mathcal{M}_{5t}^{(3)}  &  =\left(  e-1\right)  \left(  e-2\right)
\mathcal{M}_{t}^{1-e}\sum_{p=1}^{\overline{p}_{n}}\Sigma_{pt}^{e-2}\left(
\Sigma_{pt}^{(1)}\right)  ^{3},\\
\mathcal{M}_{6t}^{(3)}  &  =\mathcal{M}_{t}^{1-e}\sum_{p=1}^{\overline{p}_{n}%
}\Sigma_{pt}^{e-1}\Sigma_{pt}^{(3)}.
\end{align*}
\medskip\textbf{The third-order item(\ref{Lind3})}. Since%
\begin{align*}
&  \frac{1}{2}\int_{0}^{1}(1-x)^{2}\left\{  \sum_{t=1}^{n}\left\vert
\mathbb{E}\left[  \mathcal{I}_{t}^{(3)}(x;D_{t})-\mathcal{I}_{t}^{(3)}%
(x;\eta_{t})\right]  \right\vert \right\}  dx\\
&  \text{ }\leq\frac{1}{2}\int_{0}^{1}(1-x)^{2}\left\{  \sum_{t=1}^{n}\left(
\left\vert \mathbb{E}\left[  \mathcal{I}_{t}^{(3)}(x;D_{t})\right]
\right\vert +\left\vert \mathbb{E}\left[  \mathcal{I}_{t}^{(3)}(x;\eta
_{t})\right]  \right\vert \right)  \right\}  dx,
\end{align*}
it is sufficient to bound $\sum_{t=1}^{n}\left\vert \mathbb{E}\left[
\mathcal{I}_{t}^{(3)}(x)\right]  \right\vert $ independently of $x$ where
$\mathcal{I}_{t}^{(3)}(x)$ stands for $\mathcal{I}_{t}^{(3)}(x;\eta_{t})$ or
$\mathcal{I}_{t}^{(3)}(x;D_{t})$. We have, dropping dependence w.r.t. to $x$
for ease of notation,%
\begin{align*}
\sum_{t=1}^{n}\left\vert \mathbb{E}\left[  \mathcal{I}_{t}^{(3)}\right]
\right\vert  &  \leq C\sum_{t=1}^{n}\left\{  \mathbb{E}\left[  \left\vert
\mathcal{M}_{t}^{(1)}\right\vert ^{3}\right]  +\mathbb{E}\left[  \left\vert
\mathcal{M}_{t}^{(1)}\mathcal{M}_{1t}^{(2)}\right\vert \right]  +\mathbb{E}%
\left[  \left\vert \mathcal{M}_{t}^{(1)}\mathcal{M}_{2t}^{(2)}\right\vert
\right]  \right\} \\
&  +C\sum_{t=1}^{n}\left\{  \mathbb{E}\left[  \left\vert \mathcal{M}_{t}%
^{(1)}\mathcal{M}_{3t}^{(2)}\right\vert \right]  +\sum_{j=1}^{6}%
\mathbb{E}\left[  \left\vert \mathcal{M}_{jt}^{(3)}\right\vert \right]
\right\}  .
\end{align*}
We now study the ten items above.

\bigskip

\textbf{(1)} $\sum_{t=1}^{n}\mathbb{E}\left[  \left\vert \mathcal{M}_{t}%
^{(1)}\right\vert ^{3}\right]  $. We have for $a$, $\overline{a}\geq1$ with
$1/a=1-1/\overline{a}$,%
\begin{align*}
\mathbb{E}\left[  \left\vert \mathcal{M}_{t}^{(1)}\right\vert ^{3}\right]   &
=\mathbb{E}\left[  \left\vert \mathcal{M}_{t}^{1-e}\sum_{p=1}^{\overline
{p}_{n}}\Sigma_{pt}^{e-1}\Sigma_{pt}^{(1)}\right\vert ^{3}\right] \\
&  \leq\sum_{p_{1},p_{2},p_{3}=1}^{\overline{p}_{n}}\mathbb{E}\left[
\left\vert \mathcal{M}_{t}^{3(1-e)}\Sigma_{p_{1}t}^{e-1}\Sigma_{p_{2}t}%
^{e-1}\Sigma_{p_{3}t}^{e-1}\Sigma_{p_{1}t}^{(1)}\Sigma_{p_{2}t}^{(1)}%
\Sigma_{p_{3}t}^{(1)}\right\vert \right] \\
&  \leq\max_{p,t}\left\Vert \Sigma_{pt}^{(1)}\right\Vert _{3a}^{3}\sum
_{p_{1},p_{2},p_{3}=1}^{\overline{p}_{n}}\mathbb{E}^{1/\overline{a}}\left[
\left\vert \mathcal{M}_{t}^{3(1-e)}\Sigma_{p_{1}t}^{e-1}\Sigma_{p_{2}t}%
^{e-1}\Sigma_{p_{3}t}^{e-1}\right\vert ^{\overline{a}}\right] \\
&  \leq\frac{C}{n^{3/2}}\sum_{p_{1},p_{2},p_{3}=1}^{\overline{p}_{n}%
}\mathbb{E}^{1/\overline{a}}\left[  \left\vert \mathcal{M}_{t}^{3(1-e)}%
\Sigma_{p_{1}t}^{e-1}\Sigma_{p_{2}t}^{e-1}\Sigma_{p_{3}t}^{e-1}\right\vert
^{\overline{a}}\right]  ,
\end{align*}
by (\ref{Sig1}) for all $x\in\left[  0,1\right]  $. Now, since $t\mapsto
t^{1/\overline{a}}$, $t\mapsto t^{1-1/e}$ are concave and $\sum_{p=1}%
^{\overline{p}_{n}}t_{p}^{\overline{a}}\leq\left(  \sum_{p=1}^{\overline
{p}_{n}}t_{p}\right)  ^{\overline{a}}$, the definition of $\mathcal{M}_{t}$
gives%
\begin{align*}
&  \sum_{p_{1},p_{2},p_{3}=1}^{\overline{p}_{n}}\mathbb{E}^{1/\overline{a}%
}\left[  \left\vert \mathcal{M}_{t}^{3(1-e)}\Sigma_{p_{1}t}^{e-1}\Sigma
_{p_{2}t}^{e-1}\Sigma_{p_{3}t}^{e-1}\right\vert ^{\overline{a}}\right] \\
&  \text{ }=\overline{p}_{n}^{3}\times\frac{1}{\overline{p}_{n}^{3}}%
\sum_{p_{1},p_{2},p_{3}=1}^{\overline{p}_{n}}\mathbb{E}^{1/\overline{a}%
}\left[  \left\vert \mathcal{M}_{t}^{3(1-e)}\Sigma_{p_{1}t}^{e-1}\Sigma
_{p_{2}t}^{e-1}\Sigma_{p_{3}t}^{e-1}\right\vert ^{\overline{a}}\right] \\
&  \text{ }\leq\overline{p}_{n}^{3}\left(  \frac{1}{\overline{p}_{n}^{3}%
}\mathbb{E}\left[  \sum_{p_{1},p_{2},p_{3}=1}^{\overline{p}_{n}}%
\mathcal{M}_{t}^{3\overline{a}(1-e)}\Sigma_{p_{1}t}^{\overline{a}%
e(1-1/e)}\Sigma_{p_{2}t}^{\overline{a}e(1-1/e)}\Sigma_{p_{3}t}^{\overline
{a}e(1-1/e)}\right]  \right)  ^{1/\overline{a}}\\
&  \text{ }=\overline{p}_{n}^{3}\left(  \mathbb{E}\left[  \left(  \sum
_{p=1}^{\overline{p}_{n}}\Sigma_{pt}^{e}\right)  ^{-3\overline{a}%
(1-1/e)}\left(  \frac{1}{\overline{p}_{n}}\sum_{p=1}^{\overline{p}_{n}}%
\Sigma_{pt}^{\overline{a}e(1-1/e)}\right)  ^{3}\right]  \right)
^{1/\overline{a}}\\
&  \text{ }\leq\overline{p}_{n}^{3}\left(  \mathbb{E}\left[  \left(
\sum_{p=1}^{\overline{p}_{n}}\Sigma_{pt}^{\overline{a}e}\right)
^{-3(1-1/e)}\left(  \frac{1}{\overline{p}_{n}}\sum_{p=1}^{\overline{p}_{n}%
}\Sigma_{pt}^{\overline{a}e}\right)  ^{3(1-1/e)}\right]  \right)
^{1/\overline{a}}\\
&  \text{ }\leq\overline{p}_{n}^{3\left(  1-1/\overline{a}\right)
+3/(e\overline{a})}\leq C\overline{p}_{n}^{3/a},
\end{align*}
uniformly w.r.t. to $t$ since $\left(  \ln\overline{p}_{n}\right)  /e=o(1)$.
Hence for all $x\in\left[  0,1\right]  $
\begin{equation}
\sum_{t=1}^{n}\mathbb{E}\left[  \left\vert \mathcal{M}_{t}^{(1)}\right\vert
^{3}\right]  \leq C\frac{\overline{p}_{n}^{3/a}}{n^{1/2}}. \label{M13}%
\end{equation}

\medskip

\textbf{(2)} $\sum_{t=1}^{n}\mathbb{E}\left[  \left\vert \mathcal{M}_{t}%
^{(1)}\mathcal{M}_{1t}^{(2)}\right\vert \right]  $. We have, since
$\mathcal{M}_{t}\geq1$,%
\begin{align*}
\mathbb{E}\left[  \left\vert \mathcal{M}_{t}^{(1)}\right\vert \left\vert
\mathcal{M}_{1t}^{(2)}\right\vert \right]   &  \leq C\mathbb{E}\left[
\mathcal{M}_{t}^{2-3e}\left\vert \sum_{p=2}^{\overline{p}_{n}}\Sigma
_{pt}^{e-1}\Sigma_{pt}^{(1)}\right\vert ^{3}\right]  \leq C\mathbb{E}\left[
\mathcal{M}_{t}^{3-3e}\left\vert \sum_{p=2}^{\overline{p}_{n}}\Sigma
_{pt}^{e-1}\Sigma_{pt}^{(1)}\right\vert ^{3}\right] \\
&  \leq C\mathbb{E}\left[  \left\vert \mathcal{M}_{t}^{(1)}\right\vert
^{3}\right]  ,
\end{align*}
for all $t$, such that $\sum_{t=1}^{n}\mathbb{E}\left[  \left\vert
\mathcal{M}_{t}^{(1)}\right\vert ^{2}\left\vert \mathcal{M}_{1t}%
^{(2)}\right\vert \right]  \leq C\sum_{t=1}^{n}\mathbb{E}\left[  \left\vert
\mathcal{M}_{t}^{(1)}\right\vert ^{3}\right]  $. Hence a bound similar to
(\ref{M13}) holds.

\medskip

\textbf{(3)} $\sum_{t=1}^{n}\mathbb{E}\left[  \left\vert \mathcal{M}_{t}%
^{(1)}\mathcal{M}_{2t}^{(2)}\right\vert \right]  $. Let $\overline{a}>1$ be
such that $1/\overline{a}=1-1/a$. Arguing as for \textbf{(1)} with
(\ref{Sig1}) and (\ref{Sig2}),
\begin{align*}
\mathbb{E}\left[  \left\vert \mathcal{M}_{t}^{(1)}\mathcal{M}_{1t}%
^{(2)}\right\vert \right]   &  \leq C\sum_{p_{1},p_{2}=1}^{\overline{p}_{n}%
}\mathbb{E}\left[  \mathcal{M}_{t}^{2(1-e)}\left\vert \Sigma_{p_{1}t}%
^{e-1}\Sigma_{p_{2}t}^{e-1}\Sigma_{p_{1}t}^{(1)}\Sigma_{p_{2}t}^{(2)}%
\right\vert \right] \\
&  \leq C\max_{p,t}\left\{  \left\Vert \Sigma_{pt}^{(1)}\right\Vert
_{3a}\left\Vert \Sigma_{pt}^{(2)}\right\Vert _{3a/2}\right\}  \sum
_{p_{1},p_{2}=1}^{\overline{p}_{n}}\mathbb{E}^{1/\overline{a}}\left[
\left\vert \mathcal{M}_{t}^{2(1-e)}\Sigma_{p_{1}t}^{e-1}\Sigma_{p_{2}t}%
^{e-1}\right\vert ^{\overline{a}}\right] \\
&  \leq C\frac{\overline{p}_{n}^{1/2}}{n^{3/2}}\times\overline{p}_{n}%
^{2}\times\mathbb{E}^{1/\overline{a}}\left[  \left(  \sum_{p=1}^{\overline
{p}_{n}}\Sigma_{pt}^{e}\right)  ^{-2\overline{a}(1-1/e)}\left(  \frac
{1}{\overline{p}_{n}}\sum_{p=1}^{\overline{p}_{n}}\Sigma_{pt}^{e\overline
{a}(1-1/e)}\right)  ^{2}\right] \\
&  \leq C\frac{\overline{p}_{n}^{1/2}}{n^{3/2}}\times\overline{p}_{n}%
^{2}\times\mathbb{E}^{1/\overline{a}}\left[  \left(  \sum_{p=1}^{\overline
{p}_{n}}\Sigma_{pt}^{e\overline{a}(1-1/e)}\right)  ^{-2}\left(  \frac
{1}{\overline{p}_{n}}\sum_{p=1}^{\overline{p}_{n}}\Sigma_{pt}^{e\overline
{a}(1-1/e)}\right)  ^{2}\right] \\
&  =C\frac{\overline{p}_{n}^{1/2}}{n^{3/2}}\times\overline{p}_{n}^{2}%
\times\overline{p}_{n}^{-2/\overline{a}}=C\frac{\overline{p}_{n}^{\frac{1}%
{2}\left(  1+4/a\right)  }}{n^{3/2}}.
\end{align*}
Hence, uniformly w.r.t. $x\in\left[  0,1\right]  $,%
\begin{equation}
\sum_{t=1}^{n}\mathbb{E}\left[  \left\vert \mathcal{M}_{t}^{(1)}%
\mathcal{M}_{2t}^{(2)}\right\vert \right]  \leq C\frac{\overline{p}_{n}%
^{\frac{1}{2}\left(  1+4/a\right)  }}{n^{1/2}}. \label{M1M22}%
\end{equation}

\medskip

\textbf{(4)} $\sum_{t=1}^{n}\mathbb{E}\left[  \left\vert \mathcal{M}_{t}%
^{(1)}\mathcal{M}_{3t}^{(2)}\right\vert \right]  $. Proceeding as \textbf{(1)}
and \textbf{(3)} gives, since $\inf_{p,t}\Sigma_{pt}\geq1$,%
\[
\mathbb{E}\left[  \left\vert \mathcal{M}_{t}^{(1)}\mathcal{M}_{3t}%
^{(2)}\right\vert \right]  \leq Ce\sum_{p_{1},p_{2}=1}^{\overline{p}_{n}%
}\mathbb{E}\left[  \mathcal{M}_{t}^{2(1-e)}\left\vert \Sigma_{p_{1}t}%
^{e-1}\Sigma_{p_{2}t}^{e-1}\Sigma_{p_{1}t}^{(1)}\left(  \Sigma_{p_{2}t}%
^{(1)}\right)  ^{2}\right\vert \right]  \leq C\frac{e\overline{p}_{n}^{2/a}%
}{n^{3/2}}\leq C\frac{\overline{p}_{n}^{3/a}}{n^{3/2}},
\]
provided $e=O(\overline{p}_{n}^{1/a})$. Hence $\sum_{t=1}^{n}\mathbb{E}\left[
\left\vert \mathcal{M}_{t}^{(1)}\mathcal{M}_{3t}^{(2)}\right\vert \right]  $
can be bounded as in (\ref{M13}).

\medskip

\textbf{(5)} $\sum_{t=1}^{n}\mathbb{E}\left[  \left\vert \mathcal{M}%
_{1t}^{(3)}\right\vert \right]  $ can be bounded as in (\ref{M13}) since
$\mathcal{M}_{t}\geq1$ gives $\mathbb{E}\left[  \left\vert \mathcal{M}%
_{1t}^{(3)}\right\vert \right]  \leq C\mathbb{E}\left[  \mathcal{M}%
_{t}^{3(1-e)}\left\vert \sum_{p=2}^{\overline{p}_{n}}\Sigma_{pt}^{e-1}%
\Sigma_{pt}^{(1)}\right\vert ^{3}\right]  .$

\medskip

\textbf{(6)} $\sum_{t=1}^{n}\mathbb{E}\left[  \left\vert \mathcal{M}%
_{2t}^{(3)}\right\vert \right]  $. Arguing as in \textbf{(3)} gives that
$\sum_{t=1}^{n}\mathbb{E}\left[  \left\vert \mathcal{M}_{2t}^{(3)}\right\vert
\right]  $ can be bounded as in (\ref{M1M22}).

\medskip

\textbf{(7)} $\sum_{t=1}^{n}\mathbb{E}\left[  \left\vert \mathcal{M}%
_{3t}^{(3)}\right\vert \right]  $. Arguing as in \textbf{(4) }shows that this
item is negligible compared to (\ref{M13}).

\medskip

\textbf{(8)} $\sum_{t=1}^{n}\mathbb{E}\left[  \left\vert \mathcal{M}%
_{4t}^{(3)}\right\vert \right]  $. \ Let $\overline{a}>1$ be such that
$1/\overline{a}=1-1/a$. We have, since $\inf_{p,t}\Sigma_{pt}\geq1$,
\begin{align*}
\mathbb{E}\left[  \left\vert \mathcal{M}_{4t}^{(3)}\right\vert \right]   &
\leq Ce\mathbb{E}\left[  \mathcal{M}_{t}^{1-e}\sum_{p=1}^{\overline{p}_{n}%
}\left\vert \Sigma_{pt}^{e-2}\Sigma_{pt}^{(2)}\Sigma_{pt}^{(1)}\right\vert
\right]  \leq Ce\sum_{p=p_{o}}^{\overline{p}_{n}}\mathbb{E}^{1/\overline{a}%
}\left[  \left(  \mathcal{M}_{t}^{1-e}\Sigma_{pt}^{e-1}\right)  ^{\overline
{a}}\right]  \left\Vert \Sigma_{pt}^{(2)}\right\Vert _{3a/2}\left\Vert
\Sigma_{pt}^{(1)}\right\Vert _{3a}\\
&  \leq C\frac{e\overline{p}_{n}^{1/2}\overline{p}_{n}^{1-1/\overline{a}}%
}{n^{3/2}}\leq C\frac{\overline{p}_{n}^{\frac{1}{2}\left(  1+4/a\right)  }%
}{n^{3/2}},
\end{align*}
provided $e=O\left(  \overline{p}_{n}^{1/a}\right)  $. This gives a bound
similar to (\ref{M1M22}) for $\sum_{t=1}^{n}\mathbb{E}\left[  \left\vert
\mathcal{M}_{4t}^{(3)}\right\vert \right]  $.

\medskip

\textbf{(9)} $\sum_{t=1}^{n}\mathbb{E}\left[  \left\vert \mathcal{M}%
_{5t}^{(3)}\right\vert \right]  $ can be bounded as in (\ref{M13}) provided
$e=O(\overline{p}_{n}^{1/(2a)})$.

\medskip

\textbf{(10)} $\sum_{t=1}^{n}\mathbb{E}\left[  \left\vert \mathcal{M}%
_{6t}^{(3)}\right\vert \right]  $ can be bounded as in (\ref{M1M22}).

\medskip

\noindent Hence, collecting the dominant bounds (\ref{M13}) and (\ref{M1M22})
in \textbf{(1)}-\textbf{(10)} gives%
\begin{equation}
\frac{1}{2}\int_{0}^{1}(1-x)^{2}\left\{  \sum_{t=1}^{n}\left\vert
\mathbb{E}\left[  \mathcal{I}_{t}^{(3)}(x;D_{t})-\mathcal{I}_{t}^{(3)}%
(x;\eta_{t})\right]  \right\vert \right\}  dx\leq C\frac{\overline{p}%
_{n}^{\frac{3}{a}}+\overline{p}_{n}^{\frac{1}{2}\left(  1+4/a\right)  }%
}{n^{1/2}}\leq C\left(  \frac{\overline{p}_{n}^{1+\frac{4}{a}}}{n}\right)
^{\frac{1}{2}}. \label{Lind3bnd}%
\end{equation}

\textbf{The second-order term (\ref{Lind2})}. Note that $\mathcal{I}_{t}%
^{(2)}(0;\eta)=\eta^{\prime}A_{t}\eta$ where $A_{t}$ depends upon
$D_{1},\ldots,D_{t-1}$ and $\eta_{t+1},\ldots,\eta_{n}$. In the standard
Lindeberg method, $\left\{  D_{t},t\in\left[  1,n\right]  \right\}  $ and
$\left\{  \eta_{t},t\in\left[  1,n\right]  \right\}  $ are both independent
variables with identical mean and variance, so that the second order term,
which writes as a sum of items $\mathbb{E}\left[  D_{t}^{\prime}A_{t}%
D_{t}\right]  -\mathbb{E}\left[  \eta_{t}^{\prime}A_{t}\eta_{t}\right]  $, is
equal to $0$ in this simpler case. However this does not hold in our case. In
this step, the second order term is dealt with by removing from $\mathcal{I}%
_{t}^{(2)}(0;\eta)$ a block $\sum_{j=1}^{p}K_{jp}\sum_{s=t-\ell}^{t-1}D_{js}$
and by changing the $D_{jt}$ into $D_{jt}^{t-\ell+1}=\mathbb{E}\left[
D_{jt}\left\vert e_{t},\ldots,e_{t-\ell+1}\right.  \right]  $.

Observe that $\mathcal{I}_{t}^{(2)}(0;\eta)=\mathcal{I}_{1t}^{(2)}%
(0;\eta)+\mathcal{I}_{2t}^{(2)}(0;\eta)+\mathcal{I}_{3t}^{(2)}(0;\eta
)+\mathcal{I}_{4t}^{(2)}(0;\eta)$ with, dropping the dependence upon $0$ and
$\eta$,%

\begin{align*}
\mathcal{I}_{1t}^{(2)}  &  =\left(  \frac{1}{e}-1\right)  I_{tn}%
^{(1)}\mathcal{M}_{t}^{1-2e}\left(  \sum_{p=1}^{\overline{p}_{n}}\Sigma
_{pt}^{e-1}\Sigma_{pt}^{(1)}\right)  ^{2},\quad I_{tn}^{(1)}=I^{\prime}\left(
\mathcal{M}_{t}\right)  ,\\
\mathcal{I}_{2t}^{(2)}  &  =I_{tn}^{(1)}\mathcal{M}_{t}^{1-e}\sum
_{p=1}^{\overline{p}_{n}}\Sigma_{pt}^{e-1}\Sigma_{pt}^{(2)},\quad
\mathcal{I}_{3t}^{(2)}=\left(  e-1\right)  I_{tn}^{(1)}\mathcal{M}_{t}%
^{1-e}\sum_{p=1}^{\overline{p}_{n}}\Sigma_{pt}^{e-1}\left(  \Sigma_{pt}%
^{(1)}\right)  ^{2},\\
\mathcal{I}_{4t}^{(2)}  &  =I^{^{\prime\prime}}\left(  \mathcal{M}_{t}\right)
\left(  \mathcal{M}_{t}^{1-e}\sum_{p=1}^{\overline{p}_{n}}\Sigma_{pt}%
^{e-1}\Sigma_{pt}^{(1)}\right)  ^{2}.
\end{align*}
Observe $\mathcal{M}_{t}\left(  0;D_{t}\right)  =\mathcal{M}_{t}\left(
0;\eta_{t}\right)  $ and $\Sigma_{pt}\left(  0;D_{t}\right)  =\Sigma
_{pt}\left(  0;\eta_{t}\right)  $ and that these quantities do not depend upon
$\eta_{t}$ or $D_{t}$. We shall first focus on $\mathcal{I}_{1t}^{(2)}$. Let
$\ell\geq2\overline{p}_{n}$ be an integer number. Define, for $y\in\left[
0,1\right]  $,%
\begin{align*}
\mathfrak{S}_{pt}\left(  y;\eta\right)   &  =\frac{2\sum_{j=1}^{p}%
K_{jp}\left(  \sum_{s=j+1}^{t-\ell-1}D_{js}+y\sum_{s=t-\ell}^{t-1}D_{js}%
+\sum_{s=t+1}^{n}\eta_{js}\right)  \eta_{j}}{n\sigma^{4}V_{\Delta}\left(
p\right)  },\\
\mathfrak{S}_{pt}\left(  y\right)   &  =\mathfrak{S}_{pt}\left(
y;yD_{t}+\left(  1-y\right)  D_{t}^{t-\ell+1}\right)  ,\\
\mathfrak{T}_{pt}\left(  y;\eta\right)   &  =\check{s}_{pt}^{(2)}%
(y;\eta)=\frac{2\sum_{j=1}^{p}K_{jp}\eta_{j}^{2}}{n\sigma^{4}V_{\Delta}\left(
p\right)  },\quad\mathfrak{T}_{pt}\left(  y\right)  =\mathfrak{T}_{pt}\left(
y;yD_{t}+\left(  1-y\right)  D_{t}^{t-\ell+1}\right)  ,
\end{align*}
which are such that $\mathfrak{S}_{pt}\left(  1;\eta\right)  =\check{s}%
_{pt}^{(1)}(0;\eta)$, $\mathfrak{S}_{pt}\left(  1\right)  =\check{s}%
_{pt}^{(1)}(0;D_{t})$, $\mathfrak{T}_{pt}\left(  1\right)  =\check{s}%
_{pt}^{(2)}(0;D_{t})$. Define also%
\begin{align*}
\mathbf{M}_{jt}\left(  y\right)   &  =\sum_{s=j+1}^{t-\ell-1}D_{js}%
+y\sum_{s=t-\ell}^{t-1}D_{js}+\sum_{s=t+1}^{n}\eta_{js},\quad\mathbf{R}%
_{jt}\left(  y\right)  =\frac{\mathbf{M}_{jt}\left(  y\right)  }{n},\\
\mathbf{s}_{pt}\left(  y\right)   &  =\frac{n\sum_{j=1}^{p}K_{jp}%
\mathbf{R}_{jt}^{2}\left(  y\right)  -\sigma^{4}E_{\Delta}(p)}{\sigma
^{4}V_{\Delta}\left(  p\right)  },\quad\mathbf{\Sigma}_{pt}\left(  y\right)
=f\left(  \mathbf{s}_{pt}\left(  y\right)  \right)  ,\\
\widetilde{\Sigma}_{pt}^{\left(  1\right)  }\left(  y;\eta\right)   &
=f^{\left(  1\right)  }\left(  \mathbf{s}_{pt}\left(  y\right)  \right)
\mathfrak{S}_{pt}\left(  y;\eta\right)  ,\\
\widetilde{\Sigma}_{pt}^{\left(  2\right)  }\left(  y;\eta\right)   &
=f^{\left(  1\right)  }\left(  \mathbf{s}_{pt}\left(  y\right)  \right)
\mathfrak{T}_{pt}\left(  y;\eta\right)  +f^{\left(  2\right)  }\left(
\mathbf{s}_{pt}\left(  y\right)  \right)  \left(  \mathfrak{S}_{pt}\left(
y;\eta\right)  \right)  ^{2},\\
\widetilde{\Sigma}_{pt}^{\left(  1\right)  }\left(  y\right)   &
=\widetilde{\Sigma}_{pt}^{\left(  1\right)  }\left(  y;yD_{t}+\left(
1-y\right)  D_{t}^{t-\ell+1}\right)  ,\\
\widetilde{\Sigma}_{pt}^{\left(  2\right)  }\left(  y;\eta\right)   &
=\widetilde{\Sigma}_{pt}^{\left(  2\right)  }\left(  y;yD_{t}+\left(
1-y\right)  D_{t}^{t-\ell+1}\right)  ,\\
\mathfrak{M}_{t}\left(  y\right)   &  =\left(  \sum_{p=1}^{\overline{p}_{n}%
}\mathbf{\Sigma}_{pt}^{e}\left(  y\right)  \right)  ^{\frac{1}{e}}%
,\quad\mathfrak{I}_{tn}^{(1)}\left(  y\right)  =I^{\prime}\left(
\mathfrak{M}_{t}\left(  y\right)  \right)  ,
\end{align*}
and the counterpart of $\mathcal{I}_{1t}^{(2)}\left(  0;\eta_{t}\right)  $ and
$\mathcal{I}_{1t}^{(2)}\left(  0;D_{t}\right)  $ as%
\begin{align*}
\mathfrak{I}_{t}\left(  y;\eta\right)   &  =\left(  \frac{1}{e}-1\right)
\mathfrak{I}_{tn}^{(1)}\left(  y\right)  \mathfrak{M}_{t}^{1-2e}\left(
y\right)  \left(  \sum_{p=1}^{\overline{p}_{n}}\mathbf{\Sigma}_{pt}%
^{e-1}\left(  y\right)  \widetilde{\Sigma}_{pt}^{\left(  1\right)  }\left(
y;\eta\right)  \right)  ^{2},\\
\mathfrak{I}_{t}\left(  y\right)   &  =\mathfrak{I}_{t}\left(  y;yD_{t}%
+\left(  1-y\right)  D_{t}^{t-\ell+1}\right)  .
\end{align*}
Observe that $\mathcal{I}_{1t}^{(2)}\left(  0;\eta_{t}\right)  =\mathfrak{I}%
_{t}\left(  1;\eta_{t}\right)  $ and $\mathcal{I}_{1t}^{(2)}\left(
0;D_{t}\right)  =\mathfrak{I}_{t}\left(  1\right)  $. Hence $\mathbb{E}\left[
\mathcal{I}_{1t}^{(2)}\left(  0;D_{t}\right)  -\mathcal{I}_{1t}^{(2)}\left(
0;\eta_{t}\right)  \right]  =\mathbb{E}\left[  \mathfrak{I}_{t}\left(
1\right)  -\mathfrak{I}_{t}\left(  1;\eta_{t}\right)  \right]  $ and%
\begin{align}
\mathbb{E}\left[  \mathcal{I}_{1t}^{(2)}\left(  0;D_{t}\right)  -\mathcal{I}%
_{1t}^{(2)}\left(  0;\eta_{t}\right)  \right]   &  =\mathbb{E}\left[
\mathfrak{I}_{t}\left(  0\right)  -\mathfrak{I}_{t}\left(  0;\eta_{t}\right)
\right] \label{I12}\\
&  +\int_{0}^{1}\mathbb{E}\left[  \mathfrak{I}_{t}^{(1)}\left(  y\right)
-\mathfrak{I}_{t}^{(1)}\left(  y;\eta_{t}\right)  \right]  dy, \label{IntI12}%
\end{align}
where $\mathfrak{I}_{t}^{(1)}\left(  y\right)  =d\mathfrak{I}_{t}\left(
y\right)  /dy$ and $\mathfrak{I}_{t}^{(1)}\left(  y;\eta_{t}\right)
=d\mathfrak{I}_{t}\left(  y;\eta_{t}\right)  /dy$.

We first consider the integral item $\int_{0}^{1}\left\vert \mathbb{E}\left[
\mathfrak{I}_{t}^{(1)}\left(  y\right)  \right]  \right\vert dy$ from
(\ref{IntI12}) and first compute $\mathfrak{I}_{1t}^{(1)}\left(  y\right)  $.
Define%
\begin{align*}
\mathfrak{S}_{pt}^{(1)}\left(  y\right)   &  =\frac{d\mathfrak{S}_{pt}\left(
y\right)  }{dy}=\frac{2\sum_{j=1}^{p}K_{jp}\left(  \sum_{s=t-\ell}^{t-1}%
D_{js}\right)  \left(  yD_{jt}+\left(  1-y\right)  D_{jt}^{t-\ell+1}\right)
}{n\sigma^{4}V_{\Delta}\left(  p\right)  }\\
&  +\frac{2\sum_{j=1}^{p}K_{jp}\left(  \sum_{s=j+1}^{t-\ell-1}D_{js}%
+y\sum_{s=t-\ell}^{t-1}D_{js}+\sum_{s=t+1}^{n}\eta_{js}\right)  \left(
D_{jt}^{t-\ell+1}-D_{jt}\right)  }{n\sigma^{4}V_{\Delta}\left(  p\right)  },
\end{align*}%
\begin{align*}
\mathfrak{T}_{pt}^{(1)}\left(  y\right)   &  =\frac{d\mathfrak{T}_{pt}\left(
y\right)  }{dy}=\frac{4\sum_{j=1}^{p}K_{jp}\left(  yD_{jt}+\left(  1-y\right)
D_{jt}^{t-\ell+1}\right)  \left(  D_{jt}-D_{jt}^{t-\ell+1}\right)  }%
{n\sigma^{4}V_{\Delta}\left(  p\right)  },\\
\mathbf{s}_{pt}^{(1)}\left(  y\right)   &  =\frac{d\mathbf{s}_{pt}\left(
y\right)  }{dy}=\frac{2\sum_{j=1}^{p}K_{jp}\mathbf{M}_{jt}\left(  y\right)
\sum_{s=t-\ell}^{t-1}D_{js}}{n\sigma^{4}V_{\Delta}\left(  p\right)  },\\
\widetilde{\Sigma}_{pt}^{\left(  1,1\right)  }\left(  y\right)   &
=\frac{d\widetilde{\Sigma}_{pt}^{\left(  1\right)  }\left(  y\right)  }%
{dy}=f^{\left(  2\right)  }\left(  \mathbf{s}_{pt}\left(  y\right)  \right)
\mathbf{s}_{pt}^{(1)}\left(  y\right)  \mathfrak{S}_{pt}\left(  y\right)
+f^{\left(  1\right)  }\left(  \mathbf{s}_{pt}\left(  y\right)  \right)
\mathfrak{S}_{pt}^{\left(  1\right)  }\left(  y\right)  ,
\end{align*}%
\begin{align*}
\widetilde{\Sigma}_{pt}^{\left(  2,1\right)  }\left(  y\right)   &
=\frac{d\widetilde{\Sigma}_{pt}^{\left(  2\right)  }\left(  y\right)  }%
{dy}=f^{\left(  2\right)  }\left(  \mathbf{s}_{pt}\left(  y\right)  \right)
\mathbf{s}_{pt}^{(1)}\left(  y\right)  \mathfrak{T}_{pt}\left(  y\right)
+f^{\left(  1\right)  }\left(  \mathbf{s}_{pt}\left(  y\right)  \right)
\mathfrak{T}_{pt}^{\left(  1\right)  }\left(  y\right) \\
&  +f^{\left(  3\right)  }\left(  \mathbf{s}_{pt}\left(  y\right)  \right)
\mathbf{s}_{pt}^{(1)}\left(  y\right)  \left(  \mathfrak{S}_{pt}\left(
y\right)  \right)  ^{2}+2f^{\left(  2\right)  }\left(  \mathbf{s}_{pt}\left(
y\right)  \right)  \mathfrak{S}_{pt}\left(  y\right)  \mathfrak{S}%
_{pt}^{\left(  1\right)  }\left(  y\right)  ,
\end{align*}%
\[
\mathfrak{I}_{tn}^{(2)}\left(  y\right)  =I^{\prime\prime}\left(
\mathfrak{M}_{t}\left(  y\right)  \right)  ,
\]
and%
\begin{align*}
\mathfrak{I}_{1t}^{(1)}\left(  y\right)   &  =\left(  \frac{1}{e}-1\right)
\mathfrak{I}_{tn}^{(2)}\left(  y\right)  \mathfrak{M}_{t}^{2-3e}\left(
y\right)  \left(  \sum_{p=1}^{\overline{p}_{n}}\mathbf{\Sigma}_{pt}%
^{e-1}\left(  y\right)  \widetilde{\Sigma}_{pt}^{\left(  1\right)  }\left(
y\right)  \right)  ^{2}\sum_{p=1}^{\overline{p}_{n}}\mathbf{\Sigma}_{pt}%
^{e-1}\left(  y\right)  \mathbf{\Sigma}_{pt}^{(1)}\left(  y\right)  ,\\
\mathfrak{I}_{2t}^{(1)}\left(  y\right)   &  =\left(  \frac{1}{e}-1\right)
\left(  \frac{1}{e}-2\right)  \mathfrak{I}_{tn}^{(1)}\left(  y\right)
\mathfrak{M}_{t}^{1-3e}\left(  y\right)  \left(  \sum_{p=1}^{\overline{p}_{n}%
}\mathbf{\Sigma}_{pt}^{e-1}\left(  y\right)  \widetilde{\Sigma}_{pt}^{\left(
1\right)  }\left(  y\right)  \right)  ^{2}\sum_{p=1}^{\overline{p}_{n}%
}\mathbf{\Sigma}_{pt}^{e-1}\left(  y\right)  \mathbf{\Sigma}_{pt}^{(1)}\left(
y\right)  ,\\
\mathfrak{I}_{3t}^{(1)}\left(  y\right)   &  =2\left(  \frac{1}{e}-1\right)
\left(  e-1\right)  \mathfrak{I}_{tn}^{(1)}\left(  y\right)  \mathfrak{M}%
_{t}^{1-2e}\left(  y\right)  \left(  \sum_{p=1}^{\overline{p}_{n}%
}\mathbf{\Sigma}_{pt}^{e-1}\left(  y\right)  \widetilde{\Sigma}_{pt}^{\left(
1\right)  }\left(  y\right)  \right)  \left(  \sum_{p=1}^{\overline{p}_{n}%
}\mathbf{\Sigma}_{pt}^{e-2}\left(  y\right)  \left(  \mathbf{\Sigma}%
_{pt}^{(1)}\left(  y\right)  \right)  ^{2}\right)  ,\\
\mathfrak{I}_{4t}^{(1)}\left(  y\right)   &  =2\left(  \frac{1}{e}-1\right)
\mathfrak{I}_{tn}^{(1)}\left(  y\right)  \mathfrak{M}_{t}^{1-2e}\left(
y\right)  \left(  \sum_{p=1}^{\overline{p}_{n}}\mathbf{\Sigma}_{pt}%
^{e-1}\left(  y\right)  \widetilde{\Sigma}_{pt}^{\left(  1\right)  }\left(
y\right)  \right)  \left(  \sum_{p=1}^{\overline{p}_{n}}\mathbf{\Sigma}%
_{pt}^{e-1}\left(  y\right)  \widetilde{\Sigma}_{pt}^{\left(  1,1\right)
}\left(  y\right)  \right)  .
\end{align*}
To bound the moments of $\widetilde{\Sigma}_{pt}^{\left(  1\right)  }\left(
y\right)  $, $\widetilde{\Sigma}_{pt}^{\left(  1,1\right)  }\left(  y\right)
$ and $\mathbf{\Sigma}_{pt}^{(1)}\left(  y\right)  $, consider first
$\left\Vert \mathfrak{S}_{pt}\left(  y\right)  \right\Vert _{3a}$, $\left\Vert
\mathfrak{S}_{pt}^{(1)}\left(  y\right)  \right\Vert _{3a}$ and $\left\Vert
\mathbf{s}_{pt}^{(1)}\left(  y\right)  \right\Vert _{3a}$. For $\left\Vert
\mathfrak{S}_{pt}\left(  y\right)  \right\Vert _{3a}$ and $\left\Vert
\mathfrak{S}_{pt}^{(1)}\left(  y\right)  \right\Vert _{3a}$, (\ref{Sig1}), the
Burkholder inequality, (\ref{Shaolem52}) $\overline{p}_{n}=O\left(
n^{1/2}\right)  $, $2\overline{p}_{n}\leq\ell\leq3\overline{p}_{n}$ and
$\Theta_{6a}\left(  \ell-\overline{p}_{n}\right)  \leq C\overline{p}_{n}^{-1}$
give%
\begin{align*}
&  \left\Vert \mathfrak{S}_{pt}\left(  y\right)  \right\Vert _{3a}\\
&  \leq\left\Vert \frac{2\sum_{j=1}^{p}K_{jp}\left(  \sum_{s=j+1}^{t-\ell
-1}D_{js}+y\sum_{s=t-\ell}^{t-1}D_{js}+\sum_{s=t+1}^{n}\eta_{js}\right)
D_{jt}}{n\sigma^{4}V_{\Delta}\left(  p\right)  }\right\Vert _{3a}\\
&  +2\left\vert 1-y\right\vert \sum_{j=1}^{p}\frac{\left\vert K_{jp}%
\right\vert }{n\sigma^{4}V_{\Delta}\left(  p\right)  }\left\Vert \left(
\sum_{s=j+1}^{t-\ell-1}D_{js}+y\sum_{s=t-\ell}^{t-1}D_{js}+\sum_{s=t+1}%
^{n}\eta_{js}\right)  \right\Vert _{6a}\left\Vert D_{jt}-D_{jt}^{t-\ell
+1}\right\Vert _{6a}\\
&  \leq C\left(  \frac{1}{n^{1/2}}+\frac{\overline{p}_{n}}{n}+\left(
\frac{\overline{p}_{n}}{n}\right)  ^{1/2}\Theta_{6a}\left(  \ell-\overline
{p}_{n}\right)  \right)  \leq\frac{C}{n^{1/2}},
\end{align*}

\begin{align*}
&  \left\Vert \mathfrak{S}_{pt}^{(1)}\left(  y\right)  \right\Vert _{3a}\\
&  \text{ }\leq\left\Vert \frac{2\sum_{j=1}^{p}K_{jp}\left(  \sum_{s=t-\ell
}^{t-1}D_{js}\right)  D_{jt}}{n\sigma^{4}V_{\Delta}\left(  p\right)
}\right\Vert _{3a}\\
&  \text{ }+2\left\vert 1-y\right\vert \sum_{j=1}^{p}\frac{\left\vert
K_{jp}\right\vert }{n\sigma^{4}V_{\Delta}\left(  p\right)  }\left\Vert
\sum_{s=t-\ell}^{t-1}D_{js}\right\Vert _{6a}\left\Vert D_{jt}-D_{jt}%
^{t-\ell+1}\right\Vert _{6a}\\
&  \text{ }+2\sum_{j=1}^{p}\frac{\left\vert K_{jp}\right\vert }{n\sigma
^{4}V_{\Delta}\left(  p\right)  }\left\Vert \sum_{s=j+1}^{t-\ell-1}%
D_{js}+y\sum_{s=t-\ell}^{t-1}D_{js}+\sum_{s=t+1}^{n}\eta_{js}\right\Vert
_{6a}\left\Vert D_{jt}-D_{jt}^{t-\ell+1}\right\Vert _{6a}\\
&  \text{ }\leq C\left(  \frac{\ell^{1/2}}{n}+\frac{\ell^{1/2}\overline{p}%
_{n}^{1/2}}{n}\Theta_{6a}\left(  \ell-\overline{p}_{n}\right)  +\left(
\frac{\overline{p}_{n}}{n}\right)  ^{1/2}\Theta_{6a}\left(  \ell-\overline
{p}_{n}\right)  \right) \\
&  \text{ }\leq C\left(  \frac{\overline{p}_{n}^{1/2}}{n}+\frac{1}{\left(
n\overline{p}_{n}\right)  ^{1/2}}\right)  ,
\end{align*}%
\[
\left\Vert \mathfrak{T}_{pt}\left(  y\right)  \right\Vert _{3a}\leq
C\frac{\overline{p}_{n}^{1/2}}{n},\quad\left\Vert \mathfrak{T}_{pt}^{\left(
1\right)  }\left(  y\right)  \right\Vert _{3a}\leq\frac{C}{n\overline{p}_{n}%
}.
\]
For $\left\Vert \mathbf{s}_{pt}^{(1)}\left(  y\right)  \right\Vert _{3a}$
(\ref{Sig1}), $\overline{p}_{n}=O\left(  n^{1/2}\right)  $ and the Burkholder
inequality give%
\begin{align*}
&  \left\Vert \mathbf{s}_{pt}^{(1)}\left(  y\right)  \right\Vert _{3a}\\
&  \text{ }\leq\left\Vert 2\sum_{s_{1}=t-\ell}^{t-1}\sum_{j=1}^{p}\frac
{K_{jp}}{n\sigma^{4}V_{\Delta}\left(  p\right)  }\left(  \sum_{s_{2}%
=j+1}^{t-\ell-1}D_{js_{2}}\right)  D_{js_{1}}\right\Vert _{3a}+\left\Vert
\frac{2\sum_{j=1}^{p}K_{jp}\left(  \sum_{s=t-\ell}^{t-1}D_{js}\right)  ^{2}%
}{n\sigma^{4}V_{\Delta}\left(  p\right)  }\right\Vert _{3a}\\
&  \text{ }+\left\Vert \frac{2\sum_{j=1}^{p}K_{jp}\left(  \sum_{s=t-\ell
}^{t-1}D_{js}\right)  \left(  \sum_{s=t+1}^{n}\eta_{js}\right)  }{n\sigma
^{4}V_{\Delta}\left(  p\right)  }\right\Vert _{3a}\\
&  \text{ }\leq C\left(  \sum_{s_{1}=t-\ell}^{t-1}\left\Vert \sum_{j=1}%
^{p}\frac{K_{jp}}{n\sigma^{4}V_{\Delta}\left(  p\right)  }\left(  \sum
_{s_{2}=j+1}^{t-\ell-1}D_{js_{2}}\right)  D_{js_{1}}\right\Vert _{3a}%
^{2}\right)  ^{1/2}+C\sum_{j=1}^{p}\frac{\left\vert K_{jp}\right\vert
}{n\sigma^{4}V_{\Delta}\left(  p\right)  }\left\Vert \sum_{s=t-\ell}%
^{t-1}D_{js}\right\Vert _{6a}^{2}\\
&  \text{ }+C\left\Vert \frac{\left(  \sum_{j=1}^{p}K_{jp}^{2}\left(
\sum_{s=t-\ell}^{t-1}D_{js}\right)  ^{2}\right)  ^{1/2}}{\left(  np\right)
^{1/2}}\right\Vert _{3a}\\
&  \text{ }\leq C\left(  \ell^{1/2}\left(  \frac{1}{n^{1/2}}+\frac
{\overline{p}_{n}}{n}\right)  +\frac{\overline{p}_{n}^{1/2}\ell}{n}+\frac
{\ell^{1/2}}{n^{1/2}}\right)  \leq C\left(  \frac{\overline{p}_{n}}{n}\right)
^{1/2}.
\end{align*}
These bounds and (\ref{f}) give, uniformly in $y$, $p$ and $t$,%
\begin{align*}
\left\Vert \widetilde{\Sigma}_{pt}^{\left(  1\right)  }\left(  y\right)
\right\Vert _{3a}  &  \leq\frac{C}{n^{1/2}},\quad\left\Vert \mathbf{\Sigma
}_{pt}^{\left(  1\right)  }\left(  y\right)  \right\Vert _{3a}\leq C\left(
\frac{\overline{p}_{n}}{n}\right)  ^{1/2},\\
\left\Vert \widetilde{\Sigma}_{pt}^{\left(  1,1\right)  }\left(  y\right)
\right\Vert _{3a/2}  &  \leq C\left(  \frac{\overline{p}_{n}^{1/2}}{n}+\left(
\frac{\overline{p}_{n}}{n}\right)  ^{3/2}+\frac{\overline{p}_{n}^{1/2}%
}{n^{3/2}}+\frac{1}{n\overline{p}_{n}^{1/2}}\right)  \leq C\frac{\overline
{p}_{n}^{1/2}}{n}.
\end{align*}
Now, arguing as for the study of (\ref{Lind3}), $e=O\left(  \overline{p}%
_{n}^{1/a}\right)  $ give uniformly in $p$, $t$ and $y$,%
\[
\mathbb{E}\left[  \left\vert \mathfrak{I}_{1t}^{(1)}\left(  y\right)
\right\vert \right]  +\mathbb{E}\left[  \left\vert \mathfrak{I}_{2t}%
^{(1)}\left(  y\right)  \right\vert \right]  +\mathbb{E}\left[  \left\vert
\mathfrak{I}_{4t}^{(1)}\left(  y\right)  \right\vert \right]  \leq
C\frac{\overline{p}_{n}^{1/2+3/a}}{n^{3/2}},\quad\mathbb{E}\left[  \left\vert
\mathfrak{I}_{3t}^{(1)}\left(  y\right)  \right\vert \right]  \leq
C\frac{\overline{p}_{n}^{1+3/a}}{n^{3/2}}.
\]
It then follows $\sum_{t=1}^{n}\int_{0}^{1}\left\vert \mathbb{E}\left[
\mathfrak{I}_{t}^{(1)}\left(  y\right)  \right]  \right\vert dy\leq
C\overline{p}_{n}^{1+3/a}/n^{1/2}$. Since $\sum_{t=1}^{n}\int_{0}%
^{1}\left\vert \mathbb{E}\left[  \mathfrak{I}_{t}^{(1)}\left(  y;\eta
_{t}\right)  \right]  \right\vert dy$ satisfies a similar bound, we have for
(\ref{IntI12}),%
\[
\sum_{t=1}^{n}\left\vert \int_{0}^{1}\mathbb{E}\left[  \mathfrak{I}_{t}%
^{(1)}\left(  y\right)  -\mathfrak{I}_{t}^{(1)}\left(  y;\eta_{t}\right)
\right]  dy\right\vert \leq C\frac{\overline{p}_{n}^{1+3/a}}{n^{1/2}}.
\]
Consider now (\ref{I12}). Since $D_{jt}^{t-\ell+1}$ and $\eta_{t}$ are
independent of the $\mathfrak{J}_{tn}^{(1)}\left(  0\right)  $, $\mathfrak{M}%
_{t}^{1-2e}\left(  0\right)  $ and $\mathbf{\Sigma}_{pt}\left(  0\right)  $,
we have using (\ref{Vareta}),%
\begin{align*}
&  \mathbb{E}\left[  \mathfrak{I}_{t}\left(  0\right)  -\mathfrak{I}%
_{t}\left(  0;\eta_{t}\right)  \right] \\
&  =\frac{4}{n}\mathbb{E}\left[  \left(  \frac{1}{e}-1\right)  \mathfrak{J}%
_{tn}^{(1)}\left(  0\right)  \mathfrak{M}_{t}^{1-2e}\left(  0\right)  \right.
\\
&  \sum_{p_{1},p_{2}=1}^{\overline{p}}\mathbf{\Sigma}_{p_{1}t}^{e-1}\left(
0\right)  \mathbf{\Sigma}_{p_{2}t}^{e-1}\left(  0\right)  f\left(
\mathbf{\Sigma}_{p_{1}t}^{e-1}\left(  0\right)  \right)  f\left(
\mathbf{\Sigma}_{p_{2}t}^{e-1}\left(  0\right)  \right)  \sum_{j_{1}=1}%
^{p_{1}}\sum_{j_{2}=1}^{p_{2}}\left(  \mathbb{E}\left[  D_{j_{1}t}^{t-\ell
+1}D_{j_{2}t}^{t-\ell+1}\right]  -\mathbb{E}\left[  \eta_{j_{1}t}\eta_{j_{2}%
t}\right]  \right) \\
&  \left.  \frac{K_{j_{1}p_{1}}\left(  \sum_{s_{1}=j_{1}+1}^{t-\ell+1}%
D_{j_{1}s_{1}}+\sum_{s_{1}=t-\ell}^{n}\eta_{j_{1}s_{1}}\right)  }%
{n^{1/2}\sigma^{4}V_{\Delta}\left(  p_{1}\right)  }\frac{K_{j_{2}p_{2}}\left(
\sum_{s_{2}=j_{2}+1}^{t-\ell+1}D_{j_{2}s_{2}}+\sum_{s_{2}=t-\ell}^{n}%
\eta_{j_{2}s_{2}}\right)  }{n^{1/2}\sigma^{4}V_{\Delta}\left(  p_{2}\right)
}\right] \\
&  =0.
\end{align*}
Hence (\ref{I12}) and (\ref{IntI12}) give%
\[
\left\vert \sum_{t=1}^{n}\mathbb{E}\left[  \mathcal{I}_{1t}^{(2)}\left(
0;D_{t}\right)  -\mathcal{I}_{1t}^{(2)}\left(  0;\eta_{t}\right)  \right]
\right\vert \leq C\frac{\overline{p}_{n}^{1+3/a}}{n^{1/2}}.
\]

To study $\left\vert \mathbb{E}\left[  \mathcal{I}_{2t}^{(2)}\left(
0;D_{t}\right)  -\mathcal{I}_{2t}^{(2)}\left(  0;\eta_{t}\right)  \right]
\right\vert $, observe that, uniformly with respect to $p$, $t$ and $y$,%
\begin{align*}
\max\left(  \left\Vert \widetilde{\Sigma}_{pt}^{\left(  2\right)  }\left(
y\right)  \right\Vert _{3a/2},\left\Vert \widetilde{\Sigma}_{pt}^{\left(
2\right)  }\left(  y;\eta_{t}\right)  \right\Vert _{3a/2}\right)   &  \leq
C\frac{\overline{p}_{n}^{1/2}}{n},\\
\max\left(  \left\Vert \widetilde{\Sigma}_{pt}^{\left(  2,1\right)  }\left(
y\right)  \right\Vert _{a},\left\Vert \widetilde{\Sigma}_{pt}^{\left(
2,1\right)  }\left(  y;\eta_{t}\right)  \right\Vert _{a}\right)   &  \leq
C\left(  \frac{\overline{p}_{n}}{n^{3/2}}+\frac{1}{n\overline{p}_{n}}\right)
.
\end{align*}
Arguing as for $\sum_{t=1}^{n}\mathbb{E}\left[  \mathcal{I}_{1t}^{(2)}\left(
0;D_{t}\right)  -\mathcal{I}_{1t}^{(2)}\left(  0;\eta_{t}\right)  \right]  $
gives $\left\vert \sum_{t=1}^{n}\mathbb{E}\left[  \mathcal{I}_{2t}%
^{(2)}\left(  0;D_{t}\right)  -\mathcal{I}_{2t}^{(2)}\left(  0;\eta
_{t}\right)  \right]  \right\vert \leq C\left(  \frac{\overline{p}_{n}%
^{1+2/a}}{n^{1/2}}+\frac{\overline{p}_{n}^{1/a}}{\overline{p}_{n}}\right)  $,
and provided $e=O\left(  \overline{p}_{n}^{1/(2a)}\right)  $%
\[
\left\vert \sum_{t=1}^{n}\mathbb{E}\left[  \mathcal{I}_{3t}^{(2)}\left(
0;D_{t}\right)  -\mathcal{I}_{3t}^{(2)}\left(  0;\eta_{t}\right)  \right]
\right\vert +\left\vert \sum_{t=1}^{n}\mathbb{E}\left[  \mathcal{I}_{4t}%
^{(2)}\left(  0;D_{t}\right)  -\mathcal{I}_{4t}^{(2)}\left(  0;\eta
_{t}\right)  \right]  \right\vert \leq C\frac{\overline{p}_{n}^{1+3/a}%
}{n^{1/2}}.
\]
It then follows%
\begin{equation}
\left\vert \sum_{t=1}^{n}\mathbb{E}\left[  \mathcal{I}_{t}^{(2)}\left(
0;D_{t}\right)  -\mathcal{I}_{t}^{(2)}\left(  0;\eta_{t}\right)  \right]
\right\vert \leq C\left(  \frac{\overline{p}_{n}^{1+3/a}}{n^{1/2}}+\frac
{1}{\overline{p}_{n}^{1-1/a}}\right)  . \label{Lind2bnd}%
\end{equation}
Substituting (\ref{Lind3bnd}), (\ref{Lind2bnd}) in (\ref{Lind3}),
(\ref{Lind2}) shows that the Lemma is proved.\hfill$\square$

\subsubsection{End of the proof of Proposition \ref{Selection}}

The rest of the proof is divided in 3 steps.\smallskip

\textbf{Step 1: Martingale approximation}. Let $\widetilde{S}_{p}$ and
$\check{S}_{p}$ be as in (\ref{TildeR}) and (\ref{Pseudomax}) respectively.
Let $\mathfrak{a}=4a/3$. The Cauchy-Schwarz inequality gives%
\begin{align*}
\left\vert \check{S}_{p}-\widetilde{S}_{p}\right\vert  &  =\sum_{j=1}%
^{p}\left(  K_{jp}\frac{1}{n^{1/2}}\left\vert M_{jn}-\sum_{t=j+1}^{n}%
u_{t}u_{t-j}\right\vert \times\frac{1}{n^{1/2}}\left\vert M_{jn}+\sum
_{t=j+1}^{n}u_{t}u_{t-j}\right\vert \right) \\
&  \leq C\left(  \sum_{j=1}^{p}\frac{1}{n}\left(  M_{jn}-\sum_{t=j+1}^{n}%
u_{t}u_{t-j}\right)  ^{2}\right)  ^{1/2}\left(  \sum_{j=1}^{p}\frac{1}%
{n}\left(  M_{jn}+\sum_{t=j+1}^{n}u_{t}u_{t-j}\right)  ^{2}\right)  ^{1/2}.
\end{align*}
Hence%
\begin{align*}
&  \left\Vert \check{S}_{p}-\widetilde{S}_{p}\right\Vert _{\mathfrak{a}/2}\\
&  \text{ }\leq C\mathbb{E}^{\frac{1}{\mathfrak{a}}}\left[  \left(  \sum
_{j=1}^{p}\frac{1}{n}\left(  M_{jn}-\sum_{t=j+1}^{n}u_{t}u_{t-j}\right)
^{2}\right)  ^{\frac{\mathfrak{a}}{2}}\right]  \mathbb{E}^{\frac
{1}{\mathfrak{a}}}\left[  \left(  \sum_{j=1}^{p}\frac{1}{n}\left(  M_{jn}%
+\sum_{t=j+1}^{n}u_{t}u_{t-j}\right)  ^{2}\right)  ^{\frac{\mathfrak{a}}{2}%
}\right]  .
\end{align*}
Observe now that (\ref{Cov2M}) gives%
\begin{align*}
&  \mathbb{E}^{\frac{1}{\mathfrak{a}}}\left[  \left(  \sum_{j=1}^{p}\frac
{1}{n}\left(  M_{jn}-\sum_{t=j+1}^{n}u_{t}u_{t-j}\right)  ^{2}\right)
^{\frac{\mathfrak{a}}{2}}\right] \\
&  \text{ }\leq\left(  \frac{1}{n}\sum_{j=1}^{p}\mathbb{E}^{\frac
{2}{\mathfrak{a}}}\left[  \left\vert M_{jn}-\sum_{t=j+1}^{n}u_{t}%
u_{t-j}\right\vert ^{\mathfrak{a}}\right]  \right)  ^{1/2}\leq C\left(
\frac{p}{n}\right)  ^{1/2}.
\end{align*}
Since the Burkholder inequality and $\max_{j}\mathbb{E}\left[  \left\vert
D_{jt}\right\vert ^{\mathfrak{a}}\right]  <\infty$ give $\max_{j\in\left[
1,\overline{p}_{n}\right]  }\mathbb{E}^{1/\mathfrak{a}}\left[  \left\vert
M_{jn}\right\vert ^{\mathfrak{a}}\right]  \leq Cn^{1/2}$, we also have%
\begin{align*}
&  \mathbb{E}^{\frac{1}{\mathfrak{a}}}\left[  \left(  \sum_{j=1}^{p}\frac
{1}{n}\left(  M_{jn}+\sum_{t=j+1}^{n}u_{t}u_{t-j}\right)  ^{2}\right)
^{\frac{\mathfrak{a}}{2}}\right] \\
&  \leq\left(  \frac{1}{n}\sum_{j=1}^{p}\left(  \mathbb{E}^{\frac
{1}{\mathfrak{a}}}\left[  \left\vert M_{jn}+\sum_{t=j+1}^{n}u_{t}%
u_{t-j}\right\vert ^{\mathfrak{a}}\right]  \right)  ^{2}\right)  ^{1/2}\\
&  \leq\left(  \frac{1}{n}\sum_{j=1}^{p}\left(  2\mathbb{E}^{\frac
{1}{\mathfrak{a}}}\left[  \left\vert M_{jn}\right\vert ^{\mathfrak{a}}\right]
+\mathbb{E}^{\frac{1}{\mathfrak{a}}}\left[  \left\vert \sum_{t=j+1}^{n}%
u_{t}u_{t-j}-M_{jn}\right\vert ^{\mathfrak{a}}\right]  \right)  ^{2}\right)
^{1/2}\\
&  \text{ }\leq\left(  \frac{p(Cn^{1/2}+C)^{2}}{n}\right)  ^{1/2}\leq
Cp^{1/2}.
\end{align*}
It then follows that $\left\Vert \check{S}_{p}-\widetilde{S}_{p}\right\Vert
_{\mathfrak{a}/2}\leq Cp/n^{1/2}$ and them $\max_{p\in\left[  1,\overline
{p}_{n}\right]  }\mathbb{E}\left[  \left\vert \left(  \check{S}_{p}%
-\widetilde{S}_{p}\right)  /p^{1/2}\right\vert ^{\mathfrak{a}/2}\right]  \leq
C\left(  \overline{p}_{n}/n\right)  ^{\mathfrak{a}/4}$. Hence the Markov
inequality gives%
\begin{align*}
&  \mathbb{P}\left(  \max_{p\in\left[  1,\overline{p}_{n}\right]  }\left\vert
\frac{\check{S}_{p}-\widetilde{S}_{p}}{p^{1/2}}\right\vert \geq t\right)
\leq\sum_{p=1}^{\overline{p}_{n}}\mathbb{P}\left(  \left\vert \frac{\check
{S}_{p}-\widetilde{S}_{p}}{p^{1/2}}\right\vert \geq t\right) \\
&  \text{ }\leq\frac{\overline{p}_{n}}{t^{\mathfrak{a}/2}}\max_{p\in\left[
1,\overline{p}_{n}\right]  }\mathbb{E}\left[  \left\vert \frac{\check{S}%
_{p}-\widetilde{S}_{p}}{p^{1/2}}\right\vert ^{\frac{\mathfrak{a}}{2}}\right]
\leq\frac{C}{t^{a/2}}\left(  \frac{\overline{p}_{n}^{1+\frac{4}{\mathfrak{a}}%
}}{n}\right)  ^{\mathfrak{a}/4},
\end{align*}
and $\overline{p}_{n}=o\left(  n^{1/\left(  2\left(  1+4/\mathfrak{a}\right)
\right)  }\right)  $ gives%
\begin{equation}
\max_{p\in\left[  1,\overline{p}_{n}\right]  }\left\vert \frac{\check{S}%
_{p}-\widetilde{S}_{p}}{p^{1/2}}\right\vert =o_{\mathbb{P}}(1).
\label{Check2tilde}%
\end{equation}

\textbf{Step 2: some Gaussian approximations}. Let $\gamma_{n}^{\prime}%
=\gamma_{n}\left(  1+\epsilon/2\right)  /\left(  1+\epsilon\right)  $.
(\ref{Gam}) gives $\gamma_{n}\geq\gamma_{n}^{\prime}\geq\widetilde{\gamma}%
_{n}=\left(  2\ln\ln\overline{p}_{n}\right)  ^{1/2}\left(  1+\epsilon
/3\right)  $. Consider a three times continuously differentiable function
$\iota\left(  x\right)  $ with $\max_{j=1,2,3}\sup_{x}\left\vert
\iota^{\left(  3\right)  }\left(  x\right)  \right\vert <\infty$ and
$\mathbb{I}\left(  x\geq0\right)  \leq\iota\left(  x\right)  \leq
\mathbb{I}\left(  x\geq-\epsilon\right)  $. Let $\mathcal{I}\left(  x\right)
=\iota\left(  x-\gamma_{n}^{\prime}\right)  $. Let $\check{s}_{p}$ be as in
(\ref{Pseudomax}). Then Lemma \ref{Lindeberg} with $e=\overline{p}%
_{n}^{1/(2a)}$, (\ref{f}) and (\ref{Pseudomax}), and Assumption \ref{Reg} give%
\begin{align*}
&  \mathbb{P}\left(  \max_{p\in\left[  2,\overline{p}_{n}\right]  }\left\{
\check{s}_{p}\right\}  \geq\gamma_{n}^{\prime}\right)  \leq\mathbb{P}\left(
\mathcal{M}\geq\gamma_{n}^{\prime}\right)  \leq\mathbb{E}\left[
\mathcal{I}\left(  \mathcal{M}\right)  \right] \\
&  \text{ }\leq\mathbb{E}\left[  \mathcal{I}\left(  \mathcal{M}_{1}\left(
\eta_{1}\right)  \right)  \right]  +o\left(  1\right)  \leq\mathbb{P}\left(
\mathcal{M}_{1}\left(  \eta_{1}\right)  \geq\gamma_{n}^{\prime}-\epsilon
\right)  +o\left(  1\right)  .
\end{align*}
We now look for a more explicit expression for the RHS. Recall that
$\mathcal{M}_{1}\left(  \eta_{1}\right)  =\left(  \sum_{p=1}^{\overline{p}%
_{n}}f^{e}\left(  \check{s}_{p1}\left(  1;\eta_{1}\right)  \right)  \right)
^{1/e}$. Consider $\Omega\left(  p\right)  =\left[  \omega_{1},\ldots
,\omega_{p}\right]  ^{\prime}$ where the $\omega_{p}$'s are i.i.d. standard
normal variables,
\begin{align*}
\mathcal{K}\left(  p\right)   &  =\operatorname*{Diag}\left(  \left(
1-j/n\right)  K_{jp},j=1,\ldots,p\right)  ,\\
\mathcal{C}_{\eta}\left(  p\right)   &  =\left[  \operatorname*{Cov}\left(
\eta_{j_{1}t},\eta_{j_{2}t}\right)  ,j_{1},j_{2}=1,\ldots,p\right]  ,\\
\mathcal{V}_{\eta}\left(  p\right)   &  =\mathcal{C}_{\eta}^{1/2}\left(
p\right)  \mathcal{K}\left(  p\right)  \mathcal{C}_{\eta}^{1/2}\left(
p\right)  ,
\end{align*}
and $\mathcal{D}_{\eta}\left(  p\right)  =\operatorname*{Diag}\left(  \left(
1-j/n\right)  K_{jp}\operatorname*{Var}\left(  \eta_{jt}\right)
,j=1,\ldots,p\right)  $ the $p\times p$ diagonal matrix obtained from the
diagonal entries of $\mathcal{V}_{\eta}\left(  p\right)  $. Then the
$\check{s}_{p1}\left(  1;\eta_{1}\right)  $, $p=1,\ldots,\overline{p}_{n}$,
have the same joint distribution than%
\[
\tilde{s}_{p}=\frac{\Omega\left(  p\right)  ^{\prime}\mathcal{V}_{\eta}\left(
p\right)  \Omega\left(  p\right)  -\sigma^{4}E_{\Delta}\left(  p\right)
}{\sigma^{4}V_{\Delta}\left(  p\right)  },\quad p=1,\ldots,\overline{p}_{n},
\]
so that $\mathcal{M}_{1}\left(  \eta_{1}\right)  $ and $\widetilde
{\mathcal{M}}=\left(  \sum_{p=1}^{\overline{p}_{n}}f^{e}\left(  \tilde{s}%
_{p}\right)  \right)  ^{1/e}$ have the same distribution, and then%
\[
\mathbb{P}\left(  \max_{p\in\left[  2,\overline{p}_{n}\right]  }\left\{
\check{s}_{p}\right\}  \geq\gamma_{n}^{\prime}\right)  \leq\mathbb{P}\left(
\widetilde{\mathcal{M}}\geq\gamma_{n}^{\prime}-\epsilon\right)  +o\left(
1\right)  .
\]
Define now%
\[
\bar{s}_{p}=\frac{\Omega\left(  p\right)  ^{\prime}\mathcal{D}_{\eta}\left(
p\right)  \Omega\left(  p\right)  -\sigma^{4}E_{\Delta}\left(  p\right)
}{\sigma^{4}V_{\Delta}\left(  p\right)  }=\frac{\sum_{j=1}^{p}\left(
1-\frac{j}{n}\right)  K_{jp}\operatorname*{Var}\left(  \eta_{jt}\right)
\omega_{j}^{2}-\sigma^{4}E_{\Delta}\left(  p\right)  }{\sigma^{4}V_{\Delta
}\left(  p\right)  }.
\]
Then for all $p=1,\ldots,\overline{p}_{n}$,
\begin{align*}
&  \left\vert \tilde{s}_{p}-\bar{s}_{p}\right\vert =\left\vert \frac
{\Omega\left(  p\right)  ^{\prime}\left(  \mathcal{V}_{\eta}\left(  p\right)
-\mathcal{D}_{\eta}\left(  p\right)  \right)  \Omega\left(  p\right)  }%
{\sigma^{4}V_{\Delta}\left(  p\right)  }\right\vert \\
&  \text{ }\leq C\sum_{1\leq j_{1}\neq j_{2}\leq p}\left\vert
\operatorname*{Cov}\left(  \left(  1-\frac{j_{1}}{n}\right)  ^{1/2}K_{j_{1}%
p}^{1/2}\eta_{j_{1}t},\left(  1-\frac{j_{2}}{n}\right)  ^{1/2}K_{j_{2}p}%
^{1/2}\eta_{j_{2}t}\right)  \right\vert \left\vert \omega_{j_{1}}\right\vert
\left\vert \omega_{j_{2}}\right\vert \\
&  \text{ }\leq C\sum_{1\leq j_{1}\neq j_{2}\leq\overline{p}_{n}}\left\vert
\operatorname*{Cov}\left(  \eta_{j_{1}t},\eta_{j_{2}t}\right)  \right\vert
\left\vert \omega_{j_{1}}\right\vert \left\vert \omega_{j_{2}}\right\vert
=O_{\mathbb{P}}\left(  1\right)  ,
\end{align*}
by Lemma \ref{Sumvareta}. Hence since $f\left(  x\right)  \leq2\vee x$ by
(\ref{f}) and using (\ref{Smthmax}),%
\begin{align*}
\widetilde{\mathcal{M}}  &  \leq\left(  1+O\left(  \frac{\ln n}{\overline
{p}_{n}^{1/(2a)}}\right)  \right)  \max_{p\in\left[  2,\overline{p}%
_{n}\right]  }\left\{  2\vee\tilde{s}_{p}\right\}  \leq\left(  1+O\left(
\frac{\ln n}{\overline{p}_{n}^{1/(2a)}}\right)  \right)  2\vee\max
_{p\in\left[  2,\overline{p}_{n}\right]  }\left\{  \tilde{s}_{p}\right\} \\
&  \leq\left(  1+O\left(  \frac{\ln n}{n^{1/8a}}\right)  \right)  \max
_{p\in\left[  2,\overline{p}_{n}\right]  }\left\{  \bar{s}_{p}\right\}
+O_{\mathbb{P}}\left(  1\right)  .
\end{align*}
Define now%
\[
\mathsf{V}_{\Delta}\left(  p\right)  =\left(  2\sum_{j=1}^{p}K_{jp}%
^{2}\right)  ^{1/2},\text{\quad}\mathsf{s}_{p}=\frac{\sum_{j=1}^{p}%
K_{jp}\left(  \omega_{j}^{2}-1\right)  }{\mathsf{V}_{\Delta}\left(  p\right)
},
\]
which is such that%
\begin{align*}
&  \left\vert \bar{s}_{p}-\mathsf{s}_{p}\right\vert \leq\left\vert
\mathfrak{e}_{1p}\right\vert +\left\vert \mathfrak{e}_{2p}\right\vert \text{
where}\\
&  \text{ }\mathfrak{e}_{1p}=\left(  \frac{\sigma^{4}\mathsf{V}_{\Delta
}\left(  p\right)  }{\sigma^{4}V_{\Delta}\left(  p\right)  }-1\right)
\sigma^{4}\mathsf{s}_{p},\\
&  \text{ }\mathfrak{e}_{2p}=\frac{\sum_{j=1}^{p}\left\{  \left(  1-\frac
{j}{n}\right)  \operatorname*{Var}\left(  \eta_{jt}\right)  -\sigma
^{4}\right\}  K_{jp}\omega_{j}^{2}-\sigma^{4}\sum_{j=1}^{p}\frac{j}{n}K_{jp}%
}{\sigma^{4}V_{\Delta}\left(  p\right)  }.
\end{align*}
Since $K^{\prime}\left(  \cdot\right)  $ is continuous on $\left[  0,1\right]
$, the Weierstrass Theorem implies it can be uniformly approximated with a
sequence of polynomial function. Hence (\ref{kjdifcov}), Assumption
\ref{Kernel} and the LIL for weighted sums in Li and Tomkins (1996) gives that%
\[
\limsup_{p\rightarrow\infty}\frac{\left\vert \mathsf{V}_{\Delta}\left(
p\right)  \mathsf{s}_{p}\right\vert }{p^{1/2}\left(  2\ln\ln p\right)  ^{1/2}%
}\leq\left(  2\int K^{4}\left(  t\right)  dt\right)  ^{1/2}\text{, almost
surely.}%
\]
Since, under Assumption \ref{Kernel}, $\mathsf{V}_{\Delta}\left(  p\right)
/p^{1/2}\rightarrow\left(  2\int K^{4}\left(  t\right)  dt\right)  ^{1/2}$ by
convergence of Riemann sums, this gives%
\begin{equation}
\sup_{p\in\left[  2,\overline{p}_{n}\right]  }\left\vert \mathsf{s}%
_{p}\right\vert \leq\left(  2\ln\ln\overline{p}_{n}\right)  ^{1/2}\left(
1+o_{\mathbb{P}}\left(  1\right)  \right)  . \label{LIL}%
\end{equation}
Observe also that Lemma \ref{Ordersums}-(ii), $\overline{p}_{n}=o\left(
n^{1/2}\right)  $, and Assumption \ref{Kernel} give uniformly in $p\in\left[
1,\overline{p}_{n}\right]  $
\[
\left\vert \frac{\mathsf{V}_{\Delta}\left(  p\right)  }{V_{\Delta}\left(
p\right)  }-1\right\vert \leq C\left(  \frac{1}{p}\sum_{j=1}^{p}\frac{j^{2}%
}{n^{2}}K_{jp}^{2}\right)  ^{1/2}=o\left(  \frac{1}{n^{1/2}}\right)  .
\]
Hence%
\[
\max_{p\in\left[  2,\overline{p}_{n}\right]  }\left\vert \mathfrak{e}%
_{1p}\right\vert =o_{\mathbb{P}}\left(  \left(  \frac{2\ln\ln\overline{p}_{n}%
}{n}\right)  ^{1/2}\right)  =o_{\mathbb{P}}\left(  1\right)  .
\]
Now, for $\max_{p\in\left[  2,\overline{p}_{n}\right]  }\left\vert
\mathfrak{e}_{2p}\right\vert $, we have by Lemmas \ref{Ordersums}-(ii) and
\ref{Sumvareta}, $\overline{p}_{n}=o\left(  n^{1/2}\right)  $, and Assumption
\ref{Kernel},%
\[
\max_{p\in\left[  2,\overline{p}_{n}\right]  }\left\vert \mathfrak{e}%
_{2p}\right\vert \leq C\left\{  \sum_{j=1}^{\overline{p}_{n}}\left\vert
\operatorname*{Var}\left(  \eta_{jt}\right)  -\sigma^{4}\right\vert \omega
_{j}^{2}+\frac{1}{n}\sum_{j=1}^{\overline{p}_{n}}j\omega_{j}^{2}%
+\frac{\overline{p}_{n}^{3/2}}{n}\right\}  =O_{\mathbb{P}}\left(  1\right)
+O_{\mathbb{P}}\left(  \frac{\overline{p}_{n}^{2}}{n}\right)  =O_{\mathbb{P}%
}\left(  1\right)  .
\]
Hence $\max_{p\in\left[  2,\overline{p}_{n}\right]  }\left\vert \bar{s}%
_{p}-\mathsf{s}_{p}\right\vert =O_{\mathbb{P}}\left(  1\right)  $ and
substituting in the bounds for $\mathbb{P}\left(  \max_{p\in\left[
2,\overline{p}_{n}\right]  }\left\{  \check{s}_{p}\right\}  \geq\gamma
_{n}^{\prime}\right)  $ and $\widetilde{\mathcal{M}}$ above gives, by
(\ref{Gam}), $\gamma_{n}^{\prime}=\gamma_{n}\left(  1+\epsilon/2\right)
/\left(  1+\epsilon\right)  $, $\gamma_{n}^{\prime}\geq\left(  2\ln
\ln\overline{p}_{n}\right)  ^{1/2}\left(  1+\epsilon/3\right)  $ and
(\ref{LIL})%
\begin{align}
\mathbb{P}\left(  \max_{p\in\left[  2,\overline{p}_{n}\right]  }\left\{
\check{s}_{p}\right\}  \geq\gamma_{n}^{\prime}\right)   &  =\mathbb{P}\left(
\left(  1+O\left(  \frac{\ln n}{n^{1/8a}}\right)  \right)  \max_{p\in\left[
2,\overline{p}_{n}\right]  }\left\{  \mathsf{s}_{p}\right\}  +O_{\mathbb{P}%
}\left(  1\right)  \geq\gamma_{n}^{\prime}-\epsilon\right)  +o\left(  1\right)
\nonumber\\
&  \leq\mathbb{P}\left(  \max_{p\in\left[  2,\overline{p}_{n}\right]
}\left\{  \mathsf{s}_{p}\right\}  \geq\left(  2\ln\ln\overline{p}_{n}\right)
^{1/2}\left(  1+\epsilon/3\right)  \right)  +o\left(  1\right) \nonumber\\
&  =o\left(  1\right)  . \label{DE56}%
\end{align}

\textbf{Step 3: Conclusion}. Propositions \ref{Esti} and \ref{Covesti}, Lemma
\ref{Ordersums} and $\overline{p}_{n}=O\left(  n^{1/2}\right)  $, the
expression of $\check{S}_{p}$ and $\check{s}_{p}$ in (\ref{Pseudomax}) and
(\ref{Check2tilde}) gives%
\begin{align*}
&  \max_{p\in\left[  2,\overline{p}_{n}\right]  }\frac{\left(  \widehat{S}%
_{p}-\widehat{S}_{1}\right)  /\widehat{R}_{0}^{2}-E_{\Delta}\left(  p\right)
}{V_{\Delta}\left(  p\right)  }=\max_{p\in\left[  2,\overline{p}_{n}\right]
}\frac{\left(  \widehat{S}_{p}-\widehat{S}_{1}\right)  -\widehat{R}_{0}%
^{2}E_{\Delta}\left(  p\right)  }{\widehat{R}_{0}^{2}V_{\Delta}\left(
p\right)  }\\
&  \text{ }=\left(  1+o_{\mathbb{P}}\left(  1\right)  \right)  \max
_{p\in\left[  2,\overline{p}_{n}\right]  }\frac{\left(  \widetilde{S}%
_{p}-\widetilde{S}_{1}\right)  -R_{0}^{2}E_{\Delta}\left(  p\right)  }%
{R_{0}^{2}V_{\Delta}\left(  p\right)  }+O_{\mathbb{P}}\left(  1+\overline
{p}_{n}^{1/2}\left(  \widehat{R}_{0}^{2}-R_{0}^{2}\right)  \right) \\
&  \text{ }=\left(  1+o_{\mathbb{P}}\left(  1\right)  \right)  \max
_{p\in\left[  2,\overline{p}_{n}\right]  }\left\{  \check{s}_{p}\right\}
+O_{\mathbb{P}}\left(  1\right)  .
\end{align*}
Hence (\ref{DE56}) gives, since $\gamma_{n}-\gamma_{n}^{\prime}\rightarrow
+\infty$,%
\[
\mathbb{P}\left(  \max_{p\in\left[  2,\overline{p}_{n}\right]  }\frac{\left(
\widehat{S}_{p}-\widehat{S}_{1}\right)  /\widehat{R}_{0}^{2}-E_{\Delta}\left(
p\right)  }{V_{\Delta}\left(  p\right)  }\geq\gamma_{n}\right)  \leq
\mathbb{P}\left(  \max_{p\in\left[  2,\overline{p}_{n}\right]  }\left\{
\check{s}_{p}\right\}  \geq\gamma_{n}^{\prime}\right)  +o\left(  1\right)
=o\left(  1\right)  .
\]
This ends the proof of the Proposition.$\hfill\square$

\subsection{Proof of Propositions \ref{MeanH1} and \ref{VarH1}}

\noindent When studying the mean and variance of $\widetilde{S}_{p}$, we make
use of Theorem 2.3.2 in Brillinger (2001) which implies in particular that,
for any real zero-mean random variables $Z_{1},\ldots,Z_{4}$,
\begin{align}
&  \operatorname*{Var}\left(  Z_{1}Z_{2},Z_{3}Z_{4}\right)
=\operatorname*{Var}(Z_{1},Z_{3})\operatorname*{Var}(Z_{2},Z_{4}%
)+\operatorname*{Var}(Z_{1},Z_{4})\operatorname*{Var}(Z_{2},Z_{3})\nonumber\\
&  +\operatorname*{Cum}\left(  Z_{1},Z_{2},Z_{3},Z_{4}\right)  . \label{Brill}%
\end{align}
Note that Assumption \ref{Reg} and Theorem \ref{XW11} imply that%
\begin{equation}
\sup_{n,q\in\left[  2,8\right]  }\sum_{t_{2},\ldots,t_{q}=-\infty}^{\infty
}\left\vert \Gamma_{n}\left(  0,t_{2},\ldots,t_{q}\right)  \right\vert
<\infty. \label{Sumcum}%
\end{equation}

\subsubsection{\textbf{Proof of Proposition \ref{MeanH1}}}

(\ref{Brill}) yields
\begin{align*}
&  \mathbb{E}\left[  \widetilde{R}_{j}^{2}\right]  =\frac{1}{n^{2}}\sum
_{t_{1},t_{2}=1}^{n-j}\mathbb{E}\left[  u_{t_{1}}u_{t_{1}+j}u_{t_{2}}%
u_{t_{2}+j}\right] \\
&  \text{ }=\frac{1}{n^{2}}\sum_{t_{1},t_{2}=1}^{n-j}\left(  R_{j}%
^{2}+R_{t_{2}-t_{1}}^{2}+R_{t_{2}-t_{1}+j}R_{t_{2}-t_{1}-j}+\Gamma\left(
0,j,t_{2}-t_{1},t_{2}-t_{1}+j\right)  \right)  ,
\end{align*}
where
\begin{align*}
\sum_{t_{1},t_{2}=1}^{n-j}R_{t_{2}-t_{1}}^{2}  &  =(n-j)R_{0}^{2}+2\sum
_{\ell=1}^{n-j-1}(n-j-\ell)R_{\ell}^{2},\\
\sum_{t_{1},t_{2}=1}^{n-j}R_{t_{2}-t_{1}+j}R_{t_{2}-t_{1}-j}  &
=(n-j)R_{j}^{2}+2\sum_{\ell=1}^{n-j-1}(n-j-\ell)R_{\ell+j}R_{\ell-j},\\
\sum_{t_{1},t_{2}=1}^{n-j}\Gamma\left(  0,j,t_{2}-t_{1},t_{2}-t_{1}+j\right)
&  =\sum_{\ell=-n+j+1}^{n-j-1}\left(  n-j-|\ell|\right)  \Gamma\left(
0,j,\ell,\ell+j\right)  .
\end{align*}
Set $k_{j}=K^{2}\left(  j/p\right)  $ to prove the first equality and
$k_{j}=K^{2}\left(  j/p\right)  /\tau_{j}^{2}$ for the second. Note that
Assumptions \ref{Kernel} and \ref{Reg} give, in both case, $\max_{j\in\left[
1,n-1\right]  }k_{j}\leq C$ and $k_{j}\geq C\mathbb{I}\left(  j\leq
p/2\right)  $. The equalities above give
\begin{align}
&  \mathbb{E}\left[  \sum_{j=1}^{n-1}k_{j}\widetilde{R}_{j}^{2}\right]
-R_{0}^{2}\sum_{j=1}^{n-1}\left(  1-\frac{j}{n}\right)  k_{j}\nonumber\\
&  \text{ }=n\sum_{j=1}^{n-1}\left(  \left(  1-\frac{j}{n}\right)  ^{2}%
+\frac{1}{n}\left(  1-\frac{j}{n}\right)  \right)  k_{j}R_{j}^{2}%
\label{EtildeSp}\\
&  \text{ }+2\sum_{j=1}^{n-1}k_{j}\sum_{\ell=1}^{n-j-1}\left(  1-\frac{j+\ell
}{n}\right)  \left(  R_{\ell}^{2}+R_{\ell+j}R_{\ell-j}\right) \nonumber\\
&  \text{ }+\sum_{j=1}^{n-1}k_{j}\sum_{\ell=-n+j+1}^{n-j-1}\left(
1-\frac{j+|\ell|}{n}\right)  \Gamma\left(  0,j,\ell,\ell+j\right)  .\nonumber
\end{align}
We start with the item $R_{0}^{2}\sum_{j=1}^{n-1}\left(  1-\frac{j}{n}\right)
k_{j}$, which is equal to $R_{0}^{2}E\left(  p\right)  $ when $k_{j}%
=K^{2}\left(  j/p\right)  $, that is when proving the first equality. When
$k_{j}=K^{2}\left(  j/p\right)  /\tau_{j}^{2}$, (\ref{Tau2sig}) gives, under
Assumptions \ref{Kernel} and \ref{Reg},%
\[
\left\vert R_{0}^{2}\sum_{j=1}^{n-1}\left(  1-\frac{j}{n}\right)
k_{j}-E\left(  p\right)  \right\vert \leq C\sum_{j=1}^{p}\left\vert \tau
_{j}^{2}-R_{0}^{2}\right\vert \leq C\sum_{j=1}^{\infty}j^{-6}%
\]
so that $R_{0}^{2}\sum_{j=1}^{n-1}\left(  1-j/n\right)  k_{j}\geq E\left(
p\right)  -C^{\prime}$.

Let us now turn to the other items. The lower bound$k_{j}\geq CI(j\leq p/2)$
gives that (\ref{EtildeSp}) is larger than $Cn\sum_{j=1}^{p/2}R_{j}^{2}$. To
bound the remaining terms in (\ref{EtildeSp}), we note that by Assumptions
\ref{Kernel}, \ref{Reg} and (\ref{Sumcum}),%
\[
\left\vert \sum_{j=1}^{n-1}k_{j}\sum_{\ell=1}^{n-j-1}\left(  1-\frac{j+\ell
}{n}\right)  R_{\ell}^{2}\right\vert \leq C\sum_{j=1}^{n-1}\mathbb{I}(j\leq
p)\times\sum_{j=1}^{\infty}R_{j}^{2}\leq Cp\sum_{j=1}^{\infty}R_{j}%
^{2}=o(n)\sum_{j=1}^{\infty}R_{j}^{2},
\]%
\[
\left\vert \sum_{j=1}^{n-1}k_{j}\sum_{\ell=1}^{n-j-1}\left(  1-\frac{j+\ell
}{n}\right)  R_{\ell+j}R_{\ell-j}\right\vert \leq C\sum_{j=1}^{+\infty}%
\sum_{\ell=1}^{+\infty}\left\vert R_{\ell+j}R_{\ell-j}\right\vert \leq
C\left(  \sum_{j=0}^{\infty}|R_{j}|\right)  ^{2}\leq C,
\]%
\[
\left\vert \sum_{j=1}^{n-1}k_{j}\sum_{\ell=-n+j+1}^{n-j-1}\left(
1-\frac{j+\ell}{n}\right)  \Gamma\left(  0,j,\ell,\ell+j\right)  \right\vert
\leq C\sum_{t_{2},t_{3},t_{4}=-\infty}^{\infty}\left\vert \Gamma(0,t_{2}%
,t_{3},t_{4})\right\vert \leq C
\]
uniformly with respect to $p\in\left[  1,\overline{p}_{n}\right]  $.
Substituting these bounds in the equality above establishes the proposition.
\hspace*{\fill}$\Box$

\subsubsection{\textbf{Proof of Proposition \ref{VarH1}}}

Let $f$ be the spectral density of the alternative. Using (\ref{Sumcum}), we
obtain%
\begin{equation}
\sup_{\lambda\in\lbrack-\pi,\pi]}\left\vert f\left(  \lambda\right)
\right\vert \leq C\text{\quad and\quad}\sum_{j=1}^{\infty}R_{j}^{2}\leq C
\label{L2}%
\end{equation}
because $\sup_{\lambda\in\lbrack-\pi,\pi]}\left\vert f\left(  \lambda\right)
\right\vert \leq\left(  |R_{0}|+2\sum_{j=1}^{\infty}|R_{j}|\right)  /(2\pi)$
and $\sum_{j=1}^{\infty}R_{j}^{2}\leq\left(  \sum_{j=1}^{\infty}%
|R_{j}|\right)  ^{2}$. We recall that $\widetilde{R}_{j}=\sum_{t=1}^{n-j}%
u_{t}u_{t+j}/n$\ and define $\overline{R}_{j}=\mathbb{E}\left[  \widetilde
{R}_{j}\right]  =\left(  1-j/n\right)  R_{j}$. Set $k_{j}=K^{2}\left(
j/p\right)  $ to prove the first equality and $k_{j}=K^{2}\left(  j/p\right)
/\tau_{j}^{2}$ for the second. Note that Assumptions \ref{Kernel} and
\ref{Reg} give, in both case, $k_{j}\leq C\mathbb{I}\left(  j\leq p\right)  $.
To avoid notation burdens, redefine $\widetilde{S}_{p}$ as $\sum_{j=1}%
^{n-1}k_{j}\widetilde{R}_{j}^{2}$. Define $D_{j}=\widetilde{R}_{j}%
-\overline{R}_{j}$. We have $\mathbb{E}\left[  D_{j}\right]  =0$ and
$\widetilde{S}_{p}=n\sum_{j=1}^{n-1}k_{j}\overline{R}_{j}^{2}+2n\sum
_{j=1}^{n-1}k_{j}\overline{R}_{j}D_{j}+n\sum_{j=1}^{n-1}k_{j}D_{j}^{2}$. The
inequality $(a+b)^{2}\leq2a^{2}+2b^{2}$ implies that
\begin{equation}
\operatorname*{Var}\left(  \widetilde{S}_{p}\right)  \leq4\operatorname*{Var}%
\left(  n\sum_{j=1}^{n-1}k_{j}\overline{R}_{j}\widetilde{R}_{j}\right)
+2\operatorname*{Var}\left(  n\sum_{j=1}^{n-1}k_{j}D_{j}^{2}\right)  .
\label{VarH1.1}%
\end{equation}
By identity (\ref{Brill}),
\[
\operatorname*{Var}\left(  n\sum_{j=1}^{n-1}k_{j}\overline{R}_{j}\widetilde
{R}_{j}\right)  =\sum_{j_{1},j_{2}=1}^{n-1}k_{j_{1}}k_{j_{2}}\overline
{R}_{j_{1}}\overline{R}_{j_{2}}\sum_{t_{1}=1}^{n-j_{1}}\sum_{t_{2}=1}%
^{n-j_{2}}\operatorname*{Cov}\left(  u_{t_{1}}u_{t_{1}+j_{1}},u_{t_{2}%
}u_{t_{2}+j_{2}}\right)  \leq V_{1}+K_{1}%
\]
with
\begin{align*}
V_{1}  &  =\left\vert \sum_{j_{1},j_{2}=1}^{n-1}k_{j_{1}}k_{j_{2}}\overline
{R}_{j_{1}}\overline{R}_{j_{2}}\sum_{t_{1}=1}^{n-j_{1}}\sum_{t_{2}=1}%
^{n-j_{2}}\left(  R_{t_{2}-t_{1}}R_{t_{2}-t_{1}+j_{2}-j_{1}}+R_{t_{2}%
-t_{1}-j_{1}}R_{t_{2}-t_{1}+j_{2}}\right)  \right\vert ,\\
K_{1}  &  =\left\vert \sum_{j_{1},j_{2}=1}^{n-1}k_{j_{1}}k_{j_{2}}\overline
{R}_{j_{1}}\overline{R}_{j_{2}}\sum_{t_{1}=1}^{n-j_{1}}\sum_{t_{2}=1}%
^{n-j_{2}}\Gamma\left(  t_{1},t_{1}+j_{1},t_{2},t_{2}+j_{2}\right)
\right\vert .
\end{align*}
The second term on the right of (\ref{VarH1.1}) is, up to a multiplicative
constant, equal to
\[
\operatorname*{Var}\left(  n\sum_{j=1}^{n-1}k_{j}D_{j}^{2}\right)  =n^{2}%
\sum_{j_{1},j_{2}=1}^{n-1}k_{j_{1}}k_{j_{2}}\operatorname*{Cov}\left(
D_{j_{1}}^{2},D_{j_{2}}^{2}\right)  .
\]
Applying (\ref{Brill}) twice we obtain
\begin{align*}
\lefteqn{\operatorname*{Cov}\left(  D_{j_{1}}^{2},D_{j_{2}}^{2}\right)  }\\
&  \text{ }=\frac{1}{n^{4}}\sum_{t_{1},t_{2}=1}^{n-j_{1}}\sum_{t_{3},t_{4}%
=1}^{n-j_{2}}\operatorname*{Cov}\left[  \prod_{q=1}^{2}\left(  u_{t_{q}%
}u_{t_{q}+j_{1}}-\mathbb{E}[u_{t_{q}}u_{t_{q}+j_{1}}]\right)  ,\prod_{q=3}%
^{4}\left(  u_{t_{q}}u_{t_{q}+j_{2}}-\mathbb{E}[u_{t_{q}}u_{t_{q}+j_{2}%
}]\right)  \right] \\
&  \text{ }=\frac{1}{n^{4}}\sum_{t_{1},t_{2}=1}^{n-j_{1}}\sum_{t_{3},t_{4}%
=1}^{n-j_{2}}\left[  \operatorname*{Cov}\left(  u_{t_{1}}u_{t_{1}+j_{1}%
},u_{t_{3}}u_{t_{3}+j_{2}}\right)  \operatorname*{Cov}\left(  u_{t_{2}%
}u_{t_{2}+j_{1}},u_{t_{4}}u_{t_{4}+j_{2}}\right)  \right. \\
&  \left.
\;\;\;\;\;\;\;\;\;\;\;\;\;\;\;\;\;\;\;\;\;\;\;\;\;\;\;\;\;\;\;\;\;\;\;\;\;\;\;\;+\operatorname*{Cov}%
\left(  u_{t_{1}}u_{t_{1}+j_{1}},u_{t_{4}}u_{t_{4}+j_{2}}\right)
\operatorname*{Cov}\left(  u_{t_{2}}u_{t_{2}+j_{1}},u_{t_{3}}u_{t_{3}+j_{2}%
}\right)  \right] \\
&  +\frac{1}{n^{4}}\sum_{t_{1},t_{2}=1}^{n-j_{1}}\sum_{t_{3},t_{4}=1}%
^{n-j_{2}}\operatorname*{Cum}\left(  u_{t_{1}}u_{t_{1}+j_{1}},u_{t_{2}%
}u_{t_{2}+j_{1}},u_{t_{3}}u_{t_{3}+j_{2}},u_{t_{4}}u_{t_{4}+j_{2}}\right) \\
&  \text{ }=\frac{2}{n^{4}}\left(  \sum_{t_{1}=1}^{n-j_{1}}\sum_{t_{2}%
=1}^{n-j_{2}}\left(  R_{t_{2}-t_{1}}R_{t_{2}-t_{1}+j_{2}-j_{1}}+R_{t_{2}%
-t_{1}-j_{1}}R_{t_{2}-t_{1}+j_{2}}+\Gamma(t_{1},t_{1}+j_{1},t_{2},t_{2}%
+j_{2})\right)  \right)  ^{2}\\
&  +\frac{1}{n^{4}}\sum_{t_{1},t_{2}=1}^{n-j_{1}}\sum_{t_{3},t_{4}=1}%
^{n-j_{2}}\operatorname*{Cum}\left(  u_{t_{1}}u_{t_{1}+j_{1}},u_{t_{2}%
}u_{t_{2}+j_{1}},u_{t_{3}}u_{t_{3}+j_{2}},u_{t_{4}}u_{t_{4}+j_{2}}\right)  .
\end{align*}
Since $(a+b+c)^{2}\leq3(a^{2}+b^{2}+c^{2})$, we can write $\operatorname*{Var}%
\left(  n\sum_{j=1}^{n-1}k_{j}D_{j}^{2}\right)  \leq6V_{2}+K_{2}%
+6K_{2}^{\prime}$ with
\begin{align*}
\lefteqn{V_{2}=\frac{1}{n^{2}}\sum_{j_{1},j_{2}=1}^{n-1}k_{j_{1}}k_{j_{2}%
}\left(  \left(  \sum_{t_{1}=1}^{n-j_{1}}\sum_{t_{2}=1}^{n-j_{2}}%
R_{t_{2}-t_{1}}R_{t_{2}-t_{1}+j_{2}-j_{1}}\right)  ^{2}+\left(  \sum_{t_{1}%
=1}^{n-j_{1}}\sum_{t_{2}=1}^{n-j_{2}}R_{t_{2}-t_{1}-j_{1}}R_{t_{2}-t_{1}%
+j_{2}}\right)  ^{2}\right)  ,}\\
&  K_{2}=\left\vert \frac{1}{n^{2}}\sum_{j_{1},j_{2}=1}^{n-1}k_{j_{1}}%
k_{j_{2}}\sum_{t_{1},t_{2}=1}^{n-j_{1}}\sum_{t_{3},t_{4}=1}^{n-j_{2}%
}\operatorname*{Cum}\left(  u_{t_{1}}u_{t_{1}+j_{1}},u_{t_{2}}u_{t_{2}+j_{1}%
},u_{t_{3}}u_{t_{3}+j_{2}},u_{t_{4}}u_{t_{4}+j_{2}}\right)  \right\vert ,\\
&  K_{2}^{\prime}=\frac{1}{n^{2}}\sum_{j_{1},j_{2}=1}^{n-1}k_{j_{1}}k_{j_{2}%
}\left(  \sum_{t_{1}=1}^{n-j_{1}}\sum_{t_{2}=1}^{n-j_{2}}\Gamma\left(
t_{1},t_{1}+j_{1},t_{2},t_{2}+j_{2}\right)  \right)  ^{2},
\end{align*}
Substituting in (\ref{VarH1.1}) shows that the proposition holds if the
following inequalities hold:
\[
V_{1}\leq Cn\sum_{j=1}^{p}R_{j}^{2},\quad V_{2}\leq Cp,\quad K_{1}\leq C,\quad
K_{2}^{\prime}\leq C,\quad K_{2}\leq C\frac{p^{2}}{n}.
\]
We establish these inequalities in five steps.

\textit{Step 1: bound for }$V_{1}$\textit{.} We note that $|\overline{R}%
_{j}|\leq|R_{j}|$ and that under Assumption \ref{Kernel}, $0\leq k_{j}\leq C $
for all $j$. Using a covariance spectral representation $R_{j}=\int_{-\pi
}^{\pi}\exp(\pm ij\lambda)f(\lambda)d\lambda$, the Cauchy-Schwarz inequality
and (\ref{L2}), we obtain by Assumption \ref{Kernel}%
\begin{align*}
\lefteqn{\left\vert \sum_{j_{1},j_{2}=1}^{n-1}k_{j_{1}}k_{j_{2}}\overline
{R}_{j_{1}}\overline{R}_{j_{2}}\sum_{t_{1}=1}^{n-j_{1}}\sum_{t_{2}=1}%
^{n-j_{2}}R_{t_{2}-t_{1}}R_{t_{2}-t_{1}+j_{2}-j_{1}}\right\vert }\\
&  \text{ }=\int_{-\pi}^{\pi}\int_{-\pi}^{\pi}\left\vert \sum_{j=1}^{n-1}%
k_{j}\overline{R}_{j}\sum_{t=1}^{n-j}\text{e}^{it\lambda_{1}}\text{e}%
^{i(t+j)\lambda_{2}}\right\vert ^{2}f(\lambda_{1})f(\lambda_{2})d\lambda
_{1}d\lambda_{2}\\
&  \text{ }\leq\left(  \sup_{\lambda\in\lbrack-\pi,\pi]}|f(\lambda)|\right)
^{2}\int_{-\pi}^{\pi}\int_{-\pi}^{\pi}\sum_{j_{1},j_{2}=1}^{n-1}k_{j_{1}%
}\overline{R}_{j_{1}}k_{j_{2}}\overline{R}_{j_{2}}\sum_{t_{1}=1}^{n-j_{1}}%
\sum_{t_{2}=1}^{n-j_{2}}\text{e}^{it_{1}\lambda_{1}}\text{e}^{i(t_{1}%
+j_{1})\lambda_{2}}\text{e}^{-it_{2}\lambda_{1}}\text{e}^{-i(t_{2}%
+j_{2})\lambda_{2}}d\lambda_{1}d\lambda_{2}\\
&  \text{ }\leq C\sum_{j=1}^{n-1}(n-j)k_{j}^{2}\overline{R}_{j}^{2}\leq
Cn\sum_{j=1}^{p}R_{j}^{2},
\end{align*}%
\begin{align*}
\lefteqn{\left\vert \sum_{j_{1},j_{2}=1}^{n-1}k_{j_{1}}k_{j_{2}}\overline
{R}_{j_{1}}\overline{R}_{j_{2}}\sum_{t_{1}=1}^{n-j_{1}}\sum_{t_{2}=1}%
^{n-j_{2}}R_{t_{2}-t_{1}-j_{1}}R_{t_{2}-t_{1}+j_{2}}\right\vert }\\
&  \text{ }=\left\vert \int_{-\pi}^{\pi}\int_{-\pi}^{\pi}\sum_{j_{1}=1}%
^{n-1}k_{j_{1}}\overline{R}_{j_{1}}\sum_{t_{1}=1}^{n-j_{1}}\text{e}%
^{-i(t_{1}+j_{1})\lambda_{1}}\text{e}^{-it_{1}\lambda_{2}}\times\sum_{j_{2}%
=1}^{n-1}k_{j_{2}}\overline{R}_{j_{2}}\sum_{t_{2}=1}^{n-j_{2}}\text{e}%
^{it_{2}\lambda_{1}}\text{e}^{i(t_{2}+j_{2})}f(\lambda_{1})f(\lambda
_{2})d\lambda_{1}d\lambda_{2}\right\vert \\
&  \text{ }\leq\int_{-\pi}^{\pi}\int_{-\pi}^{\pi}\left\vert \sum_{j=1}%
^{n-1}k_{j}\overline{R}_{j}\sum_{t=1}^{n-j}\text{e}^{it\lambda_{1}}%
\text{e}^{i(t+j)\lambda_{2}}\right\vert ^{2}f(\lambda_{1})f(\lambda
_{2})d\lambda_{1}d\lambda_{2}\leq Cn\sum_{j=1}^{p}R_{j}^{2}%
\end{align*}
This establishes the bound for $V_{1}$.

\textit{Step 2: bound for }$V_{2}$\textit{.} We define $t_{2}=t_{1}%
+t_{2}^{\prime}$, $j_{2}=j_{1}+j_{2}^{\prime}$. By Assumption \ref{Kernel} and
by (\ref{Sumcum}),
\begin{align*}
\lefteqn{\frac{1}{n^{2}}\sum_{j_{1},j_{2}=1}^{n-1}k_{j_{1}}k_{j_{2}}\left(
\sum_{t_{1}=1}^{n-j_{1}}\sum_{t_{2}=1}^{n-j_{2}}R_{t_{2}-t_{1}}R_{t_{2}%
-t_{1}-j_{1}+j_{2}}\right)  ^{2}}\\
&  \text{ }\leq\frac{C}{n^{2}}\sum_{j_{1}=1}^{n-1}K^{2}(j_{1}/p)\sum
_{j_{2}\prime=-\infty}^{\infty}\left(  n\sum_{t_{2}\prime=-\infty}^{+\infty
}\left\vert R_{t_{2}\prime}R_{t_{2}\prime+j_{2}\prime}\right\vert \right)
^{2}\\
&  \text{ }\leq Cp\times\left(  \sum_{j_{2},t_{1},t_{2}=-\infty}^{\infty
}\left\vert R_{t_{1}}R_{t_{1}+j_{2}}R_{t_{2}}R_{t_{2}+j_{2}}\right\vert
\right)  \leq Cp\left(  \sum_{t=-\infty}^{\infty}|R_{t}|\right)  ^{4}\leq Cp,
\end{align*}%
\begin{align*}
\lefteqn{\frac{1}{n^{2}}\sum_{j_{1},j_{2}=1}^{n-1}k_{j_{1}}k_{j_{2}}\left(
\sum_{t_{1}=1}^{n-j_{1}}\sum_{t_{2}=1}^{n-j_{2}}R_{t_{2}-t_{1}-j_{1}}%
R_{t_{2}-t_{1}+j_{2}}\right)  ^{2}}\\
&  \leq\frac{C}{n^{2}}\sum_{j_{1}=1}^{n-1}K^{2}(j_{1}/p)\sum_{j_{2}%
\prime=-\infty}^{\infty}\left(  n\sum_{t_{2}\prime=-\infty}^{+\infty
}\left\vert R_{t_{2}\prime-j_{1}}R_{t_{2}\prime+j_{1}+j_{2}\prime}\right\vert
\right)  ^{2}\\
&  \text{ }\leq Cp\sum_{j_{2}^{\prime},t_{1},t_{2}=-\infty}^{\infty}\left\vert
R_{t_{1}-j_{1}}R_{t_{1}+j_{1}+j_{2}^{\prime}}R_{t_{2}-j_{1}}R_{t_{2}%
+j_{1}+j_{2}^{\prime}}\right\vert \leq Cp\sum_{j,t_{1},t_{2}=-\infty}^{\infty
}\left\vert R_{t_{1}}R_{t_{1}+j}R_{t_{2}}R_{t_{2}+j}\right\vert \\
&  \text{ }\leq Cp\left(  \sum_{t=-\infty}^{\infty}|R_{t}|\right)  ^{4}\leq
Cp,
\end{align*}
therefore $V_{2}\leq Cp$.

\textit{Step 3: bound for }$K_{1}$\textit{.} Define $t_{2}=t_{1}+t$.
Assumption \ref{Kernel}, and (\ref{Sumcum}) yield
\[
K_{1}\leq C\sum_{j_{1},j_{2}=1}^{p}\sum_{t=-\infty}^{\infty}\left\vert
\Gamma(0,j_{1},t,t+j_{2})\right\vert \leq\sum_{t_{1},t_{2},t_{3}=-\infty
}^{\infty}\left\vert \Gamma(0,t_{1},t_{2},t_{3})\right\vert .
\]

\textit{Step 4: bound for }$K_{2}^{\prime}$\textit{.} (\ref{Sumcum}) gives
\begin{align*}
&  K_{2}^{\prime}\leq\frac{1}{n^{2}}\sum_{j_{1},j_{2}=1}^{n-1}k_{j_{1}%
}k_{j_{2}}\left(  \sum_{t_{1}=1}^{n-j_{1}}\sum_{t_{2}=1}^{n-j_{2}}\left\vert
\Gamma\left(  0,j_{1},t_{2}-t_{1},t_{2}-t_{1}+j_{2}\right)  \right\vert
\right)  ^{2}\\
&  \leq C\sum_{j_{1},j_{2}=1}^{+\infty}\left(  \sum_{t=-\infty}^{\infty
}\left\vert \Gamma(0,j_{1},t,t+j_{2})\right\vert \right)  ^{2}\\
&  =C\sum_{j_{1},j_{2}=1}^{+\infty}\sum_{t_{1},t_{2}=-\infty}^{\infty
}\left\vert \Gamma(0,j_{1},t_{1},t_{1}+j_{2})\Gamma(0,j_{1},t_{2},t_{2}%
+j_{2})\right\vert \\
&  \leq C\left(  \sum_{t_{2},t_{3},t_{4}=-\infty}^{\infty}\left\vert
\Gamma(0,t_{2},t_{3},t_{4})\right\vert \right)  ^{2}\leq C.
\end{align*}

\textit{Step 5: bound for }$K_{2}$\textit{.}\quad Bounding $K_{2}$ requires
additional notation. First set $t_{5}=t_{1}+j_{1}$, $t_{6}=t_{2}+j_{1}$,
$t_{7}=t_{3}+j_{2}$ and $t_{8}=t_{4}+j_{2}$, and note that $t_{5},\ldots
,t_{8}$ depend upon $t_{1},\ldots,t_{4}$ and $j_{1},j_{2}$ only. For a
partition $B=\{B_{\ell},\ell=1,\ldots,d_{B}\}$ of $\{1,\ldots,8\}$, define
$d_{B}=\operatorname*{Card}B$, $\Gamma_{B}(t_{1},\ldots,t_{8})=\prod_{\ell
=1}^{d_{B}}\operatorname*{Cum}\left(  u_{t_{q}},q\in B_{\ell}\right)  $, and
recall that $\operatorname*{Cum}(u_{t})=Eu_{t}=0$. Then the largest $d_{B}$
yielding a non-vanishing $\Gamma_{B}$ is $d_{B}=4$. When $d_{B}=4$, $B$ is a
pairwise partition of $\{1,\ldots,8\}$ so that $\Gamma_{B}$ is a product of
covariances. Let $B$ be the set of indecomposable partitions of the two-way
table
\[%
\begin{array}
[c]{cc}%
1 & 5\\
2 & 6\\
3 & 7\\
4 & 8\\
&
\end{array}
,
\]
see Brillinger (2001, p. 20) for a definition. Then according to Brillinger
(2001, Theorem 2.3.2),
\begin{align*}
&  \operatorname*{Cum}\left(  u_{t_{1}}u_{t_{1}+j_{1}},u_{t_{2}}u_{t_{2}%
+j_{1}},u_{t_{3}}u_{t_{3}+j_{2}},u_{t_{4}}u_{t_{4}+j_{2}}\right) \\
&  \text{ }=\sum_{B\in\mathcal{B}}\Gamma_{B}(t_{1},\ldots,t_{8})=\sum
_{B\in\mathcal{B},d_{B}\leq3}\Gamma_{B}(t_{1},\ldots,t_{8})+\sum
_{B\in\mathcal{B},d_{B}=4}\Gamma_{B}(t_{1},\ldots,t_{8}).
\end{align*}
Some properties of partitions in $\mathcal{B}$ are as follows. Call $\{1,5\}$,
$\{2,6\}$, $\{3,7\}$ and $\{4,8\}$ fundamental pairs and say that a $B_{1}$ in
a partition $B$ breaks the pair $\{1,5\}$ if $\{1,5\}$ is not a subset of
$B_{1}$. Then partitions $B\in\mathcal{B}$ are such that each $B_{\ell}\in B$
must break a fundamental pair. Note that fundamental pairs play a symmetric
role. Since $t_{q+4}-t_{q}$ is $j_{1}$ or $j_{2}$ with vanishing $k_{j_{1}}$
or $k_{j_{2}}$ if $j_{1}$ or $j_{2}$ is larger than $p$, the indexes $t_{q}$
and $t_{q+4}$ of a fundamental pair also play a symmetric role in the
computations below. We now discuss the contribution to $K_{2}$ of partitions
of $\{1,\ldots,8\}$ according to the possible values $1,\ldots,4$ of $d_{B}$.
Due to symmetry, we only consider representative partitions for each case.

Under Assumption \ref{Kernel} and (\ref{Sumcum}), the case $d_{B}=1$ gives a
contribution to $K_{2}$ bounded by
\begin{align*}
\left\vert \frac{1}{n^{2}}\sum_{j_{1},j_{2}=1}^{n-1}k_{j_{1}}k_{j_{2}}%
\sum_{t_{1},t_{2}=1}^{n-j_{1}}\sum_{t_{3},t_{4}=1}^{n-j_{2}}\Gamma\left(
t_{1},\ldots,t_{8}\right)  \right\vert  &  \leq\frac{C}{n^{2}}\sum
_{t_{1},\ldots,t_{8}=-n}^{n}\left\vert \Gamma\left(  0,t_{2}-t_{1}%
,\ldots,t_{8}-t_{1}\right)  \right\vert \\
&  \leq\frac{C}{n}\sum_{t_{2}^{\prime},\ldots,t_{8}^{\prime}=-\infty}^{\infty
}\left\vert \Gamma\left(  0,t_{2}^{\prime},\ldots,t_{8}^{\prime}\right)
\right\vert \leq\frac{C}{n}.
\end{align*}

The case $d_{B}=2$ corresponds to $\{\operatorname*{Card}B_{1}%
,\operatorname*{Card}B_{2}\}$ being $\{2,6\}$, $\{3,5\}$ or $\{4,4\}$. These
cases are very similar and we limit ourselves to $\{2,6\}$ and $B_{1}%
=\{1,2\}$. The corresponding contribution to $K_{2}$ is bounded by
\begin{align*}
\lefteqn{\left\vert \frac{1}{n^{2}}\sum_{j_{1},j_{2}=1}^{n-1}k_{j_{1}}%
k_{j_{2}}\sum_{t_{1},t_{2}=1}^{n-j_{1}}\sum_{t_{3},t_{4}=1}^{n-j_{2}}%
\Gamma_{B}\left(  t_{1},\ldots,t_{8}\right)  \right\vert \leq\frac{C}{n^{2}%
}\sum_{t_{1},\ldots,t_{8}=-n}^{n}\left\vert \Gamma\left(  0,t_{2}%
-t_{1}\right)  \Gamma\left(  t_{3}-t_{1},\ldots,t_{8}-t_{1}\right)
\right\vert }\\
&  \text{ }\leq\frac{C}{n}\sum_{t_{2}^{\prime},\ldots,t_{8}^{\prime}=-n}%
^{n}\left\vert \Gamma\left(  0,t_{2}^{\prime}\right)  \Gamma\left(
t_{3}^{\prime},\ldots,t_{8}^{\prime}\right)  \right\vert \leq\frac{C}{n}%
\sum_{t=-n}^{n}\left\vert R_{t}\right\vert \sum_{t_{3}^{\prime},\ldots
,t_{8}^{\prime}=-n}^{n}\left\vert \Gamma\left(  0,t_{4}^{\prime}-t_{3}%
^{\prime},\ldots,t_{8}^{\prime}-t_{3}^{\prime}\right)  \right\vert \\
&  C\sum_{t=-\infty}^{\infty}\left\vert R_{t}\right\vert \sum_{t_{2}%
,\ldots,t_{6}=-\infty}^{\infty}\left\vert \Gamma\left(  0,t_{2},\ldots
,t_{6}\right)  \right\vert \leq C,
\end{align*}
by Assumption \ref{Kernel} and (\ref{Sumcum}).

The case $d_{B}=3$ corresponds to $\{ \operatorname*{Card}B_{1}%
,\operatorname*{Card}B_{2},\operatorname*{Card}B_{3}\}$ being $\{2,2,4\}$ or
$\{2,3,3\}$. We start with $\operatorname*{Card}B_{1}=2$,
$\operatorname*{Card}B_{2}=2$ and $\operatorname*{Card}B_{3}=4$. The
discussion concerns the number of fundamental pair broken by $B_{3}$. Note
that the situation where $B_{3}$ breaks only 3 or 1 fundamental pair is
impossible. The case where $B_{3}$ does not break any fundamental pairs
corresponds to partitions that are not indecomposable, so that the only
possible cases are those where $B_{3}$ breaks $4$ or $2$ fundamental pairs.

\begin{itemize}
\item $B_{3}$ breaks 4 fundamental pairs. Consider $B_{3}=\{1,2,3,4\}$,
$B_{2}=\{5,6\}$ and $B_{3}=\{7,8\}$. The corresponding contribution to $K_{2}$
is bounded by
\begin{align*}
&  \left\vert \frac{1}{n^{2}}\sum_{j_{1},j_{2}=1}^{n-1}k_{j_{1}}k_{j_{2}}%
\sum_{t_{1},t_{2}=1}^{n-j_{1}}\sum_{t_{3},t_{4}=1}^{n-j_{2}}\Gamma_{B}\left(
t_{1},\ldots,t_{8}\right)  \right\vert \\
&  =\left\vert \frac{1}{n^{2}}\sum_{j_{1},j_{2}=1}^{n-1}k_{j_{1}}k_{j_{2}}%
\sum_{t_{1},t_{2}=1}^{n-j_{1}}\sum_{t_{3},t_{4}=1}^{n-j_{2}}\Gamma\left(
0,t_{2}-t_{1},t_{3}-t_{1},t_{4}-t_{1}\right)  R_{t_{2}-t_{1}}R_{t_{4}-t_{3}%
}\right\vert \\
&  \text{ }\leq C\frac{p^{2}}{n}\sup_{j}|R_{j}|^{2}\sum_{t_{2},t_{3}%
,t_{4}=-\infty}^{\infty}\left\vert \Gamma\left(  0,t_{2},t_{3},t_{4}\right)
\right\vert \leq C\frac{p^{2}}{n}%
\end{align*}
by Assumption \ref{Kernel} and (\ref{Sumcum}).

\item $B_{3}$ breaks 2 fundamental pairs. Take $B_{3}=\{1,2,3,5\}$,
$B_{2}=\{4,6\}$ and $B_{1}=\{7,8\}$. The change of variables $t_{2}%
=t_{2}^{\prime}+t_{1}$, $t_{3}=t_{3}^{\prime}+t_{1}$ and $t_{4}=t_{4}^{\prime
}+t_{3}$ shows that contribution to $K_{2}$ is bounded by
\begin{align*}
&  \left\vert \frac{1}{n^{2}}\sum_{j_{1},j_{2}=1}^{n-1}k_{j_{1}}k_{j_{2}}%
\sum_{t_{1},t_{2}=1}^{n-j_{1}}\sum_{t_{3},t_{4}=1}^{n-j_{2}}\Gamma_{B}\left(
t_{1},\ldots,t_{8}\right)  \right\vert \\
&  \text{ }=\left\vert \frac{1}{n^{2}}\sum_{j_{1},j_{2}=1}^{n-1}k_{j_{1}%
}k_{j_{2}}\sum_{t_{1},t_{2}=1}^{n-j_{1}}\sum_{t_{3},t_{4}=1}^{n-j_{2}}%
\Gamma\left(  0,t_{2}-t_{1},t_{3}-t_{1},j_{1}\right)  R_{t_{4}-t_{2}-j_{1}%
}R_{t_{4}-t_{3}}\right\vert \\
&  \text{ }\leq\frac{C}{n}\sum_{j_{2}=1}^{n-1}K^{2}(j_{2}/p)\sum
_{t_{2}^{\prime},t_{3}^{\prime},j_{1}=-\infty}^{\infty}\left\vert
\Gamma\left(  0,t_{2}^{\prime},t_{3}^{\prime},j_{1}\right)  \right\vert
\sum_{t_{4}^{\prime}=-\infty}^{+\infty}\left\vert R_{t_{4}^{\prime}%
}\right\vert \times\sup_{j}|R_{j}|\leq C\frac{p}{n}.
\end{align*}
under Assumption \ref{Kernel} and (\ref{Sumcum}).
\end{itemize}

We now turn to the case $\operatorname*{Card}B_{3}=\operatorname*{Card}%
B_{2}=3$ and $\operatorname*{Card}B_{1}=2$. Observe that $B_{3}$ or $B_{2}$
must break 3 or 1 fundamental pair. The discussion now concerns the
fundamental pairs which are simultaneously broken by $B_{3}$ and $B_{2}$. Note
that $B_{3}$ and $B_{2}$ cannot break the same 3 fundamental pairs because if
it did, $B_{1}$ would be given by the remaining fundamental pair in which case
$B_{1}$ cannot communicate with $B_{2}$ or $B_{3}$, a fact that would
contradict the requirement that the partition $\{B_{1},B_{2},B_{3}\}$ is indecomposable.

\begin{itemize}
\item $B_{3}$ and $B_{2}$ break 3 fundamental pairs, 2 of which are the same.
Take $B_{3}=\{1,2,3\}$, $B_{2}=\{4,5,6\}$ and $B_{1}=\{7,8\}$. Using change of
variables $t_{2}=t_{1}+t_{2}^{\prime}$, $t_{3}=t_{1}+t_{3}^{\prime}$ and
$t_{4}=t_{3}+t_{4}^{\prime}$, we can see that under Assumption \ref{Kernel}
and (\ref{Sumcum}) the contribution to $K_{2}$ of this case is bounded by
\begin{align*}
&  \left\vert \frac{1}{n^{2}}\sum_{j_{1},j_{2}=1}^{n-1}k_{j_{1}}k_{j_{2}}%
\sum_{t_{1},t_{2}=1}^{n-j_{1}}\sum_{t_{3},t_{4}=1}^{n-j_{2}}\Gamma_{B}\left(
t_{1},\ldots,t_{8}\right)  \right\vert \\
&  \text{ }=\left\vert \frac{1}{n^{2}}\sum_{j_{1},j_{2}=1}^{n-1}k_{j_{1}%
}k_{j_{2}}\sum_{t_{1},t_{2}=1}^{n-j_{1}}\sum_{t_{3},t_{4}=1}^{n-j_{2}}%
\Gamma\left(  0,t_{2}-t_{1},t_{3}-t_{1}\right)  \Gamma\left(  0,t_{1}%
-t_{4}+j_{1},t_{2}-t_{4}+j_{1}\right)  R_{t_{4}-t_{3}}\right\vert \\
&  \text{ }\leq\frac{C}{n}\sum_{j_{1},j_{2}=1}^{n-1}K^{2}(j_{1}/p)K^{2}%
(j_{2}/p)\sup_{t_{2},t_{3}}\left\vert \Gamma(0,t_{2},t_{3})\right\vert
\sum_{t_{2}^{\prime},t_{3}^{\prime}=-\infty}^{\infty}\left\vert \Gamma\left(
0,t_{2}^{\prime},t_{3}^{\prime}\right)  \right\vert \sum_{t_{4}^{\prime
}=-\infty}^{+\infty}\left\vert R_{t_{4}^{\prime}}\right\vert \leq C\frac
{p^{2}}{n}%
\end{align*}
Note that the case where $B_{3}$ and $B_{2}$ break 3 fundamental pairs with
less than one in common is impossible.
\end{itemize}

The next case assumes that $B_{2}$ breaks only $1$ fundamental pair, which is
also necessarily broken by $B_{3}$ since $B_{2}$ must contain the remaining
unbroken pair.

\begin{itemize}
\item $B_{3}$ breaks 3 fundamental pairs and $B_{2}$ breaks only 1 pair. Take
$B_{3}=\{1,2,3\}$, $B_{2}=\{4,5,8\}$ and $B_{3}=\{6,7\}$ and consider a change
of variables $t_{2}=t_{1}+t_{2}^{\prime}$, $t_{3}=t_{1}+t_{3}^{\prime}$ and
$t_{4}=t_{1}+j_{1}-t_{4}^{\prime}$. Under Assumption \ref{Kernel} and
(\ref{Sumcum}), the contribution of this term to $K_{2}$ is bounded by\
\begin{align*}
&  \left\vert \frac{1}{n^{2}}\sum_{j_{1},j_{2}=1}^{n-1}k_{j_{1}}k_{j_{2}}%
\sum_{t_{1},t_{2}=1}^{n-j_{1}}\sum_{t_{3},t_{4}=1}^{n-j_{2}}\Gamma_{B}\left(
t_{1},\ldots,t_{8}\right)  \right\vert \\
&  \text{ }=\left\vert \frac{1}{n^{2}}\sum_{j_{1},j_{2}=1}^{n-1}k_{j_{1}%
}k_{j_{2}}\sum_{t_{1},t_{2}=1}^{n-j_{1}}\sum_{t_{3},t_{4}=1}^{n-j_{2}}%
\Gamma\left(  0,t_{2}-t_{1},t_{3}-t_{1}\right)  \Gamma\left(  t_{1}%
-t_{4}+j_{1},0,j_{2}\right)  R_{t_{3}-t_{2}+j_{2}-j_{1}}\right\vert \\
&  \leq\frac{C\sup_{j}|R_{j}|}{n}\sum_{j_{1}}^{n-1}K^{2}(j_{1}/p)\sum
_{t_{2}^{\prime},t_{3}^{\prime}=-\infty}^{\infty}\left\vert \Gamma
(0,t_{2}^{\prime},t_{3}^{\prime})\right\vert \sum_{t_{4}^{\prime}%
,j_{2}=-\infty}^{\infty}\left\vert \Gamma\left(  t_{4}^{\prime},0,j_{2}%
\right)  \right\vert \leq C\frac{p}{n}.
\end{align*}

\item $B_{3}$ and $B_{2}$ break only 1 pair. Note that $B_{3}$ and $B_{2}$
cannot break the same pair because $B_{1}$ must be the remaining pair and
cannot communicate, so that the partition is not indecomposable. Hence all the
partitions in this case are similar to $B_{3}=\{1,2,5\}$, $B_{2}=\{3,4,8\}$,
$B_{1}=\{6,7\}$. The change of variable $t_{2}=t_{1}+t_{2}^{\prime}$,
$t_{3}=-j_{2}+t_{2}+j_{1}+t_{3}^{\prime}$ and $t_{4}=t_{3}-t_{4}^{\prime}$
yields a contribution to $K_{2}$ bounded by
\begin{align*}
&  \left\vert \frac{1}{n^{2}}\sum_{j_{1},j_{2}=1}^{n-1}k_{j_{1}}k_{j_{2}}%
\sum_{t_{1},t_{2}=1}^{n-j_{1}}\sum_{t_{3},t_{4}=1}^{n-j_{2}}\Gamma_{B}\left(
t_{1},\ldots,t_{8}\right)  \right\vert \\
&  \text{ }=\left\vert \frac{1}{n^{2}}\sum_{j_{1},j_{2}=1}^{n-1}k_{j_{1}%
}k_{j_{2}}\sum_{t_{1},t_{2}=1}^{n-j_{1}}\sum_{t_{3},t_{4}=1}^{n-j_{2}}%
\Gamma\left(  0,t_{2}-t_{1},j_{1}\right)  \Gamma\left(  t_{3}-t_{4}%
,0,j_{2}\right)  R_{t_{3}-t_{2}+j_{2}-j_{1}}\right\vert \\
&  \text{ }\leq C\sum_{j_{1},t_{2}^{\prime}=-\infty}^{\infty}\left\vert
\Gamma(0,t_{2}^{\prime},j_{1})\right\vert \sum_{j_{2},t_{4}^{\prime}=-\infty
}^{\infty}\left\vert \Gamma(t_{4},0,j_{2})\right\vert \sum_{t_{3}^{\prime
}=-\infty}^{\infty}\left\vert R_{t_{3}^{\prime}}\right\vert \leq
C.\hfill\square
\end{align*}

\end{itemize}

\section*{Supplementary material additional references}

\textsc{Brillinger, D.R.} (2001). \textit{Time Series Analysis: Data Analysis
and Theory }. Holt, Rinehart \& Winston, New-York.

\textsc{Chow, Y.S.} and \textsc{H. Teicher} (1988). \textit{Probability
Theory. Independence, Interchangeability, Martingales }. Second Edition, Springer.

\textsc{Li, D.} and \textsc{R.J. Tomkins} (1996). Laws of the Iterated
Logarithm for Weighted Independent Random Variables. \textit{Statistics and
Probability Letters} \textbf{27}, 247--254.

\textsc{Priestley, M.B.} (1981). \textit{Spectral Analysis and Time Series}.
New York: John Wiley.

\end{document}